\documentclass[12pt]{amsart}
\usepackage{amsmath}
\usepackage{amssymb}
\usepackage{amstext}
\usepackage{amscd}
\usepackage[matrix,arrow,ps]{xy}

\usepackage[OT2,T1]{fontenc}
\DeclareSymbolFont{cyrletters}{OT2}{wncyr}{m}{n}
\DeclareMathSymbol{\Sha}{\mathalpha}{cyrletters}{"58}

\usepackage{mathtools} % for \DeclarePairedDelimiter
\usepackage[bookmarks=true,backref,colorlinks=true,citecolor=blue,urlcolor=blue,linkcolor=magenta]{hyperref}
\DeclarePairedDelimiter{\abs}{\lvert}{\rvert}
\DeclarePairedDelimiter{\norm}{\lVert}{\rVert}

\newcommand{\tens}{\otimes}
\newtheorem{theorem}{Theorem}[section]
\newtheorem{proposition}[theorem]{Proposition}

\newtheorem{lem}[theorem]{Lemma}
\newtheorem{cor}[theorem]{Corollary}
\newtheorem{conj}[theorem]{Conjecture}

\newtheorem{defini}[theorem]{Definition}

\newcommand{\sous}{\backslash}

\newcommand{\Hom}{\mbox{Hom}}

\newcommand{\End}{{\rm End}}

\newcommand{\GL}{{\rm GL}}

\newcommand{\Sh}{{\rm Sh}}

\newcommand{\dom}{ \preccurlyeq}

\newcommand{\wh}{\widehat}

\newcommand{\Gal}{{\rm Gal}}

\newcommand{\FF}{{\mathbb F}}
\newcommand{\CC}{{\mathbb C}}
\newcommand{\RR}{{\mathbb R}}
\newcommand{\ZZ}{{\mathbb Z}}
\newcommand{\QQ}{{\mathbb Q}}

\newcommand{\GG}{{\mathbb G}}
\newcommand{\SSS}{{\mathbb S}}
\newcommand{\AAA}{{\mathbb A}}

\newcommand{\lto}{\longrightarrow}

\def\FS{\mathfrak{S}}

\def\Fg{\mathfrak{g}}

\def\Fn{\mathfrak{n}}

\newcommand{\cH}{\mathcal{H}}

\newcommand{\cA}{{\mathcal A}}

\newcommand{\ol}{\overline}

\newcommand{\wt}{\widetilde}

\title[Heights and generalised Andr\'e-Pink-Zannier conjecture]{Height functions on Hecke orbits and the generalised Andr\'e-Pink-Zannier conjecture}
\author{Rodolphe Richard}
\address{UCL Department of Mathematics,
University College London, 
Gower Street, 
London, WC1E 6BT
}
\email{r.richard@ucl.ac.uk}
\author{Andrei Yafaev.}
\address{UCL Department of Mathematics,
University College London, 
Gower Street, 
London, WC1E 6BT
}
\email{yafaev@ucl.ac.uk}

\begin{document}

\maketitle

\begin{abstract}
We introduce and study the notion of a generalised Hecke orbit in a Shimura variety. We define a height function on such an orbit and study its properties.
We obtain lower bounds for the sizes of Galois orbits of points in a generalised Hecke orbit in terms of this height function, 
 assuming the ``weakly adelic Mumford-Tate hypothesis'' and
 prove the generalised Andr\'e-Pink-Zannier conjecture under this assumption, using Pila-Zannier strategy.

%This assumption is implied by the Serre's adelic form of Mumford-Tate conjecture, and, for Shimura
%varieties of abelian type, is known to be equivalent to the classical Mumford-Tate conjecture.
%by implementing the Pila-Zannier strategy.
\end{abstract}

\setcounter{tocdepth}{1}
\tableofcontents

%References

%Functoriality SerreL \S2.3

%BG Affine heights 1.5.6 1.5.7, (Proj) Weil heights \S2, 2.4.

%Calcul:
 
%For~$v=(v_1,\ldots,v_N)$ with~$v_i=n_i/d_i\in\QQ$ irreducible fractions~(if~$v_i=0$, then~$d_i=1$), let~$H^{\text{PW}}=\max_i\max\{d_i;\abs{n_i}\}$ be Pila's height and
%write~$H^{proj}=H_\infty\cdot H_f$ our height with~$H_\infty=\max\bigl\{1;\max_i\{\abs{n_i}/d_i\}\bigr\}\leq H^{\text{PW}}$ et~$H_f=gcd(d_i)\leq \prod d_i\leq (H^{\text{PW}})^d$.
%Then~$H^{PW}\leq H^{proj}\leq (H^{PW})^{N+1}$.

%Cadoret-Moonen:

%Notre weakly Adelic MT est une hypothese plus faible que ``Integral Mumford-Tate Conjecture `` [1.4],
%et que ``Adelic Mumford-Tate Conjecture `` (oprimale) Serre [2.5]

%Th. 5.3: (pour var abelinnes) MT classique (en une place) implique Adelic MT (et donc Weakly adelic)

\section{Introduction}

In this paper, we study the generalised Andr\'e-Pink-Zannier conjecture  for all Shimura varieties, 
whose statement is as follows.

\begin {conj}[Generalised André-Pink-Zannier] \label{APZ}
Let $S$ be a Shimura variety and $\Sigma$ a subset of 
a generalised Hecke orbit in $S$. Then the irreducible components of the Zariski closure of 
$\Sigma$ are weakly special subvarieties.
\end{conj}

We refer to~\cite{Del71,Del79} for notions and notations concerning Shimura data and Shimura varieties.
We refer to~\cite[Def.\,2.1]{UY_Char} for definitions and properties of weakly special subvarieties.
We refer to Def.~\ref{gho} or \S\ref{sec: main result} below for the notion of \emph{generalised} Hecke orbits.

%\emph{Generalised Hecke orbits} are defined in Def.~\ref{gho}. We believe that this notion of Hecke orbit is the 
%most general possible.
%We compare it to other notions of ``generalised Hecke orbits'' in \S\,\ref{Autre Hecke}. 
%The above version of the Andr\'e-Pink-Zannier conjecture is the most general possible.

\subsection{Main result}\label{sec: main result}
Let~$(G,X)$ be a Shimura datum, and~$K\leq G(\AAA_f)$ be a compact open subgroup, and~$S=Sh_K(G,X)=G(\QQ)\sous X\times G(\AAA_f)/K$ be the associated Shimura variety.
Let~$x_0\in X$ 
%, viewed as a morphism~$\SSS\to G_{\RR}$,
 and denote by~$M\leq G$ its Mumford-Tate group. Let~$s_0:=[x_0,1]\in S$.

The \emph{generalised Hecke orbit of~$x_0$ in~$X$} (see~\S\ref{sec: def orbits})
is the set~$\cH(x_0)$ of the~$\phi\circ x_0$, where~$\phi:M\to G$ ranges through the morphisms of~$\QQ$-algebraic groups
such that~$\phi\circ x_0\in X$. The \emph{generalised Hecke orbit of~$s_0$ in~$S$} is~$\cH(s_0):=G(\QQ)\sous \cH(x_0)\times G(\AAA_f)/K\subseteq S$.
%The generalised Hecke orbit~$\cH(x_0)$ of a point~$x_0\in X$ is the set of the~$[\phi\circ x_0,g]=G(\QQ)\sous G(\QQ)(\phi\circ x_0,g)K/K$ such that~$g\in G(\AAA_f)$ and~$\phi:M\to G$ is a morphism of~$\QQ$-algebraic groups
%such that~$\phi\circ x_0\in X$.
For a sufficiently large field~$E$ of finite type over~$\QQ$ we have the following (see~\S\ref{sec:def Galois repr}): ~$S$ and~$s_0$ are defined over~$E$ and there exists a
Galois representation~$\rho_{x_0}:Gal(\ol{E}/E)\to M(\AAA_f)\cap K$ such that
\[
\forall \sigma\in Gal(\ol{E}/E),~g\in G(\AAA_f),~ \sigma([x_0,g])=[x_0,\rho_{x_0}(\sigma)\cdot g].
\]

The main result of this paper is the following.
\begin{theorem}\label{thm intro}We consider the above situation.
We assume the weakly adelic Mumford-Tate hypothesis (see \S\ref{MT hypothesis}), which states that, with~$U:=\rho_{x_0}(Gal(\ol{E}/E))\subseteq M(\AAA_f)\cap K$:
\begin{equation}\label{waMT}
\exists C>0, \forall p, [K\cap M(\QQ_p):U\cap M(\QQ_p)]\leq C.
\end{equation}

Then, for any subset~$\Sigma\subseteq \cH(s_0)$, every irreducible component of~$\ol{\Sigma}^{Zar}$
is weakly special.
\end{theorem}

%Our main result, Theorem~\ref{main theorem}, proves the conjecture for generalised Hecke orbits of points satisfying  the ``weakly adelic Mumford-Tate hypothesis'' (see \ref{MT hypothesis}).
Our ``weakly adelic Mumford-Tate hypothesis'' is weaker than the adelic form of the Mumford-Tate conjecture~\cite[11.4?]{SerreConj} stated by Serre. 
Here are some instances in which above Theorem~\ref{thm intro} implies Conjecture~\ref{APZ} unconditionally. 

Combining Th.~\ref{thm intro}
with Theorem~\ref{lemme special est MT}, one recovers the following.
\begin{theorem}[\cite{EdYa} and~\cite{KY}]\label{thmEdYa}
 Conjecture~\ref{APZ} is true if~$\Sigma$ contains a special point.
\end{theorem}

Combining Th.~\ref{thm intro} with with~\cite[Th.~A (i)]{CM} we have the following,
which strictly contains a 2005 result of Pink~\cite[\S7]{Pink} (and~\cite[Th.~B]{CK}).
\begin{theorem}\label{th15} Conjecture~\ref{APZ} is true if~$S$ is of abelian type, and~$\Sigma$ contains a point~$s$
which satisfies the Mumford-Tate conjecture (at some~$\ell$, in the sense of~\cite{UY_General}).
\end{theorem}

The assumptions of Th.~\ref{th15} are satisfied in the case where~$S=\mathcal{A}_g$ and~$\Sigma$ contains a point~$[A]$,
where the abelian variety~$A$ satisfies the Mumford-Tate conjecture (at some prime~$\ell$).
Examples of such abelian varieties are: when~$\dim(A)\leq 3$; or when~$\dim(A)$ is odd and~$\End(A)\simeq \ZZ$.
More examples were given in~\cite{Pink98}, and many examples are mentioned in~\cite[\S2.4]{DL}. 

The assumptions of Th.~\ref{th15} are also satisfied for ``most'' points in~$S(\ol{\QQ})$ (with~$S$ of abelian type) in the following sense.
The subset consisting of the~$s\in S(\ol{\QQ})$ such that~$s$ does not satisfy the Mumford-Tate conjecture is thin in the sense of~\cite[\S 9.1]{MW}: this uses 
a combination of~\cite[\S1]{Serre133},~\cite[\S 9]{MW} and~\cite[Th.~A (i)]{CM} and~Th.~\ref{Baldi MT}.

%\begin{theorem}If~$S$ is of abelian type there exists a thin subset~$\Omega\subseteq S(\ol{\QQ})$ such that, if~$\Sigma\cap(S(\ol{\QQ}))\smallsetminus\Omega\neq\emptyset$,
%then the conjecture~\ref{APZ} is true if~$S$.
%\end{theorem}

For arbitrary Shimura varieties, the hypotheses of Th.~\ref{thm intro} are satisfied in the situation of~ Th.~\ref{Baldi MT}.
In a sense, our results apply unconditionally to ``most'' nonalgebraic points of a Shimura variety. The following are two special cases
of~ Th.~\ref{Baldi MT}.

\begin{theorem}\label{th16}
Conjecture~\ref{APZ} is true if~$\Sigma$ contains a $\ol{\QQ}$-Zariski generic point~$s$ of a special subvariety~$Z\subseteq S$, namely: for every proper 
subvariety~$V\subsetneq Z$ defined over~$\ol{\QQ}$, we have~$s\not\in V(\CC)$.
\end{theorem}
\begin{theorem}\label{th17}
Conjecture~\ref{APZ} is true if~$M^{ad}$ is~$\QQ$-simple and~$\Sigma$ contains a point~$s$ in~$S(\CC)\smallsetminus S(\ol{\QQ})$.
\end{theorem}

\subsection{History of Conjecture~\ref{APZ}}
Conjecture~\ref{APZ} is a special case\footnote{We refer to~\cite{Orr}, 
proof of Lemma 2.2, for the argument, which applies to our generalised setting.} of the Zilber-Pink conjecture, which has been 
and continues to be a subject of active research.

Conjecture \ref{APZ} was first formulated (in a special case) in 1989 by Y. Andr\'e in \cite[Ch. X, \S4.5]{Andre} (Problem 3).
U.~Zannier has considered questions of this type in the context of abelian schemes and tori in~\cite{Zannier}.
It was then stated in the introduction to the  second author's 2000 PhD thesis~\cite[bottom of p.\,12]{Y}\footnote{The statement there uses the terminology `totally geodesic subvarieties' instead of `weakly special', but Moonen had proved in \cite{MoMo} that the 
two notions are equivalent.}, following discussions with Bas Edixhoven.
Richard Pink, in his 2005 paper \cite{Pink}, has formulated and studied this question.

These authors consider the classical Hecke\footnote{Where André uses~$G(\QQ)$, Pink uses~$\operatorname{Aut}(G)(\QQ)$ instead of~$G(\QQ)$ in Def.\,\ref{cho}.} orbit as in Def.\,\ref{cho}.

Pink proves the Andr\'e-Pink-Zannier conjecture for ``Galois generic'' points of $\cA_g$.
These points are Hodge generic, by \cite[Pr. 6.2.1]{CK}. 
% This assumption certainly does not hold for all points in $\cA_g$, indeed it does not hold for all CM points.
Pink's method uses equidistribution of Hecke points (by Clozel-Oh-Ullmo: \cite{COU}; cf. also~\cite{EO}). 
This was generalised to Galois generic points in arbitrary Shimura varieties in 2016's~\cite{CK}. 
This was also contained in the 2009's first author's thesis under a weaker assumption~\cite[Ch. III \S7, p. 59, Cor.\,7.1]{R-These}.
 
%The result~Th.~\ref{thmEdYa}
In the case of generalised\footnote{They used a generalised notion of Hecke orbit, formulated using auxiliary linear representations; but using Prop.~\ref{prop:rho Hecke} and Th.~\ref{theo general union finie geometrique},
this leads to a statement equivalent to our Conj.~\ref{APZ}.} Hecke orbits \emph{of special points}, 
%Bas Edixhoven and the second author prove in \cite{EdYa} Conjecture~\ref{APZ} for 
%one dimensional components of the Zariski closure, and Bruno Klinger and the second author prove in~\cite{KY} Conjecture~\ref{APZ} for all components.
the articles~ \cite{EdYa,KY} use a method of  Edixhoven.
 %(introduced by Hindry in the context of abelian varieties) as an approach to the Andr\'{e}-Oort conjecture, which consists of exploiting the growth of Galois orbits  and intersections with images by Hecke correspondences.
 This method is inapplicable in more general cases, for instance 
 the case of the Hecke orbit of a \emph{Hodge generic} point.

A real breakthrough on this problem was the introduction of the Pila-Zannier strategy which uses o-minimality and 
functional transcendence.  
It has now become the most powerful approach to all problems of Zilber-Pink type.
This method was applied by M. Orr in \cite{Orr}, who considered the case of curves 
in $\cA_g$, the moduli space of principally polarised abelian varieties. His approach relies
on Masser-Wüstholz isogeny estimates. Therefore, it is limited to Shimura varieties of abelian type, and can't be applied to \emph{generalised} Hecke orbits.
For Shimura varieties of abelian type, Orr was able to prove the conjecture  for "$S$-adic Hecke orbits"\footnote{He considers Hecke correspondences whose level has only prime factors in~$S$. This corresponds to isogenies of abelian varieties whose degree has prime factors only from~$S$.}  for a finite set of primes~$S$, and for points which are Hodge generic (without the Galois generic assumption).

In the case of $S$-adic Hecke orbits, a stronger form of the conjecture, involving topological closure and equidistribution, was proved, in the abelian case, 
in \cite{RY} using ergodic theory approach relying on $p$-adic  Ratner's theorems. 

%This also extends to general Shimura varieties under a Shafarevich type assumption
%weaker than our Mumford-Tate assumption. This Shafarevich assumptions is satisfied in Shimura varieties of abelian type.

\subsection{Main technical results}\label{subsec:main technical}
After choosing bases of the Lie algebras~$\mathfrak{m}$ of~$M$ and~$\mathfrak{g}$ of~$G$, we associate to~$\phi\in\Hom(M,G)$ its "finite height"~$H_f(\phi)$, defined as the lowest common multiple of the denominators of the coefficients of the matrix of~$d\phi$.
More generally, for~$g\in G(\AAA_f)$, we define~$H_f(g^{-1}\cdot \phi\cdot g)$ as the smallest~$n\in\ZZ_{\geq1}$ such that the matrix of~$g^{-1}\cdot d\phi\cdot g$ has coefficients in~$\frac{1}{n}\cdot \widehat{\ZZ}$.
\subsubsection{}\label{tech1}
A first crucial result is the following. We choose the bases of~$\mathfrak{g}$ and~$\mathfrak{m}$ constructed in~\S\ref{sec: construction invariant height}. Then the function
\[
[\phi\circ x_0,g]\mapsto H_f(g^{-1}\cdot \phi\cdot g)
\]
is well defined on the generalised Hecke orbit, and $Gal(\ol{E}/E)$-invariant.

\subsubsection{}\label{tech2}
Our most important technical result is an estimate on the size of Galois orbits in a generalised Hecke orbit.

The following definition is used throughout this article.
\begin{defini}\label{notation dominates}
Let $A$ be a set and $f,g:A\to \RR_{\geq 0}$ two functions.

\begin{enumerate}
\item We say that $f$ \emph{polynomially dominates} $g$, and write $g \preccurlyeq f$, if
there exist $a,b, c \in \RR_{> 0}$ such that
$$
\forall x\in A, g(x) \leq c + a f(x)^b.
$$
\item We say that $f$ and $g$ are \emph{polynomially equivalent}, and write $f \approx g$, if
$f\dom g$ and $g\dom f$.
\end{enumerate}
\end{defini}

As functions on the generalised Hecke orbit~$\cH(s_0)$, we have the polynomial equivalence
\[
\# Gal(\ol{E}/E)\cdot [\phi\circ x_0,g] \approx H_f(g^{-1}\cdot \phi\cdot g).
\]
\subsubsection{}\label{tech3}
Another essential technical result, from~\S\ref{Section Siegel}, is the following. 
See the introduction in~\S\ref{Section Siegel} for the importance of this result in our approach
to Conjecture~\ref{APZ}.

Denote by~$\phi_0$ the inclusion monomorphism~$M\hookrightarrow G$.
Let~$W$ be the conjugacy class~$G\cdot \phi_0\subseteq \Hom(M,G)$, viewed as an algebraic variety over~$\QQ$.
The usual height of the matrix of~$d\phi$ defines an affine Weil height function~$H_W$ on~$W(\QQ)$ (cf.~\eqref{global iota height} and~\eqref{H sans iota}).
Let~$\mathfrak{S}\subseteq G(\RR)$ be a finite union of Siegel sets and~$\mathfrak{S}\cdot \phi_0$ be its image in~$W(\RR)$.

The main result~\ref{theorem type Orr} of~\S\ref{Section Siegel} is that, as functions of~$\phi\in W(\QQ)\cap \mathfrak{S}\cdot \phi_0$, we have
\[
H_f(\phi)\approx H_W(\phi).
\]
We note that every point of the geometric Hecke orbit can be written as~$[\phi\circ x_0,g]$ with~$g\in G(\AAA_f)$ and~$\phi\in W(\QQ)\cap \mathfrak{S}\cdot \phi_0$, provided~$\mathfrak{S}\subseteq G(\RR)$ is a fundamental set.

%For~$\phi\in W(\QQ)$ 
%For~$\phi\in W(\QQ)$ such that~$\phi\circ x_0\in X$, there exists~$g\in G(\RR)$ such that~$\phi=g\phi_0 g^{-1}$.
%A point of 
%There exists a finite union~$\mathfrak{S}$ of Siegel sets, such that any 
%
%For a given~$[\phi\circ x_0,1]$, we may choose~$\phi\in W(\QQ)$ such that~$\phi=g\phi_0 g^{-1}$ with~$g$ in a fundamental set, such as a finite union~$\mathfrak{S}$ of Siegel sets.
%
%There is 
%The usual height of the matrix of~$d\phi$ define an affine Weil function~$H_W$ on~$W(\QQ)$,
%and~$H_W$ is the product of~$H_f(\phi)$ above with an archimedean factor.
%
%As functions of~$\phi\in W(\QQ)\cap \mathfrak{S}/Z_G(M)(\RR)$, we have
%\[
%H_f(\phi)\approx H_W(\phi).
%\]

\subsection{Outline of the strategy}\label{Strategy outline}%To prove Theorem~\ref{thm intro}
%we implement the Pila-Zannier strategy.
 The proof of~Theorem~\ref{thm intro} is given in~\S\ref{Main proof}. The technical results of~\S\ref{subsec:main technical} play a crucial role in our approach.   Let us outline our approach.

We reduce the Conjecture~\ref{APZ} to the case where~$V:=\ol{\Sigma}=\ol{\{s_0;s_1;\ldots\}}$ is irreducible, $G$ is adjoint and~$V$ is Hodge generic in~$S$. We rely on functoriality properties~\S\ref{sec:functor Hecke} of geometric and generalised Hecke orbits.\footnote{This avoids one difficulty in the approach~\cite{Orr} of Orr.} Theorem~\ref{theo general union finie geometrique} allows us to use geometric and generalised Hecke orbits interchangeably. We also rely on the functoriality properties~\S\ref{sec:functor MT} of the assumption~\eqref{waMT}.

The final objective of the proof is to apply the geometric part of the André-Oort conjecture~\cite{UllmoAxLin} (or~\cite{RU}), and use
induction on the number of simple factors of~$M^{ad}$.
 For every $n$ large enough, we
construct a weakly special subvariety~$Z_n\subseteq V$ 
of non-zero dimension such that~$s_n\in Z_n$. Then~\cite{UllmoAxLin,RU} describes~$\ol{\bigcup Z_n}$, 
and we deduce Conjecture~\ref{APZ}.

In order to construct the non-zero dimensional~$Z_n$, we use the Pila-Zannier strategy. By~\eqref{WQ+ Hecke}, we identify~$\cH(s_0)$ with a subset of~$W(\QQ)$ where~$W=G\cdot \phi_0\simeq G/Z_G(M)$ is the algebraic variety of~\S\ref{tech3}.

Let~$\pi:G(\RR)\to X\to S$ be the uniformisation map, and~$\mathfrak{S}\subseteq G(\RR)$ is a finite union of Siegel sets such that~$S=\pi(\mathfrak{S})$. The goal is to apply the variant Th.~\ref{PilaWilkie} of Pila-Wilkie theorem, after constructing many rational points of small height in the set
\[
\wt{V}=\left(\stackrel{-1}{\pi}(V)\cap \mathfrak{S}\right)/Z_{G(\RR)}(M)\subseteq W(\RR)
\]
which is definable in the o-minimal structure~$\RR_{an,\exp}$.

%is obtained from the graph of the uniformisation map of a Shimura variety by (a Siegel set) of an 
%hermitian domain.

Let~$E$ be field of definition of~$V$. Then~$V$ contains the Galois orbits~$Gal(\ol{E}/E)\cdot s_n$. 

We introduce
\[
Q_n:=\{\phi\in \mathfrak{S}\cdot \phi_0\cap W(\QQ):[\phi\circ x_0:1]\in Gal(\ol{E}/E)\cdot s_n\}\subseteq \widetilde{V}.
\]
Denote by~$p$ the map~$G(\RR)\cdot \phi_0\to X$, where~$G(\RR)\cdot \phi_0\subseteq W(\RR)$. Each point~$s'\in Gal(\ol{E}/E)\cdot s_n$
lifts to a rational point~$\wt{s'}\in\wt{V}\cap W(\QQ)$.
We have surjections~$Q_n\to p(Q_n)\to Gal(\ol{E}/E)\cdot s_n$. Thus~$\#Q_n\geq \#Gal(\ol{E}/E)\cdot s_n$.

By~\S\ref{tech1},  the value of~$H_{f}$ is constant
as~$\phi$ ranges through~$Q_n$. By~\S\ref{tech3}, we also have~$H_f(\phi)\approx H_W(\phi)$.
By~\S\ref{tech2}, we have~$\#Q_n\geq \#Gal(\ol{E}/E)\cdot s_n\approx H_f(\widetilde{s_n})\approx H_W(\widetilde{s_n})$. 

%By~Pr.~\ref{Gal invariance}, the value of~$H_{[x_0,1]}$ is constant
%as~$s'$ ranges through~$Gal(\ol{E}/E)\cdot s_n$. It follows from~\eqref{Galois estimate} that there are~$\#Gal(\ol{E}/E)\cdot s_n \approx H_{[x_0,1]}({s_n})$ such points.\footnote{This is where the assumption~\eqref{waMT} is needed.} 
%By~\eqref{Type Orr}, we have~$H_{[x_0,1]}({s_n}) \approx H_\QQ(\widetilde{s_n})$. 
%For~$\phi\in Q_n$, we have~$H_{\AAA_f}(d\phi)=H_{[x_0,1]}([\phi\circ x_0,1])=H_{[x_0,1]}({s_n})\approx H_W(\wt{s_n})$.
%  \approx H_\QQ(\widetilde{s_n})$.

Thus~$\widetilde{V}$ contains~$\#Q_n\approx H_W(\wt{s_n})$ points of height~$\approx H_W(\wt{s_n})$.

By~Th.~\ref{PilaWilkie}, for sufficiently large~$n$, there exist~$\phi_n$ in~$Q_n$ such that~$p(\phi_n)\in Z^{alg}$, with~$Z=p(\widetilde{V})$. By Ax-Lindemann-Weierstrass theorem \cite{KUY},
it follows that~$s'_n=[\phi_n,1]\in Z_n\subseteq V$, for a non-zero dimensional weakly special subvariety~$Z_n$.
Using Galois action, we may assume~$s'_n=s_n$.

This concludes the proof of Th.~\ref{thm intro}.

\subsection{}
In Section~\ref{sect1}, we introduce and study generalised and geometric Hecke orbits.
In Section~\ref{section functoriality}, we recall properties of the representations~$\rho_{x_0}:Gal(\ol{E}/E)\to M(\AAA_f)$,
and we relate Galois orbits to orbits of~$U=\rho_{x_0}(Gal(\ol{E}/E))$.
In Section~\ref{section invariant height}, we make precise and prove~\ref{tech1}.
The Section~\ref{Section Siegel} deals with~\ref{tech2}.
In Section~\ref{section MT adelic}, we introduce and study the weakly adelic Mumford-Tate hypothesis,
and establish the estimates~\ref{tech3}. This relies on general estimates on adelic orbits, given in the appendices.
The content of Section~\ref{Main proof} was outlined in~\ref{Strategy outline}.

\subsection*{Acknowledgments}
The second named author has been contemplating the Andr\'e-Pink-Zannier problem 
and how to use the Mumford-Tate conjecture to approach it for twenty years.
He had a number of discussions on the subject with Bas Edixhoven, Richard Pink, Yves Andr\'e, Emmanuel Ullmo, Bruno Klingler, 
Martin Orr, Christopher Daw,...
Both authors would like to thank them all.

Both authors were supported by Leverhulme Trust Grant RPG-2019-180. The support of the Leverhulme Trust is gratefully acknowledged.

\section{Generalised and Geometric Hecke orbits.} \label{sect1}
In this section we define the notions of~\emph{generalised Hecke orbit} and of~\emph{geometric Hecke orbit}, and study their properties.
The heart of this section is Theorem~\ref{theo general union finie geometrique}, 
which in particular implies that generalised and geometric Hecke orbits can be used interchangeably in the statement of Conjecture \ref{APZ}.

These notions are naturally compatible with various operations on Shimura data. In particular we prove several statements which
will be important in reducing the conjecture~\ref{APZ} to the case where the Shimura variety is of adjoint type and~$\Sigma$ is Hodge 
generic in~$S$.

Finally \S\ref{Autre Hecke} compares our notions to different notions of generalised  Hecke orbits found in the literature.

\subsection{}\label{sec: def orbits}
Let $(G,X)$ be a Shimura datum.
We always assume, as in \cite{UY1}, that our Shimura datum is normalised so that $G$ is the generic Mumford-Tate group of $X$. 

Let $x_0$ be a point of $X$ and let $M \leq G$ be the Mumford-Tate group of $x_0$.
Recall that $x_0$ is a morphism $\SSS :=\operatorname{Res}_{\CC/\RR}(GL(1))\lto G_{\RR}$ and 
that~$M = x_0(\SSS)^{Zar, \QQ}$ is the smallest $\QQ$-algebraic subgroup of $G$ containing $x_0(\SSS)$.
In the rest of the paper we denote the identity monomorphism $M \hookrightarrow G$ by $\phi_0$.

%Consider $\Hom(M,G)$. This is a scheme over $\QQ$, possibly of infinite type.

In the following definition~$\Hom(M,G)$ denotes the set of algebraic group morphisms \emph{defined over~$\QQ$}.
\begin{defini}[Generalised Hecke orbit] \label{gho}
We define the \emph{Generalised Hecke orbit} $\cH(x_0)$ of $x_0$ \emph{in $X$} as
$$
\cH(x_0) := X\cap \{ \phi \circ  x_0 :  \phi \in \Hom(M,G)\}.
$$
\end{defini}
Let~$X_M=M(\RR)\cdot x_0\subset X$. Then~$(M,X_M)$ is a Shimura datum, and the~$\phi\in \Hom(M,G)$
such that~$\phi\circ x_0\in X$ are precisely those giving rise to a morphism of Shimura data~$(M,X_M)\to (G,X)$.
In particular~$\phi(X_M)\subseteq X$.
%we obtain a Shimura datum $(\phi(M), \phi(X_M))$ where $X_M = M(\RR) \cdot x_0$
%(where there is no ambiguity, we will write $M \subset G$ instead of $\phi_0(M)$).

Let $K$ be a compact open subgroup of $G(\AAA_f)$ and 
$\Sh_K(G,X)$ be the Shimura variety associated to this data.
There is a natural map 
$$
X \times G(\AAA_f) \lto \Sh_K(G,X)
$$
and we denote the image of a point $(x,g)$ by $[x,g]$.

%We define the generalised Hecke orbit in $\Sh_K(G,X)$ as follows:
\begin{defini} \label{ghoS}
We define the \emph{Generalised Hecke orbit} $\cH([x_0, g_0])$ of $[x_0, g_0]$ \emph{in $\Sh_K(G,X)$} by
$$
\cH([x_0,g_0]) := \{ [x,g] : x\in \cH(x_0), g\in G(\AAA_f). \}
$$
\end{defini}

Let~$W=G\cdot \phi_0$ be the conjugacy class of~$\phi_0$ which we view as an algebraic variety defined over~$\QQ$. 
Denoting by~$Z_G(M)$  the centraliser of $M$ in $G$,
we will identify~$G/Z_G(M)\simeq W$. 
The set ~$W(\ol{\QQ})$ is the $G(\ol{\QQ})$-conjugacy class of~$\phi_0$ in~$\Hom(M_{\ol{\QQ}},G_{\ol{\QQ}})$, and the points in~$W(\QQ)$ are the~$\QQ$-defined homomorphisms~$\phi\in\Hom(M,G)$
 which are conjugated to $\phi_0$ by elements of $G(\ol{\QQ})$.
%We denote~$\uHom(M,G)$ the scheme\footnote{If~$M$ is semisimple it is a finite disjoint union of irreducible affine algebraic varieties, which are each a~$G$-conjugacy class, and if~$M$ is reductive but not semisiple, it is an infinite countable disjoint union of such. In the latter case, it is not a scheme of finite type over~$\QQ$.} over~$\QQ$ of algebraic group homorphisms, and~$W\subset \uHom(M,G) $ the $G$-conjugacy class of~$\phi_0$. 
%Denoting~$Z_G(M)$ be the centraliser of $M$ in $G$
%we identify~$G/Z_G(M)\simeq W$. Since~$M$ is reductive,~$W$ is an affine algebraic variety over $\QQ$ (cf.~\cite{Ri-tuples}).
%The $G$-conjugacy class of $\phi_0$ in $\Hom(M,G)$ is identified with $W(\QQ)$.
%Note that $W(\QQ)$ corresponds to elements of $\uHom(M,G)(\QQ)$ which are conjugated to $\phi_0$ by elements of $G(\ol{\QQ})$.

In Def.~\ref{gho}, if we replace~$\Hom(M,G)$ by its subset~$W(\QQ)$, we obtain a more restrictive definition: that of a \emph{geometric} Hecke orbit. 
\begin{defini}\label{defi: geo orbit}
We define the \emph{geometric Hecke orbit} $\cH^g(x_0)$ of $x_0$ by
$$
\cH^g(x_0) =X\cap \{ \phi \circ x_0 : \phi \in W(\QQ)\} \subset \cH(x_0)
$$
and the \emph{geometric Hecke orbit} of $[x_0,g_0]$ by
$$
\cH^g([x_0,g_0]) %= \cH_S([x_0,g_0]) 
= \{ [x,g] : x\in \cH^g(x_0), g\in G(\AAA_f) \}.
$$
\end{defini}

The main result of this section is the following. 
\begin{theorem}\label{theo general union finie geometrique}
The generalised Hecke orbit $\cH(x_0)$ is a union of finitely many geometric Hecke orbits.
\end{theorem}
%We recall that~$M$ is the generic Mumford-Tate group of~$X_M$.% (as $M$ is the Mumford-Tate group
%of $x_0$). 

\begin{lem}  \label{biglemma}

Let~$\phi,\phi'\in \Hom(M,G)$ (defined over~$\QQ$) be such that~$\phi\circ x_0=\phi'\circ x_0$.

Then~$\phi=\phi'$.
\end{lem}
\begin{proof}One can check directly that
\[
H:=\{m\in M(\CC) : \phi(m)=\phi'(m)\}
\]
is a subgroup of~$M(\CC)$ (it is the ``equaliser'' of~$\phi$ and~$\phi'$). It is algebraic and defined over~$\QQ$ because~$\phi$ and~$\phi'$ are. It contains
 the image~$x_0(\CC)$ by hypothesis. But~$M$ is the Mumford-Tate group of~$x_0$: there is no proper $\QQ$-algebraic subgroup of~$M$ containing~$x_0(\CC)$. Therefore~$H=M$. Thus~$\phi=\phi'$.
\end{proof}
\begin{subequations}
The algebraic variety~$W$ is our central object in this article. We will use  the notations
\begin{equation}\label{eq W R plus}
\begin{aligned}
W(\RR)^+&=G(\RR)/Z_G(M)(\RR)\\
&=\{\phi\in W(\RR):\phi\circ x_0\in X\}
\end{aligned}
\end{equation}
and
\begin{equation}\label{eq W Q plus}
\begin{aligned}
W(\QQ)^+&=W(\RR)^+\cap W(\QQ)\\
&=\{\phi\in W(\QQ):\phi\circ x_0\in \cH^g(x_0)\}.
\end{aligned}
\end{equation}
\end{subequations}
The subset~$W(\RR)^+\subset W(\RR)$ is a union of some connected components of~$W(\RR)$.
% that we will call the ``neutral components of~$W(\RR)$''.
%Other connected components will be irrelevant to us.
 With these notations, Lemma \ref{biglemma} implies that we have a bijection
\begin{equation}\label{WQ+ Hecke}
\begin{aligned}
W(\QQ)^+&\stackrel{\sim}{\to}\cH^g(x_0)\\
\phi&\mapsto \phi\circ x_0.
\end{aligned}
\end{equation}
%Elements of~$\cH^g(x_0)$ can be identified with some rational points of~$W$, but will not be rational points of~$X$ for a model of~$X$ over~$\QQ$. This is consequential in our use of Pila-Wilkie theorem.generalised Hecke orbits.

\subsection{Functoriality of Generalised and Geometric Hecke orbits}\label{sec:functor Hecke}
%We will now examine some useful functorial properties of the generalised and geometric Hecke orbits.
\subsubsection{Restriction to special subvarieties}
The following is a set theoretic tautology.

\begin{proposition} \label{prop passage to Hodge generic}
Let~$(G',X')$ be a Shimura datum with~$M\leq G'\leq G$ and~$X_M\subset X'\subset X$,
and define~$K'=G'(\AAA_f)\cap K$.
\begin{enumerate}
\item Let~$\cH'(x_0)$ be the generalised Hecke orbit of~$x_0$ viewed as a point of~$X'$. 

Then
\[
\cH'(x_0)=
\cH(x_0)\cap X'.
\]
\item
Let~$\cH'([x_0,1])$ be the generalised Hecke orbit of~$[x_0,1]$ viewed as a point of~$Sh_{K'}(G',X')$,
and~$S'$ the image of~$$f:=Sh(\iota):Sh_{K'}(G',X')\to Sh_{K}(G,X)$$ where~$\iota:G'\to G$ is the inclusion.
Then
\[
\cH([x_0,1])\cap S' = f(\cH'([x_0,1]))\text{ and }\cH'([x_0,1]) =\stackrel{-1}{f}(\cH([x_0,1])).
\]
\end{enumerate}
\end{proposition}

The following corollary can be deduced by combining Lemma \ref{biglemma} with Th.~\ref{theo general union finie geometrique}
(it can also be deduced from~\cite{Ri-Conj}).  
\begin{cor}
We keep previous notations. Then
\[
\cH^g(x_0)\cap X'
\]
is a finite union of geometric Hecke orbits in~$X'$. 
\end{cor}
Accordingly,~$\stackrel{-1}{f}(\cH^g([x_0,1]))$ is the image of finitely many geometric Hecke orbits in~$Sh_{K'}(G',X')$.

\subsubsection{Compatibility to products}\label{section compatibility products}
A useful property of geometric Hecke orbits is the compatibility with respect to products of Shimura data.
\begin{lem}\label{Lem functor produit Hecke}
Let~$(G,X)$ be an adjoint Shimura datum, and factor~$G=G_1\times\ldots \times G_f$ as a product of its~$\QQ$-defined
simple normal subgroups, and assume~$K=K_1\times \ldots \times K_f$ for compact open subgroups~$K_i\leq G_i(\AAA_f)$.
We denote~$X=X_1\times \ldots \times X_f$ the corresponding factorisation, and choose~$x_0= (x_1,\ldots,x_f)\in X_1\times \ldots \times X_f$.
We denote~$\cH^g(x_i)$ the geometric Hecke orbit of~$x_i$ with respect to the Shimura datum~$(G_i,X_i)$.

With respect to the corresponding factorisation of Shimura varieties
\[
Sh_K(G,X)=Sh_{K_1}(G_1,X_1)\times \ldots \times Sh_{K_f}(G_f,X_f)
\]
we have
\[
\cH^g(x_0)=\cH^g(x_1)\times \ldots \times \cH^g(x_f).
\]
\end{lem}
It follows from Lemma~\ref{Lem functor produit Hecke} that, at the level of Shimura varieties,
\[
\cH^g([x_0,1])=\cH^g([x_1,1])\times \ldots \times \cH^g([x_f,1] ).
\]

\begin{proof}Since~$G$ is adjoint, we have a factorisation %product of Shimura data are simpler, and we have notably
\[
X=X_1\times \ldots \times X_f.
\]
Let~$M$ be the Mumford-Tate group of~$x_0$ and let~$\phi_0=(\phi_1,\ldots,\phi_f):M\to G=G_1\times \ldots\times G_f$ be the inclusion.
As the conjugacy class in a product is the product of conjugacy classes, we have
\[
G\cdot\phi_0=G_1\cdot \phi_1\times \ldots \times G_f\cdot \phi_f. 
\]
The Mumford-Tate group of~$x_i$ is~$M_i:=\phi_i(M)$.
 Because~$x_0(\SSS)$ is Zariski dense over~$\QQ$ in~$M$ so is~$x_i(\SSS)$ in~$M_i$. Let~$\phi_i':M_i\to G_i$ be the identity map. We can identify~$G_i\cdot \phi_i\simeq G_i\cdot \phi'_i$, and have
\[
\cH^g(x_i)=\{g\cdot \phi'_i\circ x_i : g\in G_i\}\cap X_i=\{g\cdot \phi_i\circ x_i : g\in G_i\}\cap X_i.
\]
 The rest follows from the definition of geometric Hecke orbits.
\end{proof}

\subsubsection{Passing to the adjoint Shimura datum}\label{subsection passing to adjoint}
 The following property is used to
reduce the proof of Conjecture~\ref{APZ} and Theorem~\ref{thm intro}
to the case where~$G$ is adjoint.
\begin{lem}\label{lem:pass to adjoint}
 Let~$ad:(G,X)\to (G^{ad},X^{ad})$ be the map of Shimura data\footnote{Where~$(G^{ad},X^{ad})$ is as in~\cite[Prop. 2.2]{EdYa}.} induced
by the natural morphism~$ad:G\to G^{ad}$ and choose a compact open subgroup~$K^{ad}\leq G^{ad}(\AAA_f)$
containing~$ad(K)$.
Let~$ad:x\mapsto x^{ad}:=ad\circ x$ be the map~$X\to X^{ad}$ and 
\[
Sh(ad): Sh_K(G,X)\to Sh_{K^{ad}}(G^{ad},X^{ad})
\]
the corresponding morphism of Shimura varieties. 

Let $x_0\in X$.
 Recall that~$\cH^g(x_0)$ and~$\cH^g(x_0^{ad})$ denote the geometric Hecke orbit of~$x_0$ and~$x_0^{ad}$ with respect to~$G$
and~$G^{ad}$. 

We have
\begin{equation}\label{geom et adjoint}
ad(\cH^g(x_0)) \subseteq ad(X)\cap \cH^g(x_0^{ad}).
\end{equation}
%\item With the associated~$\cH^g([x_0,1])$ and~$\cH^g([x_0^{ad},1])$, and denoting~$S'$ the image of~$Sh(ad)$, we have
%\[
%Sh(ad)(\cH^g([x_0,1]))=S'\cap\cH^g([x_0^{ad},1])
%\]
\end{lem}

Lemma~\ref{lem:pass to adjoint} implies the inclusion
$$
ad(\cH^g(x_0)) \times G(\AAA_f)) \subseteq \cH^g(x_0^{ad}) \times G^{ad}(\AAA_f).
$$
Passing to the quotient, we obtain the following.

\begin{cor}We have
\(
Sh(ad)(\cH^g([x_0,1]))\subseteq \cH^g([x_0^{ad},1]).
\)
\end{cor}
We now prove Lemma~\ref{lem:pass to adjoint}.
\begin{proof}

%We prove that  $ad(\cH^g(x_0)) \subset ad(X)\cap \cH^g(x_0^{ad})$.
Choose $x\in \cH^g(x_0)$.
Clearly $x':=ad(x) \in ad(X) \subset X^{ad}$.

The Mumford-Tate group of~$x'_0:=ad(x_0)$ is $M':=ad(M)$. We denote by~$\phi_0' \colon M' \to G^{ad}$
the natural injection. We can write~$x=\phi\circ x_0$ with $\phi = g \phi_0 g^{-1}$ and $g \in G(\ol{\QQ})$.
%Let $M'$ be $ad(M)$ and $\phi_0' \colon M' \lto G^{ad}$ induced by $\phi_0$ and 
Then $\phi' := ad(g) \phi'_0 ad(g)^{-1}$ is defined over $\QQ$ because the map~$G \cdot \phi_0\to G^{ad} \cdot \phi'_0$ between conjugacy classes is a morphism of varieties
defined over~$\QQ$.
One computes~$x'=ad(gx_0g^{-1})=ad(g)ad(x_0)ad(g)^{-1}=\phi'\circ x_0'$, where~$x_0'\in X^{ad}$, and~$\phi$ is defined over~$\QQ$ and conjugated to~$\phi_0'$
over~$\ol{\QQ}$; that is:~$x'\in \cH^g(x_0')$.%QED
%
%Then we have:
%$$
%ad\circ \phi = \phi' \circ ad |_{M}
%$$
%Then
%$ad\circ \phi \circ x_0 = \phi'  \circ ad \circ x_0$, hence $\phi' \circ x'_0 \in \cH^g(x_0^{ad})$.
\end{proof}

\subsubsection*{Remarks} In~\eqref{geom et adjoint}, the reverse inclusion is also true, but it is not used in this paper, and its proof is left to the interested reader.
The inclusion~\eqref{geom et adjoint} and the proof we have given also applies to general morphisms of Shimura data~$(G,X)\to (G',X')$ instead of just~$(G,X)\to (G^{ad},X^{ad})$.

\subsection{Rational conjugacy of linear representations}

The following notable fact will be used at several places in this article. We believe this property is also of independent interest.
\begin{theorem}[{\cite[\S12.3, third parag.]{BT}}]\label{rational conjugation representations}
For any algebraic group~$M$ over~$\QQ$, any two representations~$\phi,\phi':M\to GL(n)$ which are defined over~$\QQ$ and
conjugated under~$GL(n,\ol{\QQ})$ are actually conjugated under~$GL(n,\QQ)$.
\end{theorem}
 It follows from the theory of linear representations for which references are for example\,\cite[Ch.~XI]{H-LAG} for~$\ol{\QQ}$, and~\cite[\S12]{BT} over~$\QQ$. We will only need the case where~$M$ is connected and reductive, and this case can be found for instance in~\cite[\S12.3, third paragraph]{BT}. They give a Galois cohomology argument, and the  same Galois cohomology argument works in general with a reference to \cite[1.7 Ex.1 p.16, 1969]{Kn} instead. For reductive groups, it is also possible to reduce the result to Skolem-Noether theorem. For tori, it can be reduced to the fact that any matrix is rationally conjugated to its canonical companion form.

\subsection{Proof of the finiteness Theorem~\ref{theo general union finie geometrique}}
The strategy will combine an argument for semisimple groups and another for algebraic tori.

\begin{proposition}\label{prop finitude M semisimple}
Let~$M$ be a semisimple algebraic group over~$\QQ$ (resp.~$\ol{\QQ}$).
\begin{enumerate}
\item
For all~$d\in\ZZ_{\geq0}$, the set of linear representations defined over~$\QQ$ 
(resp.~$\ol{\QQ}$)
\[
\Hom(M,GL(d))
\]
is a finite union of conjugacy classes under~$GL(d,\QQ)$ (resp. under~$GL(d,\ol{\QQ})$.)
\item
Let~$G$ be a reductive linear algebraic group over~$\QQ$ (resp.~$\ol{\QQ}$). Then the set of 
homomorphisms defined over~$\QQ$ (resp.~$\ol{\QQ}$)
\[
\Hom(M,G)
\]
is contained in (resp. is equal to) a finite union of $G(\ol{\QQ})$-conjugacy classes.
\end{enumerate}
\end{proposition}
For simplicity, we will only give an argument which assumes~$M$ is Zariski connected, which is the case considered in the proof of Th.~\ref{theo general union finie geometrique}.
\begin{proof}
We prove the first assertion.
By virtue of Theorem~\ref{rational conjugation representations}, it is enough to treat the case where everything is defined over~$\ol{\QQ}$. 

Because~$M$ is connected it is enough to prove that there are finitely many conjugacy classes of Lie algebra representations~$\mathfrak{m}\to \mathfrak{gl}(d)$.
Equivalently there are finitely many isomorphisms classes of linear representations of~$\mathfrak{m}$ of dimension~$d$. For this\footnote{These representations are sums of irreducible representations. By the Theorem of the Highest Weight~\cite[\S7.2, Th.\,7.15]{Hall}, the irreducible representations is parametrised by dominant weights. The dimension of irreducible representations are given by Weyl dimension formula~\cite[\S7.6.3,  Th.\,7.43]{Hall}, from which lower bounds for dimensions are easily derived: there are finitely many isomorphism classes of irreducible representations of bounded dimension.}, we refer to~\cite[\S7]{Hall}.
% in term of their highest weight, and a lower bound on the dimension is found in~\cite{}.

For the second assertion we treat the case where everything is defined over~$\ol{\QQ}$, which implies the case where everything is defined over~$\QQ$. 
It is deduced from the first part by using~\cite[Th. 3.1]{Ri-Conj}.
\end{proof}

We prove Theorem~\ref{theo general union finie geometrique} combining~\cite[Lem. 2.6]{UY1} with Prop.~\ref{Bounding conjugacy classes}.
\begin{proof} We identify~$G$ with its image by a faithful representation~$G\to GL(d)$, and we let~$\Sigma=\{\phi\in \Hom(M,G):\phi\circ x_0\in X\}$.

Thanks\footnote{This is where the property~$\phi\circ x_0\in X$ is used. This also needs that the image of~$x_0$ is $\QQ$-Zariski dense in~$M$.} to~\cite[Lem. 2.6]{UY1}, we may use Proposition~\ref{Bounding conjugacy classes}, and deduce that~$\Sigma=\{\phi\in \Hom(M,G):\phi\circ x_0\in X\}$ is contained in finitely many $GL(d)$-conjugacy classes. Using~\cite{Ri-Conj}, we conclude that~$\Sigma$ is contained in finitely many~$G(\ol{\QQ})$-conjugacy classes, thus 
proving Theorem~\ref{theo general union finie geometrique}.
\end{proof}

\begin{proposition}[Bounding conjugacy classes]\label{Bounding conjugacy classes}
Let~$M$ be a connected reductive $\ol{\QQ}$-group,~$M^{der}$ its derived subgroup and~$T=Z_M(M)^0$ its connected centre.

A subset~$\Sigma\subseteq  \Hom(M,GL(d))$ is contained in finitely many $GL(d)$-conjugacy classes if and only if: there is a finite set of characters~$F\subset X(T)$ such that for every~$\rho\in\Sigma$, all the weights of the representation~$\rho\restriction_T:T\to GL(d)$ belong to~$F$.
\end{proposition}
\begin{proof}
Because the set of characters is invariant under conjugation, the condition is necessary. We prove that this condition is also sufficient.

We know that two representations of a torus~$T$ are conjugated if and only if they have the same weights, with same multiplicities.
As the weights belongs to~$F$, and the dimension~$d$ is fixed, there are only finitely many possibilities for these weights and multiplicities. Hence~$\{\rho\restriction_T:\rho\in\Sigma\}$ is contained in at most finitely many conjugacy classes~$GL(d)\cdot \rho_1\restriction_T,\ldots,GL(d)\cdot \rho_c\restriction_T$. Without loss of generality we may assume that there is only one conjugacy class, say~$GL(d)\cdot \rho_{1}\restriction_{T}$.

We want to prove that
\begin{equation}\label{finitely many text}
\text{ there are finitely many~$\rho\in\Sigma$, up to $GL(d)$-conjugation. }
\end{equation} Possibly after conjugating, we may assume~$\rho\restriction_T=\rho_1\restriction_T$. Because~$M$ is connected, one has~$M=M^{der}\cdot T$. Thus
\begin{equation}\label{determined text}
\text{ $\rho$ is determined by~$\rho\restriction_{M^{der}}$ and~$\rho\restriction_{T}$.}
\end{equation}
As~$M^{der}$ and~$T$ commute with each other,~$\rho\restriction_{M^{der}}:M^{der}\to GL(d)$ factors through~$G':=Z_{GL(d)}(\rho_1(T))$. As~$T$ is reductive, so is~$G'$.

By Proposition~\ref{prop finitude M semisimple}, these~$\rho\restriction_{M^{der}}$ belong to finitely many conjugacy classes~$G'\cdot \rho_{1,1}\restriction_{M^{der}} ,\ldots,G'\cdot \rho_{1,e}\restriction_{M^{der}}$. 
Possibly after conjugating~$\rho$ by some~$g\in G'$, which does not change~$\rho\restriction_T$, we have 
\[
\rho\restriction_T=\rho_1\restriction_T\text{ and }\rho\restriction_{M^{der}}
\in
\{\rho_{1,1}\restriction_{M^{der}} ;\ldots;\rho_{1,e}\restriction_{M^{der}}\}.
\]
In light of~\eqref{determined text}, this proves~\eqref{finitely many text} and the conclusion follows.
\end{proof}

\subsection{Relation to other notions of Hecke orbits} \label{Autre Hecke}

The following is not used in the rest of this article, however it clarifies the relation between different notions of Hecke orbits
and we believe it to be of independent interest.
We compare our generalised and geometric Hecke orbits to the classical Hecke orbits and another notion of ``generalised Hecke'' orbit found in the literature.

\subsubsection{Relation to the classical definition of Hecke orbit}
Let us recall the notion of the classical Hecke orbit.

\begin{defini}[classical Hecke orbit] \label{cho}
Define the \emph{classical Hecke orbit} of $x_0$ as follows:
$$
\cH^c(x_0) = \{ \phi \circ x_0 \in X : \phi \in G(\QQ)/Z_G(M)(\QQ) \} \subset \cH(x_0)
$$
and  the \emph{classical Hecke orbit} of $[x_0,1]$ as
\[
\cH^c(x_0) = \{ [x,g] : x\in\cH^c(x_0), g\in G(\AAA_f)\}.
\]
\end{defini}

We have a chain of inclusions:
\begin{align}\label{comparison classical hecke types}
\cH^c(x_0) \subset \cH^g(x_0) \subset \cH(x_0)\\
\cH^c(s_0) \subset \cH^g(s_0) \subset \cH(s_0).
\end{align}
In general,~$\cH^g(x_0)$ is not a finite union of classical Hecke orbits,
even when~$G$ is of adjoint type.

\subsubsection*{Hecke correspondences}

Recall that the classical Hecke orbit can be described using Hecke correspondences.
For~$g\in G(\QQ)$, the points~$s_0=[x_0,1]$ and~$s_g=[g\cdot x_0,1]$ have a common inverse image
by the left, resp. right, finite map in
\[
Sh_K(G,X) \xleftarrow{Sh(Ad_1)} Sh_{K\cap gKg^{-1}}(G,X) \xrightarrow{Sh(Ad_g)}  Sh_K(G,X) 
\]
where~$Sh(Ad_g)$ the right map is the Shimura morphism associated to the map of Shimura data~$AD_g:(G,X)\to (G,X)$
induced by the conjugation~$AD_g:G\to G$ and~$Sh(Ad_1)$ is induced by the identity map~$AD_1:G\to G$.

Likewise generalised Hecke orbits can be interpreted using finite correspondences between Shimura varieties. For a point~$\phi\circ x_0\in \cH(x_0)$,
the point~$s_0$ and~$s_\phi=[\phi\circ x_0,1]$ have a common inverse image in
\[
Sh_K(G,X) \xleftarrow{Sh(\phi_0)} Sh_{K\cap \stackrel{\raisebox{0pt}[0pt][0pt]{$\scriptscriptstyle{-1}$}}{\phi}(K)}(M,X_M) \xrightarrow{Sh(\phi)}  Sh_K(G,X).
\]
This time the correspondence is induced by a correspondence from the image of~$Sh(\phi_0)$ to that of~$Sh(\phi)$. These
are also the smallest special subvarieties containing~$s_0$, resp.~$s_\phi$.

\subsubsection{Relation to the usual definition of the generalised Hecke orbit}
%This section in not used in the rest of the paper.
We compare our notion of generalised Hecke to the ``generalised Hecke orbits'' used in~\cite{KY} and~\cite{EdYa,Pink,Orr,UY_General}. The latter is defined in terms of linear representations.

For any faithful representation~$\rho:G\to GL(N)$ over~$\QQ$, let the ``$\rho$-Hecke orbit'' be
\[
\cH^\rho(x_0) := \{ \phi \circ  x_0  \in X :  \phi \in \Hom(M,G)(\QQ), \rho\circ \phi \in GL(N,\QQ)\cdot \rho\circ \phi_0\}.
\]
By Theorem~\ref{rational conjugation representations}, we also have
\[
\cH^\rho(x_0) = \{ \phi \circ  x_0  \in X :  \phi \in \Hom(M,G)(\QQ), \rho\circ \phi \in GL(N,\overline{\QQ})\cdot \rho\circ \phi_0\}.
\]

\begin{proposition}\label{prop:rho Hecke}
The  $\rho$-Hecke orbit~$\cH^\rho(x_0)$ is contained in the generalised Hecke orbit~$\cH(x_0)$.

The  $\rho$-Hecke orbit~$\cH^\rho(x_0)$ is a finite union of geometric Hecke orbits~$\cH^g(x_0)\cup\ldots\cup\cH^\rho(x_k)$.
\end{proposition}
The first statement is clear from the definition of~$\cH^\rho(x_0)$. The second statement follows from the second definition of~$\cH^\rho(x_0)$ and~\cite{Ri-Conj}.

The number of geometric Hecke orbits is bounded independently from~$\rho$ thanks to Theorem~\ref{theo general union finie geometrique}.
It is unclear whether we can achieve~$\cH^\rho(x_0)=\cH(x_0)$ for a sufficiently general representation~$\rho$.

\section{Galois functoriality on the generalised Hecke orbit}
\label{section functoriality}

%\subsubsection{Galois representation}
In~\S\ref{sec:def Galois repr} and~\S\ref{subsection functoriality} we state known definitions and properties for the convenience of the reader. Details can be found, for instance, in~\cite{UY_General}. 
In~\S\ref{Galois orbits vs Adelic orbits} we relate cardinality of Galois orbits and cardinality of orbits in adelic groups. This is essential to our approach to the estimates of ~\S\ref{tech2}  through adelic methods.

\subsection{Galois representations}\label{sec:def Galois repr}
 Our statements will use the following terminology.
\begin{defini}[Galois representations]\label{terminologie galois repr}
Let~$(M,X_M)$ be a Shimura datum, let~$x_0$ be a point in~$X_M$, and let~$E\leq \CC$ be a subfield.% such that~$[x_0,1]\in Sh_K(M,X_M)$ for some compact open subgroup~$K\leq M(\AAA_f)$. 

We say that a continuous homomorphism
\begin{equation}\label{representation galoisienne}
\rho=\rho_{x_0}:\Gal(\ol{E}/E)\to M(\AAA_f)
\end{equation}
is \emph{a Galois representation (defined over~$E$) for~$x_0$ (in~$X_M$)} if: for any compact open subgroup~$K'\leq M(\AAA_f)$, denoting~$[x_0,1]'$
the image of~$(x_0,1)$ in~$Sh_{K'}(M,X_M)$, we have
\begin{equation}\label{defining repr galois}
\forall\sigma\in \Gal(\ol{E}/E), \sigma([x_0,1]')=[x_0,\rho_{x_0}(\sigma)]'.
\end{equation}
\end{defini}
In the important case of moduli spaces of abelian varieties, a representation~$\rho_{x_0}$ can be directly constructed from the linear Galois action on the Tate module (see~\cite{UY_General,CM}).

 Here we only need the existence of a~$\rho_{x_0}$.
 %For most purposes such a~$\rho_{x_0}$ is essentially unique.

\begin{proposition}[Existence of Galois representations]
\label{existence galois repr} Let~$[x_0,1]\in Sh_{K_M}(M,X_M)(E')$ be a point defined over a field~$E'\leq \CC$ in a Shimura variety.

Then there exist a finite extension~$E/E'$ and a Galois representation defined over~$E$ for~$x_0$ in~$X_M$.
\end{proposition}
%\footnote{This is essentially derived from one main property of the canonical model of~$Sh_{K_M}(M,X_M)$: the rationality of Hecke correspondances.}
The main ingredient in this proposition is the following, 
which is part of the definition of canonical models:
for any~$[x_0,m_0]$, any~$m\in M(\AAA_f)$ and~$\sigma \in Aut(\CC/E(M,X_M))$,
\begin{equation}\label{Hecke covariance}
\text{ if }\sigma([x_0,m_0])=[x',m']\text{ then }\sigma([x_0,m_0\cdot m])=[x',m'\cdot m].
\end{equation}

The continuity of~$\rho_{x_0}$ is used in the following lemma.
\begin{lem} Let~$K$ be an open subgroup of~$M(\AAA_f)$. Then, after possibly replacing~$E$ by a finite extension, we have
\begin{equation}\label{eq unique rho}
\rho_{x_0}(\Gal(\ol{E}/E))\leq K.
\end{equation}
\end{lem}
\begin{proof}
Such an extension corresponds to the open subgroup~$\stackrel{-1}{\rho_{x_0}}(K)\leq \Gal(\ol{E}/E)$.
\end{proof}
\subsubsection*{Comments}If~$K$ is sufficiently small so that~$K\cap Z_G(M_0)(\QQ)=\{1\}$, for instance if~$K$ is \emph{neat} then (see~\cite[\S4.1.4]{KY}) for any field~$E\leq \CC$, there is at most one Galois representation~$\rho_{x_0}$ satisfying~\eqref{eq unique rho}.

\subsection{Functoriality of the Galois representation}\label{subsection functoriality}
In the next statement we denote by~$E(G,X)$ the \emph{reflex field} of a Shimura datum~$(G,X)$. 
It is a number field over which $Sh(G,X)$ (and hence all the 
~$Sh_K(G,X)$)  admits a \emph{canonical model}.

\begin{proposition}[Functoriality]\label{Prop functoriality}
Let
\(
\phi:(M,X_M)\to (G,X)
\)
be a morphism of Shimura data, and~$x_0$ a point in~$X_M$. 

If~$\rho_{x_0}$ is a Galois representation 
 defined over a field $E$ for $x_0$, then
\[\phi\circ \rho_{x_0}|_{Gal(\ol{E}/E\cdot E(G,X))}\] 
is a Galois representation defined over~$E\cdot E(G,X)$ for~$\phi(x_0)$ in~$X$.
\end{proposition}
This follows from the definition and the identity
\[
\sigma([\phi\circ x,\phi(g)])=[\phi\circ x',\phi(g')]\text{ for }[x',g']=\sigma([x,g])
\]
which holds when~$\sigma\in \operatorname{Aut}(\CC/E(M,X_M)E(G,X))$. Equivalently the Shimura morphisms induced by~$\phi$ are defined over~$E(M,X_M)E(G,X)$. (See \cite[1.14, 5.1]{Del71}.)% where Deligne works with inifinite level Shimura varieties.

 The compositum field~$E\cdot E(G,X)\leq \CC$ is a finite extension of~$E$ which does not depend on the morphism~$\phi$. With our definition, it also does not depend on the compact open subgroups.
As a consequence, Galois representations for points in the same \emph{generalised} Hecke orbit can be deduced from each other, after passing to the \emph{same} finite extension $E\cdot E(G,X)/E$.

For future reference we summarise the above statements as follows.
\begin{proposition} \label{prop4.4}
We keep the same notations.
For any~$\sigma\in Gal(\ol{E}/E\cdot E(G,X))$, any~$g\in G(\AAA_f)$, and any~$\gamma\in G(\QQ)$, we have
\[
\sigma([\gamma\cdot \phi(x_0),g])=
[\gamma\cdot \phi(x_0), \rho'(\sigma)\cdot g]
\]
where
\[
\rho':=Ad_{\gamma}\circ \phi\circ \rho_{x_0}:\sigma\mapsto \gamma\cdot \phi\circ \rho_{x_0}(\sigma)\cdot \gamma^{-1}.
\]
is a Galois representation defined over~$E\cdot E(G,X)$ for~$\gamma\cdot \phi(x_0)$ in~$X$.
\end{proposition}
\begin{proof}
We may assume~$g=1$ by~\eqref{Hecke covariance}. This follow then from Prop.~\ref{Prop functoriality} applied to~$Ad_{\gamma}\circ\phi_0:M\to G\to G$.
\end{proof}
%\subsubsection*{Comments} Part of the definition of canonical models is that the Hecke correspondances are defined over~$E(G,X)$: it amouts to the property~$\sigma([x,g])=[x',g'g]$ if~$\sigma([x,1])=[x',g',g]$, for~$\sigma\in Aut(\CC/E(G,X))$. This the case~$\gamma=1$ and~$\phi=id:M\to M=G$. The observation~$[\gamma\cdot x,g]=[\gamma^{-1}\cdot g]$ and an easy computatiom implies the generalisation to~$\gamma\neq 1$. This also amount to the case where~$\phi:g\mapsto \gamma\cdot  g\cdot \gamma^{-1}:G\to G$ is a conjugation defined over~$\QQ$. To allow general morphisms~$\phi:M\to G$, one uses that the Shimura morphisms ... are defined over the compositum~$E(M,X_M)\cdot E(G,X)$. This is derived in Deligne by general properties of canonical models.
\subsection{Galois orbits vs Adelic orbits}\label{Galois orbits vs Adelic orbits}
Let~$U=\rho_{x_0}(Gal(\ol{E}/E))$. By definition we have
\begin{align*}
Gal(\ol{E}/E)\cdot [\phi\circ x_0,g]&= [\phi\circ x_0,\phi\circ\rho_{x_0}(Gal(\ol{E}/E))\cdot g]
\\&= G(\QQ)\backslash G(\QQ)\cdot \left(\{\phi\circ x_0\}\times \phi(U)\cdot g\right)\cdot K\cdot /K.
\end{align*}
%in which~$G(\QQ)$ is identified with its image in~$G(\AAA)$.

%Above identity is the key of our estimation of Galois orbits through adelic orbits~$g^{-1}\phi(U) gK_f/K_f$. The following statement let us controls the noninjectivity of the map
%\begin{equation}\label{adelic galois orbits map}
%\phi(U)\cdot g\cdot K_f/K_f\to G(\QQ)\backslash G(\QQ)\cdot \phi(U)\cdot g\cdot K_f\cdot /K_f.
%\end{equation}

The next proposition reduces the estimation of the size of the Galois orbit to that of the $\phi(U)$-orbit~$\phi(U)\cdot g\cdot K\cdot /K$.

\begin{proposition}\label{adelic galois orbits prop} There is a real number~$C\in\RR_{>0}$ such that
\[
\forall (\phi\circ x_0,g)\in \cH(x_0)\times G(\AAA_f), \frac{1}{C}\leq \frac{\abs{Gal(\ol{E}/E)\cdot [x_0,g]}}{[\phi(U):\phi(U)\cap K]}\leq 1.
\]
After possibly passing to a finite extension of~$E$, we may choose~$C=1$.
\end{proposition}
\begin{proof}We want to bound the cardinality of the fibres of the map
\begin{equation}\label{adelic galois orbits map}
\phi(U)\cdot g\cdot K/K\to
 G(\QQ)\backslash G(\QQ)\cdot
 \left(\{\phi\circ x_0\}\times \phi(U)\cdot g\right)\cdot K\cdot /K.
\end{equation}

We first describe the fibres. Let~$Z_\phi:=Z_G(\phi(M))$. The classical description of Hecke orbits gives an identity
\begin{align*}
G(\QQ)\backslash G(\QQ)\cdot \{\phi\circ x_0\}\times G(\AAA_f)
&\simeq Z_{\phi}(\QQ)\backslash \{\phi\circ x_0\}\times G(\AAA_f)\\
&\simeq \{\phi\circ x_0\}\times Z_{\phi}(\QQ)\backslash  G(\AAA_f).
\end{align*}
(This follows from~$G(\QQ)\cap Stab_{G(\RR)}(\phi\circ x_0)=Z_\phi(\QQ)$ in~$G(\RR)$.
We have embedded~$Z_{\phi}(\QQ)$ in~$G(\AAA)$ in the first line, and in~$G(\AAA_f)$
in the second line.)

%Above identification reduce the problem to bounding the fibers of
%\begin{equation}\label{galois adelic proof eq 1}
%\phi(U)\cdot g\cdot K_f/K_f\to Z_\phi(\QQ)\backslash Z_\phi(\QQ) \cdot \phi(U)\cdot g\cdot K_f/K_f
%\end{equation}
Define
\[
\Gamma= Z_\phi(\QQ)\cap \phi(U).
\]
The map~\eqref{adelic galois orbits map} can be written as a quotient map
\[
\phi(U)\cdot g\cdot K/K\to \Gamma\backslash(\phi(U)\cdot g\cdot K/K).
\] 
It will suffice to bound the order~$\abs{\Gamma}$.

The group~$Z_\phi(\QQ)$ is discrete in~$G(\AAA_f)$ because~$Z_\phi(\RR)$ is compact modulo~$Z(G)(\RR)$ and~$Z(G)(\QQ)$ is discrete in~$G(\AAA_f)$ (\cite[App. Lem. 5.13]{UY_General}),
where~$Z(G)$ is the centre of~$G$. As usual, we assume that~$G$ is the generic Mumford-Tate group on~$X$.
Therefore~$\Gamma$ is compact and discrete, and thus is finite.

We will realise~$\Gamma$ as a finite arithmetic group. We choose a faithful representation~$G\to GL(N)$ defined over~$\QQ$, and identify~$M$ and~$G$ with their images in~$GL(N)$.

%We use the same proof for both conclusions by choosing~$m=1$ for the first conclusion, and~$m=3$ for the second conclusion. 
We let~$K[m]=\ker(GL(N,\wh{\ZZ})\to GL(N,{\ZZ}/(m))$ for~$m\in\ZZ$.

There is a maximal compact subgroup~$K'$ of~$GL(N,\AAA_f)$ which contains~$K$. In~$GL(N,\AAA_f)$ all maximal compact subgroups are conjugated:~$K'$ is of the form~$h\cdot GL(N,\wh{\ZZ})\cdot  h^{-1}$ with~$h\in GL(N,\AAA_f)$. 
We may even choose~$h\in GL(N,\QQ)$ (this is a consequence of the fact that the class number of~$GL(N)/\QQ$ is one).
%, but for~$GL(N)$ is easily proved from standard theory of adelic linear lattices.) 

Conjugating the representation by~$h^{-1}$ we may assume~$h=1$:
we have
\[
U\leq K\leq GL(N,\widehat{\ZZ}).
\]
If~$m=3$ we pass to the finite extension of~$E$ corresponding to the subgroup~$\stackrel{-1}{\rho_{x_0}}(U\cap K[m])$ of $\Gal(\ol{E}/E)$.
 In any case we may assume
\[
U\leq K\cap K[m]\leq K[m].
\]

From~\ref{prop4.4}, we know that~$\phi=\gamma\phi_0 \gamma^{-1}$ for some~$\gamma\in GL(N,\QQ)$. It follows that
\[
\phi(U)\leq \gamma K[m] \gamma^{-1},
\]
and thus
\[
\Gamma = Z_\phi(\QQ)\cap \phi(U)\leq GL(N,\QQ)\cap \gamma K[m]\gamma^{-1}.
\]
Conjugating by~$\gamma^{-1}$ yields
\begin{align*}
\gamma^{-1}\cdot \Gamma\cdot \gamma
&\leq \gamma^{-1}GL(N,\QQ)\gamma\cap K[m]\\
&=
GL(N,\QQ)\cap K[m]\\
&=
\begin{cases}
GL(N,\ZZ)\text{ if~$m=1$,}\\
\ker\bigl(\,GL(N,\ZZ)\to GL(N,\ZZ/(3))\,\bigr)\text{ if $m=3$.}
\end{cases}
\end{align*}
Recall that~$\abs{\Gamma}=\abs{\gamma^{-1}\Gamma\gamma}$. We may thus conclude by applying the Lemma below to~$\gamma^{-1}\cdot \Gamma\cdot \gamma$. It follows that  for~$m=1$, $|\Gamma|$ is bounded independently of $\phi$ and  for~$m=3$, $|\Gamma|=1$.
\end{proof}
\begin{lem}For every~$N$, there is a real number~$C(N)$ such that, for every finite subgroup~$\Gamma\leq GL(N,\ZZ)$ we have
\[
\abs{\Gamma}\leq C(N),
\]
and if~$\Gamma\leq \ker\bigl(\,GL(N,\ZZ)\to GL(N,\ZZ/(3))\,\bigr)$ then~$\Gamma=1$.
\end{lem}
\begin{proof}From~\cite[Lemma 4.19.(Minkowski), p.\,232]{PR} the kernel has no nontrivial torsion. This implies the second assertion. 

This also implies that the reduction map~$GL(N,\ZZ)\to GL(N,\ZZ/(3))$ is injective on~$\Gamma$, thus inducing an embedding of $\Gamma$ in~$GL(N,\ZZ/(3))$. The first conclusion follows with
\[C(N)=\abs{GL(N,\ZZ/(3))}=\prod_{i=0}^{N-1}(3^N-3^i).\qedhere\]
\end{proof}

\section{Invariant Heights on Hecke orbits.}\label{section invariant height}

\subsection{Height functions}\label{subsection height generalities}

\subsubsection{Local affine height functions over~$\RR$ or~$\QQ_p$}
Let~$W$ be an affine variety over~$K=\RR$ or~$K=\QQ_p$. For every affine embedding defined over~$K$ 
\[
\iota_K:W\to \AAA^N_K
\]
there is an associated \emph{affine local Weil height} function~$H_{\iota_K}:W(K)\to \RR_{\geq 0}$ given by
\begin{equation}\label{height iotaK}
H_{\iota_K}(w)=\max\{1;\abs{w_1}_K;\ldots;\abs{w_N}_K\}.
\end{equation}
where~$\abs{-}_K$ is the standard absolute value on~$K$.

%We will not always require that~$W$ and~$\iota$ come form a variety and a morphism defined over~$\QQ$.
\subsubsection{Affine height functions over~$\QQ$}\label{section def HW}
When~$W$ and~$\iota:=\iota_K$ are defined over~$\QQ$, we can define, for~$w\in W(\QQ)$
\begin{align}\label{global iota height}
 H_\iota(w)&=H_{\iota\tens\RR}(w)\cdot H_{\iota,f}(w),\\
\text{with }H_{\iota,f}(w)&=\prod_{p} H_{\iota\tens\QQ_p}(w).
\end{align}
We define more generally, for~$w=(w_p)_p\in W(\AAA_f)$, 
\begin{equation}\label{notation:AF height}
H_{\iota,f}(w)=\prod_p H_{\iota\tens \QQ_p}(w_p).
\end{equation}

When~$W$ and the embedding~$\iota_{\RR}$, resp.~$\iota_{\QQ_p}$, resp.~$\iota$ are clear from the context, we will simply write
\begin{equation}\label{H sans iota}
H_{\RR}=H_{\iota_\RR},\quad H_p=H_{\iota_{\QQ_p}},\quad H_W=H_\iota\text{ and }H_f=H_{\iota,f}.
\end{equation}
Then~\eqref{global iota height} becomes
\begin{equation}\label{factorisation HW}
H_W=H_\RR\cdot H_f.
\end{equation}
%when~$\iota$ is defined over~$\RR$, $\QQ_p$ and~$\QQ$ respectively. 

%We have, on~$W(\QQ)$,
%\begin{equation}\label{factorisation HW}
%H_W=H_{\RR}\cdot H_{f}
%\end{equation}
%and we call~$H_{\RR}$ and~$H_{f}$ the ``archimedean part'' and ``finite part'' of the height function.

%When we say that ''we choose a height function on~$W$'', we will mean that we consider the above 
%functions for some implicit~$\iota$ that is chosen. 
\subsection{Polynomial equivalence and functoriality of Heights}\label{sec: funct heights}
We recall the \emph{functoriality} properties of heights. See~\cite{MW} or \cite{BG} for corresponding statements about projective Weil heights.
See~Def.~\ref{notation dominates} for the symbols~$\dom$ and~$\approx$.  
\begin{theorem}[Functoriality of Heights]\label{prop:functoriality heights}
 Let~$\phi:V\to V'$ be a morphism of affine varieties over~$\RR$, resp.~$\QQ_p$, resp.~$\QQ$,
and let
\[
\text{~$\iota_{\RR}:V\to \AAA^N_\RR$, resp.~$\iota_{\QQ_p}:V\to \AAA^N_{\QQ_p}$, resp.~$\iota:V\to \AAA^N_\QQ$}
\] be an affine embedding 
of~$V$, and let~$\iota'_{\RR}:V'\to \AAA^{N'}_\RR$, resp.~$\iota'_{\QQ_p}:V'\to \AAA^{N'}_{\QQ_p}$, resp.~$\iota':V'\to \AAA^{N'}_\QQ$ be an affine embedding of~$V'$.

Then, as functions on~$V(\RR)$, resp.~$V(\QQ_p)$, resp.~$V(\QQ)$ and~$V(\AAA_f)$,
\begin{multline*}
H_{\iota'_\RR}\circ \phi\dom H_{\iota_\RR}
,~\text{resp. }~
H_{\iota'_{\QQ_p}}\circ \phi\dom H_{\iota_{\QQ_p}}
,\\~\text{resp. }~
H_{\iota'}\circ \phi\dom H_{\iota}\text{ and }H_{\iota',f}\circ \phi\dom H_{\iota,f}.
\end{multline*}
\end{theorem}
\begin{cor}\label{coro:functoriality heights}
Let~$V$ be an affine algebraic variety over~$\RR$, resp.~$\QQ_p$, resp.~$\QQ$.
Let
\begin{multline*}
\text{$\iota_{\RR}:V\to \AAA^{N}_\RR$ and~$\iota'_{\RR}:V\to \AAA^{N'}_\RR$,}\\
\text{resp.~$\iota_{\QQ_p}:V\to \AAA^N_{\QQ_p}$ and~$\iota'_{\QQ_p}:V\to \AAA^{N'}_{\QQ_p}$,}\\
\text{resp.~$\iota:V\to \AAA^N_\QQ$ and~$\iota:V'\to \AAA^{N'}_\QQ$}
\end{multline*}
%  resp.~$\iota:V\to \AAA^N_\QQ$ and~$\iota:V'\to \AAA^{N'}_\QQ$
be affine embeddings of~$V'$.

Then, as functions on~$V(\RR)$, resp.~$V(\QQ_p)$, resp.~$V(\QQ)$ and~$V(\AAA_f)$,
\[
H_{\iota'_\RR}\approx H_{\iota_\RR}
,~\text{resp. }~
H_{\iota'_{\QQ_p}}\approx H_{\iota_{\QQ_p}}
,~\text{resp. }~
H_{\iota'}\approx H_{\iota}
~\text{ and }~
H_{\iota',f}\approx H_{\iota,f}.
\]
\end{cor}

\subsection{Galois invariant Height on the Hecke orbit}\label{subsection specific height function}
\label{sec: construction invariant height}
Let~$S=Sh_K(G,X)$ and~$x_0$ be as in~\S\ref{sec: main result}
and~$\rho_{x_0}:\Gal(\ol{E}/E)\to M(\AAA_f)$ be as in~\eqref{representation galoisienne}.
Let~$W=G\cdot \phi_0\subseteq \Hom(M,G)$ be the algebraic variety defined in~\S\ref{sec: def orbits}.
We have~$W\simeq G/Z_G(M)$. Let~$\Hom(\mathfrak{m},\mathfrak{g})$ be affine algebraic variety of linear maps $\mathfrak{m}\to\mathfrak{g}$.
As~$M$ is connected, we have an embedding
\[
\phi\mapsto d\phi:W\hookrightarrow \Hom(\mathfrak{m},\mathfrak{g}).
\]
As~$M$ is reductive, the image is closed, by~\cite{Ri-Conj}.

% and~$W$ is a closed subvariety of~$\Hom(\mathfrak{m},\mathfrak{g})$.

%In this section we construct a specific affine embedding~$\iota:W\to \AAA_\QQ^{\dim(G)\cdot\dim(M)}$ such that the associated height~\eqref{global iota height} satisfies the Galois invariance properties of Prop.~\ref{prop inv galois height}. This Galois invariance will be used in the step~\S\ref{preuve height bounds} of the proof of Th.~\ref{thm intro}.

We choose a lattice~$\mathfrak{g}_\ZZ\leq \mathfrak{g}$ such that
\[
\mathfrak{g}_\ZZ\tens\wh{\ZZ}\leq \mathfrak{g}\tens\AAA_f
\]
is stable under the action of~$K\leq G(\AAA_f)$. We define~$\mathfrak{m}_\ZZ=\mathfrak{g}_\ZZ\cap \mathfrak{m}$.
We choose a basis of~$ \mathfrak{g}$ which generates~$\mathfrak{g}_\ZZ$ and a basis of~$ \mathfrak{m}$ which generates~$\mathfrak{m}_\ZZ$. This choice induces an isomorphism
\[
j:\Hom(\mathfrak{m},\mathfrak{g})\xrightarrow{\sim} \QQ^{\dim(M)\cdot \dim(G)}.
\]
This induces an affine embedding
\[
\iota:=j\circ d:
W=G\cdot \phi_0\hookrightarrow \Hom(\mathfrak{m},\mathfrak{g})\to \AAA^{\dim(M)\cdot \dim(G)}_\QQ
\]
by first mapping~$\phi$ to~$d\phi$ and then to its matrix with respect to the bases we have chosen.

We denote by~$H_f:W(\AAA_f)\to \ZZ_{\geq1}$ and~$H'_f:\Hom(\mathfrak{m}_{\AAA_f},\mathfrak{g}_{\AAA_f})\to \ZZ_{\geq1}$
the functions given by~\eqref{notation:AF height} with respect to the embeddings~$\iota$ and~$j$.

\begin{proposition}[Galois invariance]
%Let~$S=Sh_K(G,X)$ and~$x_0$ be as in~\S\ref{sec: main result} and~$\rho_{x_0}:\Gal(\ol{E}/E)\to M(\AAA_f)$ be as in~\eqref{representation galoisienne}.
\label{prop inv galois height}

Let~$\phi_1,\phi_2\in W(\QQ)$ be such that~$s_1=[\phi_1\circ x_0,g_1]$ and~$s_2=[\phi_2\circ x_0,g_2]$
define points in~$\cH^g(s_0)$, where~$g_1,g_2\in G(\AAA_f)$, and assume that there exists a~$\sigma\in \Gal(\ol{E}/E)$ such that
\begin{equation}\label{eq prop invariant height}
\sigma(s_1)=s_2.
\end{equation}
Then
\[
H_{f}({g_1}^{-1}\phi_1 g_1)=H_{f}({g_2}^{-1}\phi_2 g_2).
\]
\end{proposition}
We first remark, from the formula
\begin{equation}\label{formula H'f}
\forall \psi\in \Hom(\mathfrak{m}_{\AAA_f},\mathfrak{g}_{\AAA_f}),~
H'_{f}(\psi)=\min\{n\in\ZZ_{\geq1} : n\cdot \psi(\mathfrak{m}_{\widehat{\ZZ}})\subset \mathfrak{g}_{\widehat{\ZZ}}\},
\end{equation}
that for every~$\widehat{\ZZ}$-module automorphism~$u:\mathfrak{m}_{\widehat{\ZZ}}\to \mathfrak{m}_{\widehat{\ZZ}}$ and~$k:\mathfrak{g}_{\widehat{\ZZ}}\to\mathfrak{g}_{\widehat{\ZZ}}$, we have
\begin{equation}\label{eq invariance gal et K}
H'_{f}(k\circ \psi\circ u)=H'_{f}(\psi).
\end{equation}
When~$\psi=d\phi$, and~$k=ad_{k'}:\mathfrak{g}\to\mathfrak{g}$ with~$k'\in K$ and~$u=ad_{u'}:\mathfrak{m}\to\mathfrak{m}$ with~$u'\in K\cap M(\AAA_f)$, this gives, with~$AD_{k'}\circ \phi\circ AD_{u'}:m\mapsto k'\phi(u'mu'^{-1})k'^{-1}$,
\begin{equation}\label{eq invariance gal et K 2}
H_{f}(AD_{k'}\circ \phi\circ AD_{u'})=H_{f}(\psi).
\end{equation}
\begin{proof}[Proof of Prop.~\ref{prop inv galois height}] We define~$u=\rho_{x_0}(\sigma)\in M(\AAA_f)\cap K$. From \eqref{eq prop invariant height} and the functoriality of Galois action Prop.~\ref{Prop functoriality} and~\ref{prop4.4}, we have
\[
[\phi_1\circ x_0,\phi_1(u)\cdot g_1]=[\phi_2\circ x_0,g_2].
\]
Hence, there exists~$\gamma\in G(\QQ)$ and~$k\in K$ such that
\[
(\gamma\cdot \phi_1\circ x_0,\gamma\cdot \phi_1(u)g_1 k)=(\phi_2\circ x_0,g_2).
\]
By Lemma~\ref{biglemma}, we also have~$\gamma\cdot \phi_1\cdot \gamma^{-1}=\phi_2$. Thus
\[
{g_2}^{-1}\phi_2 g_2=
\left(k^{-1}{g_1}^{-1}\phi_1(u)^{-1}\gamma^{-1}\right)
\cdot 
\left(\gamma\cdot \phi_1\cdot \gamma^{-1}\right)
\cdot 
\left(\gamma\cdot \phi_1(u)g_1 k\right)
\]
\[
=
k^{-1}{g_1}^{-1}\phi_1(u^{-1})\cdot \phi_1 \cdot \phi_1(u)g_1 k.
\]
We have
\[
\forall m\in M(\AAA_f),
\phi_1(u)^{-1}\cdot \phi_1(m) \cdot \phi_1(u)= \phi_1(u^{-1}mu)=\phi_1\circ AD_{u^{-1}}(m)
\]
and hence
\[
k^{-1}{g_1}^{-1}\phi_1(u)^{-1}\cdot \phi_1 \cdot \phi_1(u)g_1 k=AD_{k^{-1}}\circ({g_1}^{-1}\cdot \phi_1\cdot {g_1})\circ AD_{u^{-1}}.
\]
We finally have, using~\eqref{eq invariance gal et K 2},
\[
H_{f}({g_2}^{-1}\phi_2 g_2)=H_{f}(AD_{k^{-1}}\circ({g_1}^{-1}\cdot \phi_1\cdot {g_1})\circ AD_{u^{-1}})
=H_{f}({g_1}^{-1}\phi_1 g_1).\qedhere
\]
\end{proof}
%\[
%W(\RR)\times W(\AAA_f)\supseteq W(\QQ)^+\times (G(\AAA_f)\cdot \phi_0)\xleftarrow{} W(\QQ)^+\times G(\AAA_f) \to \cH^g(s_0).
%\]
\subsubsection{Height function on the generalised Hecke orbit~$\cH([x_0,1])$}
The function~$H'_f$ on~$\Hom(\mathfrak{m}_{\AAA_f},\mathfrak{g}_{\AAA_f})$ induces, at the level of the~$Sh_K(G,X)$, a function~$H_{s_0}$
on the generalised Hecke orbit~$\cH(s_0)$ of~$s_0:=[x_0,1]$, given as follows. For~$\phi\in \Hom(\mathfrak{m}_\QQ,\mathfrak{g}_\QQ)$ such that~$\phi\circ x_0\in X$ and~$g\in G(\AAA_f)$, we define
\[
H_{s_0}([\phi\circ x_0, g])=H'_{f}(d(g^{-1}\phi g)).
\]
The function~$H_{s_0}$ depends on the choices we have made, but different choices will produce the same function,
up to a bounded factor.

The case~$\sigma=1$ of Prop.~\ref{prop inv galois height} implies that~$H_{s_0}$ is well defined. 
Proposition~\ref{prop inv galois height} can then be rephrased as follows.
\begin{lem}For every~$\sigma\in \Gal(\ol{E}/E) $ and~$s\in \cH(s_0)$ we have
\begin{equation}\label{Galois invariance hecke height}
H_{s_0}(\sigma(s))=H_{s_0}(s).
\end{equation}
\end{lem}

\section{Height comparison on Siegel sets}\label{Section Siegel}
The main result of this section, Th.~\ref{theorem type Orr}, compares, for rational points of~$W=G/Z_G(M)$ contained in a given Siegel set (as in Def.~\ref{defi Siegel in W}), 
 the global height~$H_W$ of~\eqref{section def HW},
 with its factor~$H_{f}$ in~\eqref{factorisation HW}  (coming from the finite places).
The height~$H_W$  is the one appearing in our variant (Th.~\ref{PilaWilkie}) of the Theorem of Pila-Wilkie, and~$H_{f}$
 is the height appearing in our Galois bounds (see~Theorem~\ref{Galois lower bound}).

Our Th.~\ref{theorem type Orr} extends a result of Orr, in~\cite{OrrS}, which is only applicable to elements in~$G(\QQ)$. We work with elements of~$W(\QQ)$ instead. 
This is crucial to us as, in our strategy~\S\ref{Strategy outline}, we are working with geometric Hecke orbits.\footnote{
In general, when~$Z_G(M)\neq\{1\}$,  the height of an element~$g\in G(\QQ)$ is not bounded by the height of its image~$g\phi_0g^{-1}$ in~$W(\QQ)$ and
 not every element of~$W(\QQ)$ is the image of an element of~$G(\QQ)$.
}

This section develops different arguments than those of~\cite{OrrS}. 
They are more flexible, which allows us to obtain a more general result.

\subsection{Polynomial equivalence and archimedean height} 
We use Def.~\ref{notation dominates}.

\begin{lem}\label{Lemma properness real height} 
Let~$A\subset \RR^n$ be a semialgebraic subset, and let $f,g \colon A \lto \RR_{\geq 0}$ be semialgebraic and continuous functions.
Assume that $f$ is a proper map.

Then
$$
g \dom f.
$$
\end{lem}
\begin{proof}
We claim that the following function 
\begin{align*}
h \colon ]\inf_A(f), \infty[ &\to \RR_{\geq 0}\\
t &\mapsto \sup \{g(a): a\in A, f(a) \leq t\}
\end{align*}
is well defined. Fix an arbitrary~$t$ be in~$]\inf_A(f), \infty[$. The set $\{ a\in A: f(a) \leq t\}$ is compact since $f$ is proper.  It is nonempty since $t> \inf_A(f)$. As~$g$ is continuous,~$\{g(a): a\in A, f(a) \leq t\}$ is compact and nonempty, and its maximum belongs to~$\RR_{\geq 0}$, which proves the claim.

The function $h$ is also semialgebraic (see~\cite[Prop 2.2.4.]{BCR}). By~\cite{VDD} (\S~4.1 `Notes and comments' and references therein), $h$ is polynomially bounded.
The conclusion follows.
\end{proof}
The following uses Lemma~\ref{Lemma properness real height} for~$f$ and~$g$, and  again after swapping~$f$ and~$g$.
\begin{cor}\label{Coro properness real height}
 On a semialgebraic subset~$A\subset \RR^n$, two proper semialgebraic continuous functions~$f,g:A\to \RR$ are polynomially equivalent.
\end{cor}
We will also encounter the following situation.
\begin{lem}\label{Lemme equiv avec propre}
Let~$A\subset \RR^n$ and~$B\subset \RR^m$ be semialgebraic subsets, and~$f\colon A \to \RR_{\geq 0}$ and $g \colon B \to \RR_{\geq 0}$ be two proper semialgebraic continuous functions, and $p \colon A \to B$ be a proper and continuous semialgebraic function.
Then $g \circ p \approx f$.
\end{lem}
\begin{proof}We note that~$g\circ p$ is proper and continuous because~$g$ and~$p$ are. We can apply the Corollary~\ref{Coro properness real height} to~$f$ and~$g\circ p$.
\end{proof}
\subsubsection{}
%We first recover a well known property of heights. 
%For an affine algebraic variety~$V$ over $\RR$ and an affine embedding~$V\to\AAA^N$ defined over~$\RR$, we define the \emph{archimedean (affine Weil) height function associated to $\phi$} as
%\begin{equation}\label{associated archimedean Weil affine height}
%H_\phi(v) = \max\{1, \norm{\phi(v)}\}.
%\end{equation}
\begin{lem}
Let $V$ be an affine algebraic variety over $\RR$. 
Let $\phi:V\to \AAA^N$, and~$\phi':V\to\AAA^M$ be two closed embeddings, 
and let~$H_\phi$ and~$H_{\phi'}$ be defined as in~\eqref{height iotaK}.

Then~$H_\phi$ and~$H_{\phi'}$ are semialgebraic continuous proper functions, and
$$
H_{\phi} \approx H_{\phi'}.
$$
\end{lem}
\begin{proof} The real algebraic map $V(\RR)\to \RR^N$ induced by the Zariski closed embedding~$\phi$ is a closed embedding for the real topology. The functions~$\norm{~}_{\infty}:\RR^N\to\RR_{\geq0}$ and~$t\mapsto \max\{1;t\}:\RR_{\geq0}\xrightarrow{} \RR_{\geq0}$ are semialgebraic continuous proper maps. The composite map~$H_\phi$, and likewise~$H_{\phi'}$, are thus semialgebraic continuous and proper on~$V(\RR)$. We conclude with Cor.~\ref{Coro properness real height}.
\end{proof}
%We also apply our argument to a more complicated situation.
%\begin{cor}\label{coro points reels ouvert}
%Let $i \colon U \subset V$ be two affine algebraic varieties over~$\RR$ such that $U(\RR) = V(\RR)$.
%
%Let $\phi_U \colon U \lto \AAA^N$ and $\phi_V\colon V \lto \AAA^M$ be two closed embeddings and denote
%\begin{align*}
%H_{U,\RR}(x)&= \max(1, \norm{\phi_U(x)}):U(\RR)\to \RR_{\geq0}\\
%\text{ and }
%H_{V,\RR}(x)&= \max(1, \norm{\phi_V(x)}):V(\RR)\to \RR_{\geq0}
%\end{align*}
%the corresponding archimedean height functions.
%
%Then, as functions~$U(\RR)\to\RR_{\geq0}$, we have
%$$
%H_{U,\RR} \approx H_{V,\RR}.
%$$
%\end{cor}
%\begin{proof}From previous lemma,~$H_{U,\RR}$ is semialgebraic continuous proper on~$U(\RR)$, and~$H_{V,\RR}$ on~$V(\RR)$.
%But~$U(\RR)=V(\RR)$. We can use Corollary~\ref{Coro properness real height}.
%\end{proof}
\subsubsection{} 
%For an algebraic affine variety over $\RR$, \emph{an archimedean height function} will be a map~$H_{V,\RR}:V(\RR)\to \RR_{\geq0}$ %polynomially equivalent to a map of the form~\eqref{associated archimedean Weil affine height}.

\begin{lem}\label{Lemme propre sur sousensemble def}
Let~$p:U\to V$ be a morphism of affine algebraic varieties over~$\RR$,
and~$\phi_U:U\to\AAA_\RR^N$ and~$\phi_V:V\to\AAA^M_\RR$  be closed embeddings.
Let~$H_{\phi_U}$ and~$H_{\phi_V}$ be defined as in~\eqref{height iotaK}.
%and let $\phi_U \colon U \to \AAA^N$ and $\phi_V\colon V \to \AAA^M$ be two closed embeddings and 

%and let~$H_{U,\RR}$ and~$H_{V,\RR}$ be archimedean functions on~$U$ and~$V$.
%\begin{align*}
%H_{U,\RR}(x)&= \max(1, \norm{\phi_U(x)}):U(\RR)\to \RR_{\geq0}\\
%\text{ and }
%H_{V,\RR}(x)&= \max(1, \norm{\phi_V(x)}):V(\RR)\to \RR_{\geq0}
%\end{align*}
%the corresponding archimedean height functions.

Let~$A\subset U(\RR)$ be a semialgebraic subset such that~$p|_A:A\to V(\RR)$ is proper.
Then, as functions~$A\to\RR_{\geq0}$,
\[
H_{\phi_U}|_A\approx H_{\phi_V}\circ p|_A.
\]
\end{lem}
\begin{proof}We know that~$H_{\phi_U}$ and~$H_{\phi_V}$ are proper continuous and semialgebraic. As~$p|_A$ is proper,~$\iota:A\hookrightarrow U(\RR)$ is closed. It follows that~$H_{\phi_U}|_A=H_{\phi_U}\circ \iota$ is continuous, proper and semi-algebraic. We apply
Lemma~\ref{Lemme equiv avec propre} with~$A=U(\RR)$, and~$B=V(\RR)$, and~$p|_A$ as~$p$, and~$f=H_{\phi_U}|_A$ and~$g=H_{\phi_V}$.
\end{proof}
\subsection{Comparison of archimedean and finite height}
\begin{lem}\label{Lemma heights Gm}Let~$\iota:V\to\AAA^M$ be a \underline{closed} embedding with~$V = {\GG_m}^N$.
Then
\(
H_{\iota,\RR} \dom H_{\iota,f}\text{ on }\QQ^{\times N},
\)
where~$H_{\RR}$ and~$H_{f}$ are as in~\S\ref{section def HW}~\eqref{H sans iota}.
\end{lem}
\begin{proof}Thanks to Cor.~\ref{coro:functoriality heights}, we may substitute~$\iota$ with
the closed embedding
\begin{equation}\label{iotaN}
\iota_N: {\GG_m}^N\xrightarrow{(t_1,\ldots,t_N)\mapsto (t_1,{t_1}^{-1},\ldots,t_N,{t_N}^{-1})} \AAA^{2N}.
\end{equation}

We start with the case~$N=1$. We write an element~$t\in \QQ^\times$ as a reduced fraction~$n/m$. We can compute
\[
H_{\iota_1\tens\RR}(t)=\max\{|t|;|1/t|\}\text{ and }H_{\iota_1,f}(t)=|n\cdot m|.
\]
It follows~$H_{\iota_1\tens\RR}(t)\leq H_{\iota_1,f}(t)$.

For general~$\vec{t}=(t_1,\ldots,t_N)\in {\QQ^\times}^N$ there is some~$1\leq i\leq N$ such that
\[
H_{\iota_N\tens\RR}(\vec{t}\,)=\max\{|t_1|;|1/t_1|;\ldots;|t_N|;|1/t_N|\}=\max\{|t_i|;|1/t_i|\}.
\]
By the previous computation we have
\[
H_{\iota_N\tens\RR}({\vec{t}\,})=H_{\iota_1\tens\RR}(t_i)\leq H_{\iota_1,f}(t_i).
\]
We conclude by observing that
\[
H_{\iota_1,f}(t_i)\leq H_{\iota_N,f}(\vec{t}\,),
\]
as can be seen prime by prime.
\end{proof}
\begin{lem}
For $V = \GG_m^N\subset W=\AAA^N $, and affine embeddings~$\iota_V:V\to\AAA^M$, resp.~$\iota_W:W\to\AAA^{M'}$,
we have
\(
H_{\iota_V,f} \dom H_{\iota_W}\text{ on }{\QQ^\times}^N
\).
\end{lem}
\begin{proof}By Cor.~\ref{coro:functoriality heights}, we may assume that~$\iota_V$ is~$\iota_N$ of~\eqref{iotaN}, and that~$\iota_W$ is the identity map.
We can again reduce the problem to the case~$N=1$.
%$H_f(t_i)\leq H_f(\vec{t})$ to reduce to the case~$N=1$ and~$\vec{t}=t_i$.
 We write~$t_i=n/m$ as an irreducible fraction and then we have
\[
H_{\iota_V,f}(n/m)=\abs{n\cdot m}\leq \max\{\abs{n},\abs{m}\}^2=H_{\iota_W}(n/m)^2.\qedhere\]
\end{proof}
\begin{cor}\label{coro  heights Gm A1}
Let~$C\in\RR_{>0}$. We have
\[
H_{\iota_V,f} \dom H_{\iota_W,f}\text{ on }(\QQ^\times\cap [-C;C])^N.
\]
\end{cor}
\begin{proof} In the Lemma, we decompose~$H_{\iota_W}= H_{\iota_W\tens\RR}\cdot  H_{\iota_W,f}$. By hypothesis,~$ H_{\iota_W\tens\RR}\leq C$,
hence
\(
H_{\iota_W}\dom H_{\iota_W,f}\text{ on }(\QQ^\times\cap [-C;C])^N
\) which allows us to conclude.
\end{proof}

We  establish the following.
\begin{proposition}\label{Proposition equivalence heights}
Let $W$ be an affine variety over $\QQ$ and
let $p\colon W \lto \AAA^r$ be an algebraic map and
$\FS \subset W(\RR)$ be a semialgebraic closed subset such that
\begin{enumerate}
\item we have $p(\FS) \subseteq (\RR^\times)^r$;
\item the restriction $p |_{\FS} \colon \FS \lto {\RR^\times}^r$ is a proper map;
\item \label{hypo borne}		 the image~$p(\FS)$ is bounded in~$\RR^r$.
\end{enumerate}
We fix an affine embedding~$\iota$ of~$W$ and use notations~\eqref{H sans iota}. Then
$$
H_{\RR}|_{\FS \cap W(\QQ)} \dom H_{f} |_{\FS \cap W(\QQ)}.
$$
In particular 
$$
H_{W} |_{\FS \cap W(\QQ)}  \approx H_{f} |_{\FS \cap W(\QQ)}.
$$
\end{proposition}
\begin{proof}
We denote by~${\GG_m}^r\subset \AAA^r$ the affine open subset on which every
coordinate is invertible.

We fix affine embeddings~$\iota_W$ of~$W$, and~$\iota_{{\GG_m}^r}$ of~${\GG_m}^r$ and~$\iota_{\AAA^r}$ of~$\AAA^r$.

% Define
%\[
%U=\stackrel{-1}{p}({\GG_m}^r).
%\]
%Then~$U$ is an affine open subset of~$W$ defined over~$\QQ$: it is the intersection of the principal affine open subset where each of the coordinates of~$p$ are invertible.
%
%By the first hypothesis we have
%\begin{equation}\label{eq 1 du 711} 
%U(\RR)=W(\RR).
%\end{equation}
%By Corollary~\ref{coro points reels ouvert} we have
%\begin{equation}\label{eq 2 du 711} 
%H_{U,\RR}\approx H_{W,\RR}
%\end{equation}
%as functions on~$U(\RR)=W(\RR)$.
Because~$p|_{\FS}$ is continuous real algebraic, and (as a function to~$\RR^r$) is proper, by Lemma~\ref{Lemme propre sur sousensemble def} we have 
\begin{equation}\label{eq 12 du 711} 
H_{\iota_W\tens\RR}|_{\FS}\approx H_{\iota_{{\GG_m}^r}\tens\RR}\circ p|_\FS.
\end{equation}

By functoriality of heights Th.~\ref{prop:functoriality heights}, we have, on~$W(\QQ)$,
\begin{equation}\label{eq 3 du 711} 
H_{\iota_{{\AAA}^r},f}\circ p\dom H_{\iota_W,f}.
\end{equation}
As~$p(\FS)$ is bounded in~$\RR^r$, we have, 
by Lemma~\ref{Lemma heights Gm}, 
\begin{equation}\label{eq 4 du 711} 
H_{\iota_{{\GG_m}^r}\tens\RR}\dom H_{\iota_{{\GG_m}^r},f}.
\end{equation}
By hypothesis~\ref{hypo borne} we can use Corollary~\ref{coro  heights Gm A1} and get
\begin{equation}\label{eq 5 du 711} 
H_{\iota_{{\GG_m}^r},f}|_{p(\FS)\cap {\QQ^\times}^r}\dom H_{\iota_{{\AAA}^r},f}|_{p(\FS)\cap {\QQ^\times}^r}.
\end{equation}

%We now use the second hypothesis, apply Lemma~\ref{Lemme propre sur sousensemble def}, and get
%\begin{equation}\label{eq 6 du 711} 
%H_{\iota_W\tens\RR}|_{\FS}\approx H_{\iota_{{\GG_m}^r}\tens\RR}\circ p|_{\FS}
%\end{equation}
Combining these we get, using~\eqref{eq 12 du 711},
%~\eqref{eq 6 du 711},
~\eqref{eq 4 du 711},~\eqref{eq 5 du 711} and then~\eqref{eq 3 du 711},
\begin{multline*}
H_{\iota_W\tens\RR}|_{\FS\cap W(\QQ)}
\approx
H_{\iota_{{\GG_M}^r}\tens\RR}\circ p|_{\FS\cap W(\QQ)}\\
\dom 
H_{\iota_{{\GG_M}^r},f}\circ p|_{\FS\cap W(\QQ)}
\dom 
H_{\iota_{\AAA^r},f}\circ p|_{\FS\cap W(\QQ)}
\dom
H_{\iota_W,f}|_{\FS\cap W(\QQ)}.\qedhere
\end{multline*}
\end{proof}

\subsection{Construction of Siegel sets.}

We start by recalling some facts about parabolic subgroups and Siegel sets. A general reference is~\cite{BJ}.

Let $G_\QQ$ be a semisimple $\QQ$-algebraic group of adjoint type. We fix a minimal~$\QQ$-defined parabolic\footnote{Non necessarily proper: we have~$P_\QQ=G_\QQ$ when~$G_\QQ$ is of~$\QQ$-rank zero.} subgroup~ $P_\QQ$. Let~$G_\RR$ and~$P_\RR$ be the corresponding~$\RR$-algebraic groups.

Let~$X$ be the associated symmetric space\footnote{The space~$X$ in this section is of the form~$G(\RR)/K$ with~$K$ a \emph{non necessarily connected} maximal compact subgroup. This~$G(\RR)/K$ is connected and is a quotient of the space~$X=G(\RR)/K_\infty$ from  other sections of this article: when~$x$ is the image of~$x_0\in G(\RR)/K_\infty$, we have~$K_\infty=K^+$, where for simplicity we assume~$G$ is of adjoint type. The point~$x_0$ determines and is determined by a Hodge cocharacter~$h:\SSS\to G_\RR$, and the image point~$x\in G(\RR)/K$ determines and is determined by the corresponding Cartan involution~$\Theta=Ad_{h(i)}$.}, and choose~$x \in X$ and let $\Theta:G_{\RR}\to G_{\RR}$ be the Cartan involution associated with $x$. The orbit map~$g\mapsto g\cdot x$ induces the identification~$G(\RR)/K\simeq X$ where~$K$ is the maximal compact subgroup~$\{g\in G(\RR): g=\Theta(g)\}$. We denote~$K_\infty=K^+$ the neutral component.

We let~$N_\QQ$ be the unipotent radical of~$P_\QQ$: thus~$P_\QQ/N_\QQ$ is the maximal reductive quotient of~$P_\QQ$. The~$\RR$-algebraic group
$$
L:= P_{\RR} \cap \Theta(P_{\RR})
$$
is a maximal~$\RR$-algebraic reductive subgroup of~$P_{\RR}$ (cf. \cite[\S III.1.9]{BJ}), not necessarily defined over~$\QQ$, and the map~$L\to P_\RR\to (P_\QQ/N_\QQ)_\RR$ is an isomorphism. We denote by~$A'_\QQ$ the maximal central $\QQ$-split torus of~$P_\QQ/N_\QQ$, and define~$A\leq L$ as the inverse image of~$A_\RR$ in~$L$. We denote by~$A^+=A(\RR)^+$ the neutral component as a real Lie group.

We denote by~$\Phi$ the set of non-zero weights of the adjoint action of~$A$ on~$\mathfrak{g}\tens\RR$ (the ``(relative) roots''), and~$\Phi^+$ the subset of weights of the action on~$\mathfrak{n}\tens\RR$ (the ``positive'' ones). The eigenspaces are not necessarily defined over~$\QQ$.
 There exists a unique subset~$\Delta=\{\alpha_1;\ldots;\alpha_r\}\subset \Phi^+$ such that~$\alpha_1,\ldots,\alpha_r$ is a basis of~$X(A)=\Hom(A,{\GG_m}_\RR)$ and~$\Phi^+\subset \alpha_1\cdot\ZZ_{\geq 0} + \cdots +\alpha_r\cdot \ZZ_{\geq 0}$. The~$\alpha_i$ are known as the (relative) simple roots, and ~$r$ is equal to the~$\QQ$-rank of~$G_\QQ$.

The positive Weyl chamber in~$A^+$ is
\begin{equation}\label{defi Apos}
A^+_{\geq 0}=\{a\in A^+~:~ \forall 1\leq i\leq r, \alpha_i(a)\geq 1\}.
\end{equation}
We define
\begin{equation}\label{defi HP}
H_P:= \cap_{\chi \in X(P_\QQ)} \ker(\chi^2).
\end{equation}

We note that, for every one-dimensional representation~$\QQ\cdot \eta$ of~$H_P$, we have
\begin{equation}\label{HP pm}
\forall~h\in H_P(\RR),~h\cdot \eta\in \{+\eta;-\eta\}. 
\end{equation}

We first define Siegel sets in~$G_\QQ(\RR)$.
\begin{defini}[Siegel set]\label{defi Siegel in G}
A \emph{$\QQ$-Siegel set $\FS$ in~$G_\QQ(\RR)$ with respect to $P_\QQ$ and~$x$} is a subset~$\FS\subseteq G(\RR)$ of the following form.

There is a nonempty open and relatively compact subset~$\Omega\subseteq P_\QQ(\RR)$ and an element~$a\in A^+$ such that
$$
\FS = \Omega \cdot A^+_{\geq 0} \cdot a \cdot K_{\infty}.
$$

We will always assume that ~$\Omega\subseteq H_\QQ(\RR)$ and that $\Omega$ is semialgebraic.
\end{defini}

Usually Siegel sets are defined in~$G_\QQ(\RR)$ or in~$X=G_\QQ(\RR)/K$. Working with Geometric Hecke orbits as defined in~Def.~\ref{defi: geo orbit}, we use the variety~$W(\RR)^+=G(\RR)/Z_G(M)(\RR)$. We can view~$W(\RR)^+$ as an intermediary space in the sequence of maps~$G(\RR)\to W=G(\RR)/Z_G(M)(\RR) \to X$. The following definition allows us to work with Siegel sets in a greater generality.
\begin{defini}\label{defi Siegel in W}
Let~$Z$ be a compact subgroup of~$K$, and~$W=G_{\RR}/Z$. We define
\emph{a $\QQ$-Siegel set $\FS_W$ with respect to $P_\QQ$ and~$x$ in~$W^+:=G_{\RR}(\RR)/Z(\RR)$} to be the image of $\QQ$-Siegel set $\FS$ in~$G_\QQ(\RR)$ with respect to $P_\QQ$ and~$x$.

\end{defini}
We note that if~$Z$ is defined over~$\QQ$ then so is~$W$ and we can consider the subset~$W^+(\QQ)\cap \FS_W$.
\subsection{Divergence in Siegel sets} In the rest of this section we use the notation~$G=G_\QQ$.
\subsubsection{}\label{divergent and proper}
We say that an infinite sequence, in an appropriate topological space, is \emph{divergent} if it does not contain an infinite convergent subsequence. 
A continuous map is \emph{proper} if it maps divergent sequences to divergent sequences.
\subsubsection{}
We will use the closure of a Siegel set.
\begin{proposition}Consider~$\FS=\Omega\cdot A^+_{\geq 0}\cdot a\cdot K_{\infty}$ as in Def.~\ref{defi Siegel in G}.

Then its closure in~$G(\RR)$ is given by
\begin{equation}\label{FS bar}
\ol{\FS}=\ol{\Omega}\cdot A^+_{\geq 0}\cdot a\cdot K_{\infty}.
\end{equation}
and~$\ol{\FS}$ is contained in a~$\QQ$-Siegel set $\FS'$ in~$G(\RR)$ with respect to $P$ and~$x$.
\end{proposition}
\begin{proof}
The set~$\FS$ is obviously dense in the right-hand side of~\eqref{FS bar}. It is  the image of the proper map in Lemma~\ref{lemma proper mult},
and thus it is a closed subset in~$G(\RR)$. This proves the first assertion. Let~$U$ be a nonempty relatively compact semialgebraic open neighbourhood of~$1$ in~$H(\RR)$, for instance a small euclidean open ball in a faithful representation~$H\to GL(N)$. Then~$\Omega'=U\cdot \ol{\Omega}$ is an open relatively compact semialgebraic open neighbourhood of~$\ol{\Omega}$ in~$H(\RR)$, and the Siegel set~$\Omega'\cdot A^+_{\geq 0}\cdot a\cdot K$ contains~$\ol{\FS}$.
\end{proof}
We used the following.
\begin{lem}\label{lemma proper mult}
The map
\[
(\omega, a, k)\mapsto \omega\cdot a\cdot k:
\ol{\Omega}\times A^+_{\geq 0}\cdot a\times K\to G(\RR)
\]
is proper.
\end{lem}
\begin{proof}It suffices to prove that the image of every divergent sequence in the left hand-side is not a convergent sequence in the right hand-side.
We prove the contrapositive.

Let~$(\omega_n, a_n, k_n)_{n\in\ZZ_{\geq1}}$ be a sequence in the left-hand-side such that~$(\omega_n\cdot a_n\cdot k_n)_{n\in\ZZ_{\geq1}}$ is a convergent sequence in~$G(\RR)$. Because~$\ol{\Omega}$
and~$K$ are compact, after possibly extracting a subsequence we may assume that~$(\omega_n)_{n\in\ZZ_{\geq1}}$ and~$(k_n)_{n\in\ZZ_{\geq1}}$ are convergent sequences. It follows that~$(a_n)_{n\in\ZZ_{\geq1}}$ is a convergent subsequence. We
recall that~$(\alpha_1,\ldots,\alpha_r):A^+\to{\RR_{>0}}^r$ is an isomorphism. It follows that~$A^+_{\geq 0}$ is closed in~$A^+$.
Because~$A(\RR)$ is a closed subgroup of~$G(\RR)$, that~$A^+$ is a closed subgroup of~$A(\RR)$, this~$A_{\geq0}^+$ is closed in~$G(\RR)$. One deduces that~$A_{\geq0}^+\cdot a$ is closed in~$G(\RR)$ and that the limit of~$(a_n)_{n\in\ZZ_{\geq1}}$ in~$G(\RR)$ 
belongs to~$A_{\geq0}^+\cdot a$. 

We proved that the original sequence~$(\omega_n, a_n, k_n)_{n\in\ZZ_{\geq1}}$ contains a convergent infinite subsequence. Thus it is not a divergent sequence.
\end{proof}
These results have the following consequence.
\begin{cor}\label{coro chi criterion}
 A sequence~$(\omega_n\cdot a_n\cdot k_n)_{n\in\ZZ_{\geq1}}$ is divergent in~$\ol{\FS}$ if and only if~$a_n$
is divergent in~$A^+_{\geq 0}$.
\end{cor}
\begin{proof}Because~$\ol{\FS}$ is closed,~$(\omega_n\cdot a_n\cdot k_n)_{n\in\ZZ_{\geq1}}$ is also divergent in~$G(\RR)$. It 
follows that the sequence~$(\omega_n, a_n, k_n)_{n\in\ZZ_{\geq1}}$ contains no convergent subsequence.
Because~$\ol{\Omega}$ and~$K$ are compact, the projection map~
$$\ol{\Omega}\times A^+_{\geq 0}\cdot a\times K\to A^+_{\geq 0}\cdot a$$
is proper. It follows that the image sequence~$(a_n)_{n\in\ZZ_{\geq1}}$ is divergent in~$ A^+_{\geq 0}\cdot a$.
\end{proof}

%\subsubsection{A DEPLACER}
%A fundamental property of $\QQ$-Siegel sets is that if $\Omega$ is large enough and the~$\alpha_i(a)$ are small enough, there exists a finite set 
%$F \subset G(\QQ)$ such that
%$$
%\Gamma F \FS = G(\RR)
%$$
%For~$f\in G(\QQ)$, each~$f \cdot \FS$ is also a $\QQ$-Siegel set in~$G$, with respect to~$P'=fPf^{-1}$, and for~$f\Omega f^{-1}$ and~$faf^{-1}$.

\subsubsection{}\label{numero eta}
Let~$P_{1},\ldots,P_{r}$ be the maximal $\QQ$-defined proper\footnote{There are none if~$r=0$.} parabolic subgroups of~$G$ containing~$P$. We denote by~$N_{i}$ their unipotent radicals,
and~$\mathfrak{n}_i$ the~($\QQ$-linear) Lie algebra of~$N_{i}$. The adjoint representation of~$G$ induces an action of~$G$ on the~$\QQ$-vector space~$V_i=\bigwedge^{\dim(N_i)}\mathfrak{g}$. Then the $\QQ$-vector subspace~$\bigwedge^{\dim(N_i)}\mathfrak{n}_i\leq V_i$ is of dimension~$1$, and we choose a generator~$\eta_i$ of this~$\QQ$-line.

Then the $\RR$-line~$\RR \cdot \eta_i$ is an eigenspace of~$A$ acting on~$V_i\tens\RR$, and this eigenspace is defined over~$\QQ$. Let~$\chi_i$ be the corresponding eigencharacters of~$A$: we have
\begin{equation}\label{eigencharacter}
\forall a\in A(\RR),\forall 1\leq i\leq r, a\cdot \eta_i =\chi_i(a)\cdot \eta_i.
\end{equation}
For~$1\leq i\leq r$ the~$\chi_i$ are positive multiples~$k_1\cdot \omega_1,\ldots,k_r\cdot \omega_r$ of the (relative) fundamental weights\footnote{The ``weights lattice''~$\omega_1\cdot\ZZ+\ldots+\omega_r\cdot \ZZ\supset X(A)$ can be identified with~$X(\wt{A})$ where~$\wt{A}$ is the torus in a simply connected cover~$\wt{G}_\RR\to G_\RR$ which maps onto~$A$.}~$\omega_1,\ldots,\omega_r\in X(A)\tens \QQ$. In particular
\begin{equation}\label{chi positivity}
\forall \alpha\in A^+_{\geq 0}, \forall 1\leq i\leq r, \, \chi_i(\alpha)^{-1}\leq 1.
\end{equation}

One knows that the fundamental weights are positive $\QQ$-linear combinations of the~$\alpha_i$ and that they form a basis of~$X(A)\tens\QQ$. The same holds for the~$\chi_i$. We deduce the following.
\begin{lem}\label{lem chi criterion}	
 Let~$(a_n)_{n\in\ZZ_{\geq0}}$ be a sequence in~$A^+_{\geq 0}\cdot a$. Then the sequence is divergent (no infinite subsequence is convergent) if and only if
\begin{equation}\label{omega criterion}
\lim_{n\to\infty} \min_{1\leq i\leq r} \chi_i(a_n)^{-1}=0.
\end{equation}
\end{lem}	
\begin{proof} If~$(a_n)_{n\in\ZZ_{\geq0}}$ contains a convergent infinite subsequence, then the sequence $(\min_{1\leq i\leq r} \chi_i(a_n)^{-1})_{n\in\ZZ_{\geq0}}$ contains a convergent infinite subsequence in~$\RR_{>0}$ and we cannot have~\eqref{omega criterion}.

This proves one implication and we now prove the other.

Assume that~\eqref{omega criterion} fails. Equivalently
\[
L:=\limsup_{n\to\infty} \min_{1\leq i\leq r} \chi_i(a_n)^{-1}>0.
\] 
After possibly extracting a subsequence, we have
\begin{equation}\label{chi criterion proof}
\lim_{n\to\infty} \min_{1\leq i\leq r} \chi_i(a_n)^{-1}=L\neq  0.
\end{equation}
Because the~$a_n$ belong to~$A^+_{\geq 0}\cdot a$ we have~$\sup_{n\in\ZZ_{\geq0}}\alpha_i(a_n)^{-1}\leq \alpha_i(a)^{-1}$ for every~$1\leq i\leq r$.
Because the~$\chi_i$ are positive linear combination of the~$\alpha_i$, the~$\chi_i(a_n)^{-1}$ are bounded above. 
According to~\eqref{chi criterion proof} they are bounded below in~$\RR_{>0}$. Because the~$\chi_i$ form a basis of~$X(A)\tens\QQ$,
the~$\alpha_i$ are linear combination of the~$\chi_i$. Hence the~$\alpha_i(a_n)$ are bounded above and below in~$\RR_{>0}$.
Equivalently~$(a_n)_{n\in\ZZ_{\geq0}}$ is bounded in~$A^+$.
Hence~$(a_n)_{n\in\ZZ_{\geq0}}$ is not divergent.

This proves the other implication.
\end{proof}

\subsection{Height on Siegel sets}
The following statement is the main objective of~\S\,\ref{Section Siegel}.

\begin{theorem}\label{theorem type Orr}
Let~$\FS_W$ as in Def.~\ref{defi Siegel in W} with~$Z$ defined over~$\QQ$, let~$\iota:G/Z\to \AAA_\QQ^N$ be an affine embedding and let~$H_W=H_{\RR}\cdot H_{f}$ be as in~\eqref{factorisation HW}. 
Then, as functions~$\ol{\FS_W}\cap W(\QQ)\to \RR_{\geq 0}$  we have
$$
H_W \approx H_{f}.
$$
\end{theorem}
%In proving the equivalence, we may replace~$\FS_W$ by a bigger set. We will prove the equivalence on~$  \ol{\FS_W}\cap W(\QQ)$.
%
%Let~$\Omega\subseteq H_P(\RR)$ and~$a$ be as in the definition such that~$\FS_W=\FS/Z(\RR)= \Omega \cdot A_{\geq 0} \cdot a \cdot K_{\infty}/Z(\RR)$.
%
%Because~$a\cdot A_{\geq 0}$ is closed and~$K_\infty\times \ol{\Omega}$ is compact the subset
%\[
%\]

%We observe that the map
%\[
%(\omega,\alpha,k)\mapsto w(\omega,\alpha,k):\ol{\Omega}\times a\cdot A_{\geq 0} \times K_\infty 
%\to G(\RR)\to W(\RR)
%\]
%is proper because~$a\cdot A_{\geq 0}$ is closed and~$K_\infty\times \ol{\Omega}$ is compact. It follows
%\begin{itemize}
%\item that~$\ol{\FS_W}=\ol{\Omega}\cdot a\cdot A_{\geq 0} \cdot K_\infty/Z$
%\item that~$w_n=w(\omega_n,\alpha_n,k_n)$
%diverges to infinity in~$\ol{\FS_W}$ if and only if~$\alpha_n$ diverges to infinity in~$A_{\geq0}$.
%\end{itemize}
%We note that the latter means 
%\begin{equation}\label{corner divergence}
%\lim_{n\in\ZZ_{\geq0}}\min_{i=1,\ldots,r}\chi_i(a_n)=0.
%\end{equation}

This will be deduced from Proposition~\ref{Proposition equivalence heights}. We first construct the map~$p$ to which we apply the Proposition, and then verify the assumptions of the Proposition.

\subsubsection{Construction of the morphism~$p$} Let the~$\eta_i\in V_i=\bigwedge^{\dim{\Fn_i}} \Fg$ and the~$\chi_i$ be as in \S\ref{numero eta}.

%We use the minimal~$\QQ$-parabolic subgroup~$P$ of~$G$ and the associated~$\eta_i$ 
%from above. 
%
%We construct our map 
%$$
%p\colon W \lto \AAA^r_{\QQ}
%$$
%as follows.

For each $1\leq i\leq r$, we choose a positive definite quadratic form
$$
Q'_i \colon V_i \lto \QQ.
$$
We denote by~$dz$ the Haar probability measure on~$Z(\RR)$ we define the real quadratic form 
\begin{equation}\label{construction Q'}
Q_i(v) = \int_{Z(\RR)} Q'_i(z\cdot v)~ dz:V_i\tens\RR\to \RR.
\end{equation}
%where we use the induced representation on the space of quadratic forms.
The following is central in our argument.
\begin{lem}\label{lem Qi}
The quadratic form~$Q_i$ is invariant under~$Z(\RR)$, is  positive definite, and is defined over~$\QQ$.
\end{lem}
\begin{proof}The two first properties are immediate from~\eqref{construction Q'}. We prove that~$Q_i$ is defined over~$\QQ$.
Let~$V$ be the $\QQ$-vector space of quadratic forms, as a representation of~$Z$, and~$V^Z$ be the subspace of elements fixed by~$Z$.
As~$Z(\RR)$ is compact, the~$\QQ$-group~$Z$ is reductive, and  there is a~$Z$-stable
$\QQ$-subspace~$W$ such that we can decompose
\[
V=V^Z\oplus W.
\]
Let us write correspondingly
\[
Q'_i=Q'_Z+Q'_W
\]
with~$Q'_Z$ in~$V^Z$ and~$Q'_W$ in $W$. 

%We similarly decompose, in~$V\tens\RR$,
%\[
%Q_i=Q_Z+Q_W
%\]
%with~$Q_Z\in V^Z\tens\RR$ and~$Q_W\in W\tens\RR$.

Because~$Q_Z$ is invariant under~$Z(\RR)$ and~$W\tens\RR$ is stable under~$Z(\RR)$,
\[
\int_{Z(\RR)} z \cdot Q_Z\,dz=Q_Z\text{ and }\int_{Z(\RR)} z \cdot Q_W\,dz\in W\tens\RR.
\]
By construction~$\int_{Z(\RR)} z \cdot Q_W\,dz$ is fixed by~$Z(\RR)$,
and thus it belongs to~$W\tens\RR\cap V^Z\tens\RR=\{0\}$. We compute
\[
Q_i
=\int_{Z(\RR)} z \cdot Q'_i\,dz
=\int_{Z(\RR)} z \cdot Q_Z\,dz
+
\int_{Z(\RR)} z \cdot Q_W\,dz
=Q_Z+0.
\]
Because~$Q_Z$ is defined over~$\QQ$, so is~$Q_i$.
%We denote~$\pi:V\tens \RR\to V\tens\RR$ the operator~$Q\mapsto \int_{Z(\RR)} z \cdot Q dz$.
%It preserves the subrepresentations~$V^Z\tens\RR$ and~$W\tens\RR$. On~$V^Z\tens\RR$,
%one readily see~$\pi(Q^Z)=Q^Z$. We also observe that~$\pi(Q^W)$ is in~$W\tens\RR$ and is
%invariant under~$Z$. Hence~$Q^W\in W\tens\RR\cap V^Z\tens\RR=\{0\}$. That is~$Q^W=0$,
%and indeed~$Q_i=Q^Z$, which is defined over~$\QQ$.
\end{proof}

We can now define~$p:W\to \AAA_\QQ^r$ by 
\begin{equation}\label{defi p Siegel}
p(gZ)= (Q_1(g^{-1}\cdot \eta_1) , \dots , Q_r(g^{-1}\cdot \eta_r)).
\end{equation}
As the~$Q_i$ are defined over~$\QQ$ so is~$p$, and as the~$Q_i$
are~$Z$-invariant,~$p$ is well defined.

\subsubsection{Properties of the morphism~$p$}
Our next task is to verify the assumptions of Proposition~\ref{Proposition equivalence heights}.
\begin{enumerate}
%\item \label{first point} The application~$p:$ is defined over $\QQ$.
\item \label{second point} We have~$p(\ol{\FS_W}) \subset (\RR^\times)^r$.
\item \label{third point} As a map~$\ol{\FS_W}\to(\RR^\times)^r$, 
$
p|_{\ol{\FS_W}} 
$
is proper with respect to the real topologies.
\item \label{fourth point} the image~$p(\ol{\FS})$ is bounded in~$\RR^r$.
\end{enumerate}

%\begin{proof}[Proof of~\ref{first point}]
%We already proved that~$p$ is defined over~$\QQ$.
%\end{proof}

\begin{proof}[Proof of~\ref{second point}]
%For the first point we need to check that every coordinate~$p_i=Q_i(g^{-1}\eta_1)$ is non zero on~$\ol{\FS_W}$. 
Every point of~$\ol{\FS_W}$ is of the form~$gZ$ with~$g\in G(\RR)$.
The vector~$g\cdot \eta_i$ is thus in~$V_i\tens\RR$. As~$\eta_i\neq0$ and~$g$ is invertible,~$g^{-1}\cdot \eta_i\neq 0$.
As~$Q_i$ is positive definite by Lem.~\ref{lem Qi}, we have~$p_i(g):=Q_i(g^{-1}\eta_i)\in\RR_{>0}$.
%we deduce that~$Q_i$ is definite positive: it takes non zero value on non zero real vectors.
%Every point of~$\ol{\FS_W}$ is of the form~$Zg$ with~$g\in G(\RR)$.
%Because~$g\in G(\RR)$, th vector~$g\cdot \eta_i$ is in~$V_i\tens\RR$, and because th action of~$g$ is invertible,~$g\cdot \eta_i$ is non zero. Finally~$Q_i(g^{-1}\eta_i)\in\RR_{>0}\subset \RR^\times$.
\end{proof}
We will use that there exists~$C\in\RR_{>0}$ such that for every~$1\leq i\leq r$,
\begin{equation}\label{Siegel claim}
\forall (h,\alpha,k)\in H_P\times A^+_{\geq 0}\times K_\infty,\,
0\leq p_i(h\cdot a\cdot \alpha\cdot k)\leq C\cdot \chi_i(\alpha)^{-2}.
\end{equation}
%We first prove this claim.
\begin{proof}[Proof of~\eqref{Siegel claim}] %We write~$\sigma=h\cdot a\cdot \alpha\cdot k$.
We write~$\sigma=h\cdot a\cdot \alpha\cdot k$. 
By~\eqref{HP pm}, we have~$h^{-1}\cdot \eta_i=\pm\eta_i$. Thus
\[
\sigma^{-1}\cdot \eta_i=\pm k^{-1}\cdot \alpha^{-1}\cdot  a^{-1}\cdot \eta_i=\pm k^{-1} \cdot \chi_i(\alpha)^{-1}\cdot \eta_i.
\]
Because~$K$ is compact there exists a~$K$-invariant euclidean norm~$\norm{-}$ on~$V_i\tens\RR$. The two norms~$\sqrt{Q_i}$ and~$\norm{-}$ on~$V_i\tens\RR$ are comparable: there is~$C_i\in\RR_{>0}$ such that for any~$v\in V_i\tens \RR$,
\[
Q_i(v)\leq C_i\cdot \norm{v}^2.
\]
We deduce~\eqref{Siegel claim} with~$C=\max_{i\in\{1;\ldots;r\}} C_i\cdot \norm{a^{-1}\cdot \eta_i}^2$ from
\begin{align*}
p_i(\sigma)&=Q_i(\sigma^{-1}\cdot  \eta_i)\\
&\leq C_i\cdot \norm{\pm k^{-1} \cdot \chi_i(\alpha)^{-1}\cdot a^{-1}\cdot \eta_i}^2\\
&= C_i\cdot \norm{\pm\chi_i(\alpha)^{-1}\cdot a^{-1}\cdot \eta_i}^2\\
&=C_i\cdot \chi_i(\alpha)^{-2} \cdot \norm{a^{-1}\cdot \eta_i}^2\geq 0.\qedhere
\end{align*}
\end{proof}

\begin{proof}[Proof of~\ref{third point}] %This point uses the explicit description of Siegel sets.
%Write~$\FS_W=\FS/Z(\RR)$ and~$\sigma\in \FS$
%as
%\[
%\sigma = \omega a k\text{ with }\omega\in \Omega, a\in A_{\geq a_0}, k\in K_\infty.
%\]
%Because~$\Omega \subset H_P\subset \ker {\chi_i}^2$ we have~$\sigma^{-1}\cdot \eta_i=\pm k^{-1}\cdot a^{-1}\cdot \eta_i$
%and~$Q_i(\sigma^{-1}\cdot \eta_i)=Q_i(k^{-1}\cdot a^{-1}\cdot \eta_i)$.
%

For a divergent sequence~$\sigma_n=\omega_n a\cdot \alpha_n k_n$ in~$\ol{\FS_W}$, Cor.~\ref{coro chi criterion} and Lem.~\ref{lem chi criterion} imply that
\(
\min_{1\leq i \leq r} \chi_i(\alpha_n) \to 0.
\)
Using~\eqref{Siegel claim}, we deduce that
\begin{equation}\label{p diverge vers 0}
\min_{1\leq i \leq r} p_i(\omega_n a_n a k_n) \to 0.
\end{equation}
and thus~$p(\omega_n a_n k_n)$ is divergent in~${\RR^\times}^r$. This proves the properness.
\end{proof}

\begin{proof}[Proof of~\ref{fourth point}] For~$\sigma=h\cdot \alpha\cdot a\cdot k\in\ol{\FS}$, we have~$0\leq \chi_1(\alpha),\ldots,\chi_r(\alpha)\leq 1$ by~\eqref{chi positivity}, and deduce from~\eqref{Siegel claim} that
\[
\max_{1\leq i \leq r} \abs{p_i(\omega_n a_n a k_n)}\leq C.\qedhere
\]
\end{proof}
%\subsubsection{Conclusion}
We now use Proposition~\ref{Proposition equivalence heights} with~$\FS=\ol{\FS_W}$. This concludes the proof of Theorem~\ref{theorem type Orr}.

\section{Weak adelic Mumford-Tate hypothesis and Lower bounds on Galois orbits.}\label{section MT adelic}

This section is central to the proof of the Andr\'e-Pink-Zannier conjecture under our assumptions (Theorem~\ref{thm intro}). 
In this section, we state a precise form of the ``weak adelic Mumford-Tate hypothesis''.

We then translate lower and upper bounds on adelic orbits of Appendices~\ref{section adelic bounds} and~\ref{app c} into estimates for sizes of the  Galois orbits in terms of the height functions of~\S\ref{section invariant height} (when the 
 Mumford-Tate hypothesis holds).

For simplicity we will in the sequel refer to the `Mumford-Tate hypothesis' or simply `MT hypothesis'.

	Finally we provide some natural functoriality properties of the Mumford-Tate hypothesis, which will be needed for the reduction steps in the proof of our main theorem.
\subsection{The Mumford-Tate hypothesis}
We start with a property applicable in more general situations.
%general compact adelic subgroups of $M(\AAA_f)$.
\begin{defini}\label{defi MT in M}	
Let~$M$ be a linear algebraic group over~$\QQ$,
let~$K_M\leq M(\AAA_f)$ be a compact open subgroup, and let~$U\leq M(\AAA_f)$ be a compact subgroup. 

We say that~$U$ is \emph{MT in~$M$} if the indices
\begin{equation}\label{eq def MT}
[K_M \cap M(\QQ_p):U\cap K_M\cap M(\QQ_p)]
\end{equation}
are finite and uniformly bounded as~$p$ ranges through primes, where~$M(\QQ_p)\leq M(\AAA_f)$ is understood as a factor subgroup of~$M(\AAA_f)$.
%We will often say that ``$U$ is MT in $M$''.
\end{defini}
Note that the definition does not depend on the choice of~$K_M$, as any two compact open subgroups are commensurable. We may always enlarge~$K_M$ so that it takes the product form~$K_M=\prod_p K_p$ in which case the indices become
\[
[K_p:U\cap K_p].
\]
Likewise if~$U'\leq M(\AAA_f)$ is a compact subgroup commensurable to~$U$, then~$U$ is MT in~$M$ if and only if~$U'$ is MT in~$M$. 
Note (and this is very important) that the condition that~$U$ is MT in~$M$ \emph{does not imply} that~$U$ is open in~$M(\AAA_f)$. 

%This last requirement is too strong to apply to Galois representations in general: there are known obstruction and precise conjecture can be found in~\cite{SerreConj}. Our hypothesis is weaker than these conjectures.
The following observation is an immediate consequence of the definition.
\begin{lem}\label{MT in product form}
 In the Definition~\ref{defi MT in M}, let
\begin{equation}\label{eq defi U' Up}
	\text{$U_p:=U\cap K_M\cap M(\QQ_p)$ and~$U'=\prod_p U_p\leq \prod_p M(\QQ_p)$. }
\end{equation}
Then~$U$ is MT in~$M$ if and only if~$U'$ is MT in~$M$.
\end{lem}

We specialise the above definition to the context of images of Galois representations.
\begin{defini}\label{MT hypothesis}
Let~$(G,X)$ be a Shimura datum, let~$x_0\in X$ and let~$M$ be the Mumford-Tate group of~$x_0$, let~$\rho_{x_0}$ be a Galois representation for~$x_0$ defined over a field~$E$ in the sense of Definition~\ref{terminologie galois repr}, and let~$U=\rho_{x_0}(Gal(\ol{E}/E))\leq M(\AAA_f)$ be the image of~$\rho_{x_0}$. 
\begin{enumerate}
\item We say that~$x_0$ \emph{satisfies the MT hypothesis}, if~$U$ is MT in~$M$.
\item Let~$K\leq G(\AAA_f)$ be a compact open subgroup, and~$s_0=[x_0,1]\in Sh_K(G,X)$.
We say that~\emph{$s_0$ satisfies the MT hypothesis} if~$U$ is MT in~$M$.
\end{enumerate}
\end{defini}

\subsection{Lower bounds for Galois orbits in terms of finite heights under the MT hypothesis} The following statement is an essential ingredient  in the proof of the main theorem~\S\ref{Main proof}.
We also believe it to be of independent interest.

\begin{theorem}\label{Galois lower bound}
 Let~$M$ be a connected reductive group over~$\QQ$ and~$U\leq M(\AAA_f)$ 
be a compact subgroup which is MT in~$M$ in the sense of Def.~\ref{defi MT in M}. We use the notations of Def.~\ref{notation dominates}.

\begin{enumerate}
\item Let~$\phi_0:M\to GL(N)$ be a representation over~$\QQ$, and let~$W$ be the~$GL(N)$-conjugacy class of~$\phi_0$. 
We consider an affine embedding~$\iota:W\to \AAA^N_{\QQ}$ and the corresponding function~$H_f:W(\QQ_f)\to\ZZ_{\geq1}$ defined by~\eqref{notation:AF height}. Then, as $\phi$ ranges through $W(\AAA_f)$, we have
\begin{equation}\label{Galois bounds linear}
[\phi(U):\phi(U) \cap GL(N,\widehat{\ZZ})] \approx H_f(\phi).
\end{equation}
\item Let~$\phi_0:M\to G$ be a morphism of algebraic groups over~$\QQ$ and let~$W$ 
be the~$G$-conjugacy class of~$\phi_0$.We consider an affine embedding~$\iota:W\to \AAA^N_{\QQ}$ and the corresponding function~$H_f:W(\QQ_f)\to\ZZ_{\geq1}$ defined by~\eqref{notation:AF height}. We also consider an open compact subgroup~$K\leq G(\AAA_f)$. Then, as $\phi$ ranges through $W(\AAA_f)$, we have
\begin{equation}\label{Galois bounds}
[\phi(U):\phi(U) \cap K] \approx H_f(\phi).
\end{equation}
\end{enumerate}
\end{theorem}
Let us first reduce the second assertion to the first one.
\begin{proof}  
We identify~$G$ with its image by a faithful representation~$G\to GL(N)$. We may replace~$K$
by a commensurable group, and assume~$K$ is a  maximal compact subgroup of~$G(\AAA_f)$.
For any maximal compact subgroup~$K'$ of~$GL(N,\AAA_f)$ such that~$K\leq K'\leq GL(N,\AAA_f)$, we
then have
\begin{equation}\label{eq K Kf}
K=K'\cap G(\AAA_f).
\end{equation}
We choose such a~$K'$, and, possibly conjugating by an element of~$GL(N,\QQ)$, 
we may assume~$K'=GL(N,\widehat{\ZZ})$.

Consider~$\phi:M\to G$ in~\eqref{Galois bounds}.
From~$\phi(U)\leq G(\AAA_f)$ and~\eqref{eq K Kf}, we deduce
\begin{equation}\label{eq K Kf new}
[\phi(U):\phi(U) \cap K]=[\phi(U):\phi(U) \cap K']=[\phi(U):\phi(U) \cap GL(N,\widehat{\ZZ})].
\end{equation}
We have identified the left hand side of~\eqref{eq K Kf new} with the left hand side of~\eqref{Galois bounds}.

It will be enough to identify the right hand sides.
We will show that a Height function~$H_f$ on the $GL(N)$-conjugacy class of~$\phi$,
when restricted to the~$G$-conjugacy class, is a Height function on this~$G$-conjugacy class.

If~$H_f:GL(N,\AAA_f)\cdot \phi\to\ZZ_{\geq1}$ is associated to~$\iota:GL(N)\cdot \phi\to \AAA^N_\QQ$,
then its restriction to~$GL(N,\AAA_f)\cdot \phi$
is associated to~$\iota':G\cdot\phi\to GL(N)\cdot \phi\to \AAA^N_\QQ$, \emph{provided~$\iota'$ is a closed embedding}.

It is equivalent to proving that~$G\cdot \phi\subseteq GL(N,\AAA_f)\cdot \phi$ is Zariski closed.

To do this, we choose the map
\[
\iota:GL(N)\cdot \phi\xrightarrow{\phi'\mapsto d\phi'} \Hom(\mathfrak{m},\mathfrak{gl}(N)).
\]
Because~$G$ and~$GL(N)$ are Zariski connected, this map is injective. As~$M$ is reductive, according to~\cite[Th. 3.6]{Ri-Duke},
the image of~$G\cdot\phi$ is closed in~$\Hom(\mathfrak{m},\mathfrak{gl}(N))$, and thus~$G\cdot \phi\subseteq GL(N,\AAA_f)\cdot \phi$ is Zariski closed.
%Because height functions are polynomially equivalent, we may choose which height
%functions to use. We use~$H_f:GL(N,\AAA_f)\cdot \phi\to \ZZ_{\geq1}$ in\eqref{Galois bounds linear}; and~$H_f:G(\AAA_f)\cdot \phi\to \ZZ_{\geq1}$
%in~\eqref{Galois bounds}. 
%
%We note that~$G\cdot \phi_0$ is an affine closed subvariety in~$GL(N)\cdot \phi_0$. Indeed both can be viewed as closed affine subvarieties in the vector space~$\Hom(\mathfrak{m},\mathfrak{gl}(N))$. As a consequence we may use as height function on~$G\cdot \phi_0$ the restriction of the height function we use on~$GL(N)\cdot \phi_0$. 
%
%As a result we the inequality~\eqref{Galois bounds linear} is identical to the restriction of~\eqref{Galois bounds linear}
%to the~$\phi:M\to GL(N)$ such that~$\phi(M)\leq G$. It will be enough to prove~\eqref{Galois bounds linear}.
\end{proof}
We now reduce the first assertion to Cor.~\ref{B1}, Th.~\ref{general global bounds} and~\ref{general local bounds}.
\begin{proof}

Writing~$K=GL(N,\widehat{\ZZ})$ for short, we may rewrite the left hand-side of~\eqref{Galois bounds linear} as
\[
[\phi(U):\phi(U) \cap K]=\abs{\phi(U)\cdot K/K}=\abs{\left.U\bigl/\stackrel{-1}{\phi}(K)\right.}=
[U:\stackrel{-1}{\phi}(K)].
\]
Theorem~\ref{C1} implies~$[\phi(U):\phi(U) \cap K]\dom H_f(\phi)$. We now prove~$H_f(\phi)\dom[\phi(U):\phi(U) \cap K]$.

It is enough to obtain a lower bound after replacing~$U$ by the smaller group~$U'\leq U$ as in
Lem.~\ref{MT in product form}: without loss of generality we may assume~$U=U'$. We thus assume that~$U=\prod_p U_p$ as in ~\eqref{eq defi U' Up}.

The left hand-side is the product of the~$K_p=GL(N,\ZZ_p)$, hence we have
\[
[\phi(U):\phi(U) \cap K]
=\prod_p
[\phi(U_p):\phi(U_p) \cap K_p].
\]

We apply definition~\ref{defi MT in M} for~$K_M=M(\widehat{\ZZ})=M(\AAA_f)\cap GL(d,\widehat{\ZZ})$:
the upper bound~$C=\sup_p[M(\ZZ_p):U_p]$ is finite. Using~\eqref{explicit local adelic bound} we have 
\begin{equation}\label{explicit local adelic bound Up}
[\phi(M(U_p)):\phi(M(U_p))\cap GL(N,{\ZZ}_p)]\geq\frac{H_p(d\phi)}{c\cdot C}.
\end{equation}
As in the proof of~\eqref{explicit global adelic bound} of Theorem~\ref{general global bounds}, we can deduce
\begin{align}
[\phi(U):\phi(U)\cap GL(N,\widehat{\ZZ})]
\geq \frac{1}{(c\cdot C)^{\omega(H_f(d\phi))}}\cdot H_f(d\phi).
\end{align}
Arguing as in the proof of~\eqref{coro global bounds o} and~\eqref{coro global bounds dom} of Corollary~\ref{B1}, we obtain
\begin{equation}
H_{W,f}(\phi)\approx H_f(d\phi) \dom [\phi(U):\phi(U)\cap GL(N,\widehat{\ZZ})].
\end{equation}

\end{proof}

\subsection{Functoriality properties of the MT hypothesis}\label{sec:functor MT}
The following uses general properties of adelic topologies on algebraic groups. 
A good reference is~\cite{PR}.
\begin{lem}\label{lem MT properties}
Let~$\phi:M\to G$ be a morphism of connected linear algebraic groups over~$\QQ$,
and let~$U\leq M(\AAA_f)$ be a compact subgroup.
\begin{enumerate}
\item \label{MT properties 1}	
 If~$U$ is MT in~$M$ then~$\phi(U)$ is MT in~$\phi(M)$.
\item \label{MT properties 2}
If~$\phi$ is an isogeny onto its image (i.e.,~$\ker(\phi)$ is finite), then~$U$ 
is MT in~$M$ if~$\phi(U)$ is MT in~$\phi(M)$.
%\item Assume we have a direct product factorisation
%\[
%M=\prod_{i=1}^f M_i.
%\]
%Then~$U$ is MT in~$M$ if and only if~$U_i:=U\cap M_i(\AAA_f)$ is MT in~$M_i$. 
\item \label{MT properties 3}We assume~$M$ is reductive.
Let~$ad_M:M\to M^{ad}=M/Z_M(M)$ be the adjoint map, and~$ab_M:M\to M^{ab}=M/[M,M]$ be the abelianisation map. Then~$U$ is MT in~$M$ if and only if:~$ad_M(U)$ is MT in~$M^{ad}$ and~$ab_M(U)$ is MT in~$M^{ab}$.
\end{enumerate}
\end{lem}

The proof of Lemma~\ref{lem MT properties} will rely on the following.
\begin{theorem}\label{PR} Let~$\phi:H\to G$ be an epimorphism of $\QQ$-algebraic groups and~$C$ be the number of components of~$\ker(\phi)$
for the Zariski topology.
\begin{enumerate}
\item \label{PR fini}	
 Let~$K\leq H(\AAA_f)$ and~$K'\leq G(\AAA_f)$ be compact open subgroups of the form~$\prod_p K_p$ and~$\prod_p K'_p$. Then
\[
\forall p\gg0, \phi(K_p)\leq K'_p\text{ and }[K'_p:\phi(K_p)]\leq C.
\]
\item \label{PR connexe}
If~$C=1$ then the map~$p:H(\AAA_f)\to G(\AAA_f)$ is open: for any open subgroup~$K\leq H(\AAA_f)$, the image~$\phi(K)$
is open in~$G(\AAA_f)$.
\end{enumerate}
\end{theorem}
The second assertion, which is~\cite[p.\,296, \S6.2, Prop. 6.5]{PR},  is a corollary of the first assertion.
The first assertion follows from~\cite[p.,296, \S6.2, Prop. 6.4]{PR} and~\cite[p.\,296, \S6.2, Prop. 6.5]{PR} (using their exact sequence (6.9.)
under conditions of their Lemma 6.6.). 

Let us prove~\ref{MT properties 1} of Lemma~\ref{lem MT properties}.
\begin{proof} We choose a maximal compact subgroup~$K_M\leq M(\AAA_f)$, and a maximal compact subgroup~$K'\leq \phi(M)(\AAA_f)$ containing~$\phi(K_M)$.
By maximality, they have a product form~$K_M=\prod_p K_p$ and~$K'=\prod_p K'_p$. According to Definition~\ref{defi MT in M}, there exists~$c\in\RR_{>0}$
 such that for all primes~$p$, we have~$c\geq [K_p:U_p\cap K_p]$. Applying~$\phi$ we deduce
\[
c\geq [\phi(K_p): \phi(U_p\cap K_p)] \geq [\phi(K_p): \phi(U_p)\cap \phi(K_p)].
\] 
%On the other hand, because~$\phi(U_p)\leq \phi(U)_p$ and~$\phi(K_p)\leq K'_p$ in the first inequality, and~$\phi(U_p)\cap \phi(K_p)$ deduce
%\begin{multline}
%[K'_p:\phi(U)_p\cap K'_p]\leq [K'_p:\phi(U_p)\cap \phi(K_p)] \\\leq 
%[K'_p:\phi(K_p)]\cdot [\phi(K_p): \phi(U_p)\cap \phi(K_p)]
%\leq c .
%\end{multline}
Let~$C\in\RR_{>0}$ be given by Th.~\ref{PR}.
Using natural inclusions~$\phi(U_p)\subseteq \phi(U)_p$ and~$\phi(K_p)\subseteq K'_p$, we have
\begin{multline}
[K'_p:\phi(U)_p\cap K'_p]
%\leq[K'_p:\phi(U_p)\cap K'_p]
\leq [K'_p:\phi(U_p)\cap \phi(K_p)]
\\ = [K'_p:\phi(K_p)]\cdot [\phi(K_p):\phi(U_p)\cap \phi(K_p)]\leq C\cdot [\phi(K_p):\phi(U_p)\cap \phi(K_p)].
\end{multline}
Thus, for every prime~$p$, we have~$[K'_p:\phi(U)_p\cap K'_p]\leq c\cdot C$.
\end{proof}

We now prove~\ref{MT properties 2} of Lemma~\ref{lem MT properties}.
\begin{proof}We write~$K_M=\prod K_p$ and~$K'=\prod K'_p$ as before.

We choose a set of generators~$\phi(u_1),\ldots,\phi(u_k)$ for~$\phi(U)_p$ and let~$U'\leq U$ be the compact subgroup
topologically generated by the~$u_i$. Let us prove that~$k$ can be chosen independently of~$p$.
\begin{proof}
For a fixed~$p$ we use~\ref{lem gen 1} of~Lem.~\ref{lem generators} with~$V=\phi(U)_p$. For large~$p$, the group~$V':=\exp(p\phi(\mathfrak{m}_{\ZZ_p}))$ and 
the reduction map~$M(\ZZ_p)\to M(\FF_p)$ are well defined and, by~\ref{lem gen 3} of~Lem.~\ref{lem generators}, we have~$V'\leq V\leq M(\ZZ_p)$. Applying the remark
from the proof of Prop.~\ref{prop generators} to the exact sequence~$1\to V'\to V \to M(\FF_p)$, it follows from  Prop.~\ref{prop generators} for the image of~$V$ and~\ref{lem gen 2} of~Lem.~\ref{lem generators} for~$V'$.
\end{proof}
%Thanks to  Prop.~\ref{prop generators} and Lem.~\ref{}, we may bound~$k$ independantly of~$p$.

Let~$F$ be the kernel of~$\phi$. This is a finite algebraic group by hypothesis.
We define~$U'_p=U'\cap M(\QQ_p)$. Then~$U'_p$ is also the kernel of the map
\begin{multline}
U'\hookrightarrow \stackrel{-1}{\phi}\left(\phi(M(\AAA_f))\right)\to \stackrel{-1}{\phi}\left(\phi(M(\QQ_p))\right)\\
\to \stackrel{-1}{\phi}\left(\phi(M(\QQ_p))\right)/M(\QQ_p)\xleftarrow{\sim} F(\AAA_f)/(F\cap M)(\QQ_p). 
\end{multline}
The last group is a commutative group isomorphic to a subgroup of~$(\ZZ/(C))^\infty$ where~$C=\abs{F(\ol{\QQ})}$.
Because~$U'$ is generated by~$k$ elements, the size of the image of~$U'$ is bounded by~$C^{k}$. 

We deduce
\begin{multline}
[\phi(U)_p:\phi(U_p)]\leq [\phi(U'):\phi(U_p)]\\
\leq [\phi(U'):\phi(U'_p)]\leq [U':U'_p]\leq C^{k}.\qedhere
\end{multline}

%Assume~$\ker(\phi)$ is finite and let~$d$ be its number of geometric components. 
%
%For every~$p$, the group~$\ker(\phi)(\QQ_p)$
%is finite of size at most~$d$. We also have~$\ker(\phi)(\AAA_f)=\prod_p \ker(\phi(\QQ_p))$. Thus it is compact, and
%\(
%M(\AAA_F)\to\phi(M)(\AAA_f) 
%\)
%is proper. It follows that~$\stackrel{-1}{\phi}(K_{\phi(M)})$ is a compact subgroup of~$M(\AAA_f)$. As it contains the maximal compact subgroup~$K_M$,
%we have~$K_M=\stackrel{-1}{\phi}(K_{\phi(M)})$.
%...
\end{proof}

\begin{proposition} \label{prop generators}
For all~$N\in\ZZ_{\geq 0}$ there exists~$k=k(N)$ such that for every prime~$p$ and every subgroup~$U\leq GL(N,\FF_p)$, there exist~$u_1,\ldots,u_{k}$ in~$U$ which generate~$U$.
\end{proposition}

\begin{proof}We fix~$N$. There exists~$p(N)\in\ZZ_{\geq0}$ such that~$p\geq p(N)$, so that Nori applies~\cite{Nori}. 

For~$p\leq p(N)$ we have~$\#U\leq \#GL(N,\FF_p)\leq p(N)^{N^2}$ and we take~$k(N)=p(N)^{N^2}$.

We assume that~$p\geq p(N)$ and apply Nori theory~\cite{Nori}. 

According to Jordan theorem~\cite[Th. C]{Nori} there exist normal subgroups~$U^+\leq U'\leq U$ with~$U^+$ generated by the unipotent elements of~$U$, and~$U'/U^+$ abelian of order prime to~$p$, and~$[U:U']\leq d(N)$, where~$d(N)$ is as in~\cite[Th. C]{Nori}.  According to~\cite{Nori}, there exists~$\tilde{U}\leq GL(N)_{\FF_p}$ such that~$\tilde{U}(\FF_p)^+=U^+$. Define~$U''=\tilde{U}(\FF_p)\cap U$. Moreover, one knows\footnote{See~\cite[no 134, p. 25 and no 137,  p.\, 38--39, bottom of p. 44]{Serre--IV}.} that there exists an injective morphism~$U'/U''\hookrightarrow GL(N',\FF_p)$, where~$N'$ is bounded in terms of~$N$.

We will use the following remark. For every exact sequence~$1\to K\to G \to Q\to 1$, if~$K$ and~$Q$ are generated by~$k_N$ and~$k_Q$ elements, then~$G$ is generated by~$k_K+k_Q$ elements.
Thus, in order to bound the size of a generating subset of~$G$, it suffices to do it for~$K$ and for~$Q$.

Using the remark, it will be enough to prove that~$U/U'$, $U'/U''$, $U''/U^+$ and~$U^+$ can be generated by~$k_1(N),k_2(N),k_3(N),k_4(N)$ elements. Then the proposition will be satisfied with~$k(N)=\max\{k_1(N)+k_2(N)+k_3(N)+k_4(N);p(N)^{N^2}\}$.

As~$\#U/U' \leq d(N)$, we can take~$k_1(N)=d(N)$.

As~$\tilde{U}$ is generated by unipotent subgroups, we can write~$\tilde{U}=\tilde{S}\cdot \tilde{N}$ where~$\tilde{S}$ is semisimple and~$\tilde{N}$ is the unipotent radical. According to~\cite[Rem. 3.6, 3.6(v)]{Nori}, we have~$\tilde{S}(\FF_p)/\tilde{S}(\FF_p)^+\leq 2^N$. We deduce that~$\# U''/U^+\leq \# \tilde{U}(\FF_p)/\tilde{U}(\FF_p)^+=\# \tilde{S}(\FF_p)/\tilde{S}(\FF_p)^+\leq 2^N$.

We can thus take~$k_3(N)=2^N$.

The factor~$U'/U''$ is isomorphic to an abelian subgroup of~$GL(N',\FF_p)$ of order prime to~$p$. It is thus diagonalisable over some extension~$\FF_q$. Because~${\FF_q}^\times$ is cyclic, every subgroup of~$({\FF_q}^\times)^{N'}$ is generated by at most~$N'$ elements. 

We can thus take~$k_2(N)=N'$.

Let~$\tilde{U}\leq GL(N)_{\FF_p}$ be the algebraic group associated to~$U$ and let~$\tilde{\mathfrak{u}}\leq \mathfrak{gl}(N,\FF_p)$ be its Lie algebra. 
By \cite{Nori},~$\tilde{\mathfrak{u}}\leq \mathfrak{gl}(N,\FF_p)$ is linearly generated by nilpotents.
Let~$X_1,\ldots,X_d$, with~$d\leq N^2$ be a linear basis of nilpotent elements. Denote by~$U'=\langle\exp(X_1),\ldots,\exp(X_d)\rangle$  the group generated by their exponentials, and consider the associated~$\tilde{\mathfrak{u}}'\leq\tilde{\mathfrak{u}}$ and~$\tilde{U}'\leq \tilde{U}$. We have~$X_1,\ldots,X_d\in\tilde{\mathfrak{u}}'$. Thus~$\tilde{\mathfrak{u}}'=\tilde{\mathfrak{u}}$ and~$\tilde{U}'=\tilde{U}$. From~\cite[Th.~B]{Nori}, we get~$U=U^+=\tilde{U}(\FF_p)^+=\tilde{U}'(\FF_p)^+=U'^+=U'$. Thus~$U$ is generated by at most~$N^2$ elements~$\exp(X_1),\ldots,\exp(X_d)$.

We can thus take~$k_4(N)=N^2$.
%By \cite{Nori},~$U$ is generated by~$\exp(\mathfrak{u})$, for a Lie algebra~$\mathfrak{u}\leq \mathfrak{gl}(N,\FF_p)$ which is linearly generated by nilpotents, and of dimension at most~$N^2$. Let~$X_1,\ldots,X_{N^2}$ be a nilpotent basis of~$\mathfrak{u}$. Let~$U'\leq U$ be the subgroup generated by~$\exp(X_1),\ldots,\exp(X_{N^2})$, and~$\mathfrak{u}'\leq \mathfrak{u}$ be its Lie algebra. Then~$X_1,\ldots,X_{N^2}\in\mathfrak{u}'$ and thus~$\mathfrak{u}'=\mathfrak{u}$.
% and , and we take the exponentials of a basis of nilpotent generators.
\end{proof}

We used the following.
\begin{lem}\label{lem generators}
\begin{enumerate}Let~$M\leq GL(N)$ be an algebraic subgroup defined over~$\QQ_p$ and let~$\mathfrak{m}\leq \mathfrak{gl}(N,\QQ_p)$ be its Lie algebra.
\item \label{lem gen 1} Let~$V\leq GL(N,\ZZ_p)$ a compact subgroup. Then~$V$ is topologically finitely generated.
\item \label{lem gen 2}
Then~$V':=\exp\bigl(\mathfrak{m}\cap  2p\mathfrak{gl}(N,\ZZ_p)\bigr)$
is topologically generated by at most~$N^2$ elements if~$p$ is large enough.
\item \label{lem gen 3}
Let~$M(\ZZ_p):=M(\QQ_p)\cap GL(N,\ZZ_p)$ and~$V\leq M(\ZZ_p)$ an open subgroup such that~$C:=[M(\ZZ_p):V]\in\ZZ_{\geq1}$.
Then for~$p>C$, we have
\[V'\leq V.\]
% see PR page 139
\end{enumerate}
\end{lem}
\begin{proof}
The first assertion is~\cite[Prop. 2]{Serre-62}.

Let~$G=\exp(2p\mathfrak{gl}(N,\ZZ_p))= 1+2p\mathfrak{gl}(N,\ZZ_p)$ and~$H=V'\leq G$. According to~\cite[Th.~5.2]{dS} the pro-$p$ group
is powerful and~$d(G)=N^2$, where~$d(G)$ is
the minimal cardinality of a set of generators for $G$ as in~\cite[p.\,41]{dS}. We can thus apply~\cite[Th.~2.9]{dS}. This proves the second assertion.

As~$G$ is a pro-$p$ group, by~\cite[Lem 4.8, p.\,138]{PR},~$V'$ is a pro-$p$ group. We also have
\[
[V':V'\cap V]\leq [M(\ZZ_p):V]=C.
\]
Assume~$p>C$ and 
assume by contradiction that there exists~$w\in V'\smallsetminus V$. We denote by~$w^\ZZ$ the subgroup generated by~$w$. Then~$c:=[w^\ZZ:w^\ZZ\cap V]\neq 1$ and~$c\leq C$.
But~$c$ is a power of~$p$ because~$V'$ is a pro-$p$ group: thus~$c\geq p$. We deduce that~$
C\geq c\geq p$. This contradicts our assumptions.
\end{proof}

We prove~\ref{MT properties 3} of Lemma~\ref{lem MT properties}.
We will make use of Goursat's Lemma.
\begin{proof}As~$M$ is reductive, the map~$(ad_M,ab_M):M\to M':=M^{ad}\times M^{ab}$ is an isogeny. From~\eqref{MT properties 2} of Lemma~\ref{lem MT properties} it 
follows that it is enough to prove that the image~$V$ of~$U$ in~$M'(\AAA_f)$
is MT in~$M'$. We may thus assume~$M=M^{ad}\times M^{ab}$.

Using Lemma~\ref{MT in product form} we may assume~$U=\prod_p U_p$. 
Let~
\[
K_M=\prod K_{M,p}=\prod_p K_{M^{ad},p}\times K_{M^{ab},p}\leq M(\AAA_f)
\]
be a maximal compact subgroup containing~$U$.

By assumption there is an upper bound~$C\in\ZZ_{\geq1}$ for~$[K_{M^{ab},p}:ab_M(U_p)]$ and~$[K_{M^{ad},p}:ad_M(U_p)]$,
independent of~$p$.

Let~$H_1=ad_M(U_p)$ and~$H_2=ab_M(U_p)$ and~$\Gamma=(ad_M,ab_M)(U_p)\leq H_1\times H_2$.
Let~$N_1=\Gamma\cap H_1$ and~$N_2=\Gamma\cap H_2$. By Goursat's Lemma,~$N_1$ and~$N_2$ are normal subgroups in~$H_1$ and~$H_2$
and there is an isomorphism (whose graph is~$\Gamma/(N_1\times N_2)$)
\begin{equation}\label{Goursat isomorphism}
H_1/N_1\xrightarrow{\sim} H_2/N_2.
\end{equation}
Because~$H_2$ is abelian,~$N_1$ contains the derived subgroup~$[H_1,H_1]$. 

By the first part of Lemma~\ref{Lemma commutator subgroup},~$[H_1:N_1]$ is finite for every prime~$p$, 
and by the second part of Lemma~\ref{Lemma commutator subgroup},~$[H_1:N_1]$ is bounded by~$C(M^{ad})$ for almost every prime~$p$.

As a result there exists~$C'\in\ZZ_{\geq1}$ such that~$[H_1:N_1]\leq C'$ for every prime~$p$.
Using~\eqref{Goursat isomorphism}, we also have~$[H_2:N_2]\leq C'$ for every prime~$p$.

Recall that~$N_1\times N_2\leq \Gamma$. 
It follows
\[
[H_1\times H_2:\Gamma]\leq[H_1:N_1]\cdot[H_2:N_2]=C'^2.
\]
By definition of~$C$, 
\[
[K_{M,p}:H_1\times H_2]\leq C^2.
\]
We deduce
\[
[K_{M,p}:(ad_M,ab_M)(U_p)]= [K_p:H_1\times H_2]\cdot [H_1\times H_2:\Gamma]\leq C^2C'^2.
\]
The bound is independent of~$p$, which concludes.
\end{proof}
\begin{lem}\label{Lemma commutator subgroup}
 Let~$G$ be a semisimple algebraic group over~$\QQ$, and for every prime~$p$,
let~$U_p,K_p\leq G(\QQ_p)$ be compact open subgroups such that~$K=\prod_p K_p\leq G(\AAA_f)$ is open.
Let~$[U_p,U_p]$ be the subgroup generated by commutators. 
%We recall that~$G(\ZZ_p)$ is well defined for almost all~$p$, associated for example to a linear representation.
\begin{enumerate}
\item For every prime~$p$, the quotient~$U_p/[U_p,U_p]$ is finite.
\item There exists~$C(G)\in\ZZ_{\geq1}$ such that, for almost all~$p$, if~$[K_p:U_p\cap K_p]<p$ then~$U_p/[U_p,U_p]<C(G)$.
\end{enumerate}
\end{lem}
%\begin{proof}
%The first statement follows from the fact that~$[U_p,U_p]$ is open, because~$G$
%is semisimple.
%
%Let us identify~$G$ with its image by a faithful linear representation~$G\to GL(N)$.
%For~$p$ large enough, we have~$K_p\leq G(\ZZ_p):=G(\QQ_p)\cap GL(N,\ZZ_p)$.
%Let~$K_p^1$ be the kernel of the reduction map~$G(\ZZ_p)\to G(\FF_p)\leq GL(N,\FF_p)$.
%
%For~$g\in K^1_p$, the series~$X:=\log(g)$ converges, and we have~$g=\exp(X)$.
%(Prop.~\ref{Proposition bornes log matrice}).
%We have
%\[
%[g^{\ZZ}:g^{\ZZ}\cap U_p]\leq [K_p:K_p\cap U_p]<p.
%\]
%But~$[g^{\ZZ}:g^{\ZZ}\cap U_p]$ is a power of~$p$. Thus~$[g^{\ZZ}:g^{\ZZ}\cap U_p]=1$,
%and~$g\in U_p$.
%
%
%Let~$p$ be large enough, so that the subgroup~$K_p$ is hyperspecial.
%We have~$K_p=G(\ZZ_p)$, and~$G(\FF_p)$ is the image of~$K_p$.
%Let~$U(p)$ be the image of~$U_p$ in~$G(\FF_p)$.
%Then the derived subgroup~$[U(p),U(p)]$ is the image of~$[U_p,U_p]$.
%
%We have~$[G(\FF_p):G(\FF_p)\cap U(p)]< p$. Thus, using Nori theory~\cite{Nori}, we have
%\[
%G(\FF_p)^+\leq U(p).
%\]
%
%Let~$\ol{H}=[U(p),U(p)]$
%
%
%Let us first prove that~$[G(\FF_p):G(\FF_p)\cap U(p)]$
%is bounded. 
%\begin{proof}
%We have~$[G(\FF_p):G(\FF_p)\cap U(p)]< p$. Thus, using Nori theory~\cite{Nori}, we have
%\[
%G(\FF_p)^+\leq U(p),
%\]
%and~$[G(\FF_p):G(\FF_p)\cap U(p)]\leq [G(\FF_p):G(\FF_p)^+]\leq 2^N$. (\cite[p. 270]{N}) 
%\end{proof}
%Let~$[U,U]\leq [G(\FF$. 
%\end{proof}
\begin{proof} The first assertion follows from the fact that~$[U_p,U_p]$ is open, because~$G$
is semisimple.

We prove the second assertion. We may replace~$U_p$ by~$K_p\cap U_p$ and  assume~$U_p=K_p\cap U_p\leq K_p$. Thus~$[K_p:K_p\cap U_p]=[K_p:U_p]<p$.

Let us identify~$G$ with its image by a faithful linear representation~$G\to GL(N)$.
For~$p$ large enough, we have~$K_p= G(\ZZ_p):=G(\QQ_p)\cap GL(N,\ZZ_p)$. 

Let~$G(\ZZ_p)^+$ and~$G(\FF_p)^+$ be as in Lemma~\ref{CK Fact} below.

Then~$U_p\cap G(\ZZ_p)^+$ is an open subgroup of $G(\ZZ_p)^+$ and, 
\[
[G(\ZZ_p)^+:U_p\cap G(\ZZ_p)^+]\leq [G(\ZZ_p):U_p]<p.
\]
(Recall the assumption~$[K_p:U_p]<p$.)

As~$G(\ZZ_p)^+$ is generated by pro-$p$-groups,
we have, for every subgroup~$L\leq G(\ZZ_p)^+$,
\[
[G(\ZZ_p)^+:L]>1\Rightarrow
[G(\ZZ_p)^+:L]\geq p.
\]
Therefore, with~$L=U_p$, 
\[
[G(\ZZ_p)^+:U_p\cap G(\ZZ_p)^+]=1.
\]
At the level of derived subgroups, we have
\[
[G(\ZZ_p)^+,G(\ZZ_p)^+]\subseteq [U_p,U_p].
\]
We deduce
\[
[G(\ZZ_p):[U_p,U_p]]\leq
[G(\ZZ_p):G(\ZZ_p)^+]
\cdot
 [G(\ZZ_p)^+:[G(\ZZ_p)^+,G(\ZZ_p)^+]].
\]
We note that~$G(\ZZ_p)^+\leq G(\ZZ_p)$ is an open subgroup of index prime to~$p$. It follows that the image of~$G(\ZZ_p)^+$ in~$G(\FF_p)$ contains~$G(\FF_p)^+$. Thus
\[
[G(\ZZ_p):G(\ZZ_p)^+]\leq [G(\FF_p):G(\FF_p)^+].
\]
We have,~ by~\cite[p. 270]{Nori},
\begin{equation}\label{CK1}
[G(\FF_p):G(\FF_p)^+]\leq 2^N.
\end{equation}

For~$p$ large enough:
\begin{itemize}
\item We have~$G(\FF_p)=\widetilde{G}(\FF_p)$ for a connected  semisimple~$\FF_p$-algebraic subgroup~$\widetilde{G}\leq GL(N)_{\FF_p}$;
\item We have~$[\widetilde{G}(\FF_p)^+,\widetilde{G}(\FF_p)^+]=[\widetilde{G},\widetilde{G}](\FF_p)^+=\widetilde{G}(\FF_p)^+$, using Lem.~\ref{lem:parfait}.
\end{itemize}
Thus~$[G(\ZZ_p)^+,G(\ZZ_p)^+]$ maps surjectively onto
\[
[G(\FF_p)^+,G(\FF_p)^+]=G(\FF_p)^+.
\]

We apply Lemma~\ref{CK Fact} to~$H=[G(\ZZ_p)^+,G(\ZZ_p)^+]$. We deduce
\[
[G(\ZZ_p)^+,G(\ZZ_p)^+]=G(\ZZ_p)^+.
\]
This implies
\begin{equation}\label{CK2}
[G(\ZZ_p)^+:
[G(\ZZ_p)^+,G(\ZZ_p)^+]]=1.
\end{equation}

The second assertion of Lemma~\ref{Lemma commutator subgroup} follows from~\eqref{CK1} and~\eqref{CK2}.
\end{proof}
\begin{lem}[{\cite[Fact 2.4 and its proof]{CK}}]\label{CK Fact}
Let~$G\leq GL(N)_{\QQ}$ be a connected semisimple algebraic subgroup. For every prime~$p$, define~$G(\ZZ_p):=G(\QQ_p)\cap GL(N,\ZZ_p)$ and denote by~$G(\FF_p)$ the image of~$G(\ZZ_p)$ in~$GL(N,\FF_p)$. We denote by~$G(\FF_p)^+\leq G(\FF_p)$ and~$G(\ZZ_p)^+\leq G(\ZZ_p)$ the subgroups generated by $p$-Sylow subgroups, resp. pro-$p$-Sylow.

Then, for~$p$ large enough: if~$H\leq G(\ZZ_p)^+$ maps surjectively onto~$G(\FF_p)^+$, then~$H=G(\ZZ_p)$.
\end{lem}
%
%
%For~$p$ large enough, the group~$K_p$ is hyperspecial and we have~$K_p=G(\ZZ_p)$
%with respect to a smooth reductive model of~$G$ over~$\ZZ_p$ 
%We prove the second statement. Because~$K$ is open, we have~$U_p\leq K_p$ for almost all primes.
%For any unipotent subgroup~$N\leq G_{\QQ_p}$ defined over~$\QQ_p$, the order 
%of the quotient~$N\cap K_p/N\cap U_p$ is bounded by~$[K_p:U_p\cap K_p]<p$. 
%This order is also a power of~$p$ because~$N$ is unipotent. This order is then~$1$.
%Equivalently,~$N\cap K_p = N\cap U_p$.
%
%Let~$K^+_p$ be the group generated by the~$N\cap K_p$ as~$N$ ranges through unipotent 
%subgroups~$N\leq G_{\QQ_p}$.  It is enough to bound~$[K_p:K^+_p]$
%for almost all primes. Using ... we may assume~$G$ is simply connected.
%For almost primes,~$G_{\QQ_p}$ is quasi-split, and~$K_p$ is hyperspecial.
%\end{proof}
We used the following in the proof of~Lemma~\ref{Lemma commutator subgroup}.
\begin{lem}\label{lem:parfait}
	For every~$n\in\ZZ_{\geq0}$, there exists~$c(n)$ such that the following holds.
 Let~$p\geq c(n)$ be a prime, and let~$G\leq GL(n)$ be a semisimple algebraic group over~$\FF_p$.
 
Then~$[G(\FF_p)^+,G(\FF_p)^+]=G(\FF_p)^+$.% where~$[
\end{lem}
\begin{proof}
Let~$\pi:G^{sc}\to G$ be the simply connected cover. According to \cite[24.15]{MaTe}, we have~$G^{sc}(\FF_p)^+=G^{sc}(\FF_p)$.

It follows that~$\pi(G^{sc}(\FF_p))\leq G(\FF_p)^+$. Since~$G(\FF_p)^+$ is generated by elements of order a power of~$p$,
we have the following alternative
\begin{itemize}
\item either~$\pi(G^{sc}(\FF_p))= G(\FF_p)^+$
\item or~$\# G(\FF_p)^+/\pi(G^{sc}(\FF_p))\geq p$.
\end{itemize}

Let~$Z=\ker \left(\pi: G^{sc}(\ol{\FF_p})\to G(\ol{\FF_p})\right)$.

By~\cite[24.21]{MaTe}, we have~$\# G(\FF_p)^+/\pi(G^{sc}(\FF_p))\leq \#Z$.

On the other hand, there exists an integer $c(n)$ (depending only on $n$) such that
we have~$\# Z\leq c(n)$. Thus, for~$p>c(n)$, the second
case of the alternative does not happen.

By~\cite[24.17]{MaTe}, for~$p\geq 4$, the group~$G^{sc}(\ol{\FF_p})$ is perfect. It implies
that its quotient~$\pi(G^{sc}(\FF_p))= G(\FF_p)^+$ is perfect, namely, that~$[G(\FF_p)^+,G(\FF_p)^+]=G(\FF_p)^+$.
\end{proof}

\subsection{MT hypothesis for Images of Galois representations}
We use the notations of Def.~\ref{MT hypothesis}. We assume furthermore that~$E$ is of finite type over~$\QQ$.
In this case we have the following.
\begin{lem}\label{lem mt si special}
If~$x_0$ is a special point (i.e.,~$M=M^{ab}$), then~$x_0$ satisfies the MT hypothesis.
\end{lem}
The Galois representation~$Gal(\ol{E}/E)\to M^{ab}(\AAA_f)$ is prescribed by Deligne-Shimura reciprocity law, which is part of the
definition of a Canonical model~\cite[2.2.5]{Del79}. In this case, we know that~$M=M^{ab}$ is the 
Zariski closure of the image of~$x_0$. It follows that the morphism~\cite[2.2.2.1]{Del79} is 
an epimorphism, and we can apply Th.~\ref{PR}.\footnote{If the kernel of~$\mu_h:GL(1)_E\to T_E$ is connected,
then the Galois image is actually open for the topology induced by the adelic topology on~$T(\AAA_f)$. 
This is also the~$H$-maximality condition. See~\cite{CM}.}

%The following reduction applies only to Galois images and relies on Shimura-Deligne reciprocity law.

Using Lemma~\ref{lem mt si special} and \eqref{MT properties 3} of Lemma~\ref{lem MT properties} we have the following.\begin{lem}\label{MT product form}
The point~$x_0$ satisfies the MT hypothesis if and only if~$ad_M(x_0)$
satisfies the MT hypothesis. 
\end{lem}\label{lemme special est MT}

The following is not needed but can help relate our MT hypothesis to other notions found in literature.
\begin{theorem}\label{thm:not needed} Assume~$M$ is a semisimple and simply connected algebraic group over~$\QQ$.
Then a compact subgroup~$U\leq M(\AAA_f)$ is MT in~$M$ if and only if it is an open subgroup.
\end{theorem}
Theorem~\ref{thm:not needed} is a consequence of the following.
\begin{lem} Let~$M\leq GL(n)_{\QQ}$ be a simply connected semisimple $\QQ$-algebraic subroup,
and, for every prime~$p$, define~$M(\ZZ_p):=M(\QQ_p)\cap GL(n,\ZZ_p)$.

There exists~$C$ such that for every prime~$p\geq C$, every~$U\leq M(\ZZ_p)$ satisfies
\[
U=M(\ZZ_p)\text{ or }[M(\ZZ_p):U]\geq p.
\]
\end{lem}
\begin{proof}This is a consequence of the following claim: for~$p\gg0$, the group~$M(\ZZ_p)$ is generated by topologically~$p$-nilpotent elements.

Let us prove the claim. For every prime~$p$, every element in the kernel of the morphisms~$red_p:M(\ZZ_p)\to GL(n,\FF_p)$ belongs to~$1+p\mathfrak{gl(n,\ZZ_p)}$  and is topologically~$p$-nilpotent. It will be enough to prove that~$red_p(M(\ZZ_p))$ is generated by elements of order power of~$p$. 
For~$p\gg 0$, the group~$M(\ZZ_p)$ is hyperspecial and the model of~$M$ induced by~$GL(n)_{\ZZ}$ is smooth over~$\ZZ_p$ with semisimple fibre~$M_{\FF_p}$. This implies that the map~$M(\ZZ_p)\to M_{\FF_p}(\FF_p) $ is surjective.
 For~$p\gg0$ the algebraic group~$M_{\FF_p}$ is semisimple and simply connected\footnote{This is~\cite[\S6.5]{MVW} and here is an argument. Passing to a finite extension of~$\QQ$ we may assume that~$M$ is simply connected and hyperspecial at~$p$. Applying~\cite[\S3.5.4]{Tits}, we deduce that the special fibre is simply connected. The case~$C-BC_n$ of~\cite[p.\,61]{Tits} is excluded in the hyperpecial case.}. By~\cite[24.15]{MaTe} we have~$M(\FF_p)=M(\FF_p)^+$. This proves the claim.
\end{proof}

This relies on strong approximation, Hasse principle and Kneser-Tits properties for~$M$. See~\cite{Del71} for related discussions.
\subsubsection{}\label{Cadoret abelian}
For moduli spaces of abelian varieties or more generally for Shimura varieties of abelian type, a Galois representation associated to a point~$x_0\in X$ can be deduced from the Galois representation on the Tate module of an abelian variety.
%We relate here our MT hypothesis with other properties dealing 

We have the following.
\begin{theorem}[{\cite[Th.~A(i)]{CM} \cite[Th.~10.1]{HR}}]
Let~$S$ be a Shimura variety of Hodge type, let~$s\in S$ be a point.

If the abelian variety~$A$ associated to~$s$ satisfies the classical Mumford-Tate conjecture at some prime~$\ell$,
then~$s$ satisfies the weakly adelic Mumford-Tate hypothesis.
\end{theorem}
Using Lemma~\ref{MT product form} we can deduce
\begin{theorem}
Let~$S$ be a Shimura variety of abelian type, let~$s\in S$ be a point.

If~$s$ satisfies the Mumford-Tate conjecture at some prime~$\ell$ in the sense of~\cite{UY_General},
then~$s$ satisfies the weakly adelic Mumford-Tate hypothesis.
\end{theorem}

\subsubsection{}\label{Greg}
As observed in~\cite{G}, the combination of a theorem of Deligne-André and with a theorem of Weisfeiler \cite{MVW} and Nori \cite{Nori} produces, in \emph{any} Shimura variety, \emph{many} examples of (non algebraic) points for which the MT hypothesis is satisfied. With our terminology it is stated as follows.
\begin{theorem}[{\cite[Th\,1.2]{G}}]\label{Baldi MT}
Let~$M$ be the Mumford-Tate group of a point~$x_0\in X$ for a Shimura datum~$(G,X)$. We decompose the adjoint datum~$(M^{ad},X_{M^{ad}})$ of~$(M,X_M):=(M,M(\RR)\cdot x_0)$ as a product
\[
(p_1,\dots,p_f):(M^{ad},X_{M^{ad}})\simeq (M_1,X_1)\times \ldots \times (M_f,X_f)
\]
with respect to the~$\QQ$-simple factors~$M_i$ of~$M^{ad}$. 

Assume that for some compact open subgroups~$K_i\leq M_i(\AAA_f)$
\[
\forall i\in\{1;\ldots;f\},\,[p_i\circ ad_M (x_0)]\in Sh_{K_i}(M_i,X_i)(\CC)\smallsetminus Sh_{K_i}(M_i,X_i)(\ol{\QQ}).
\]
Then~$x_0$ satisfies the MT hypothesis.
\end{theorem}

\section{Proof of the main result}\label{Main proof}
In this section we prove the  Theorem~\ref{thm intro}, following the strategy outlined in~\S\ref{Strategy outline}.
We then give in~\S\ref{PW} a variant of Pila-Wilkie Theorem.

%, which relies on a combination of the Pila-Wilkie theorem~\cite{Pila},
%the definability of uniformisation maps of Shimura varieties restricted to Siegel sets~\cite{KUY}, functional transcendence 
%from the Hyperbolic Ax-Lindemann-Weierstrass theorem~\cite{KUY} and the Geometric part of André-Oort conjecture.
%The latter was achieved through the Pila-Zannier method~\cite{UllmoAxLin}, or more recently via dynamics based
%on Ratner's theorems~\cite{RU}. 
\subsection{Reduction steps}% Proof of the Generalised André-Pink Zannier conjecture under Mumford-Tate conjecture}
We put ourselves in the situation of Theorem~\ref{thm intro} and Conjecture~\ref{APZ}.

Let~$Z$ be an irreducible component of~$\ol{\Sigma}^{Zar}$. The aim is to prove that~$Z$ is weakly special.
We may replace~$\Sigma$ by~$\cH(x_0)\cap Z$.

\subsubsection{Reduction to the Hodge generic case}
We will reduce the theorem to the case where~$Z$ is Hodge generic in~$Sh_K(G,X)$. For convenience we will assume that~$s_0=[x_0,1]\in Z$.
We choose a Hodge generic point~$z$ in~$Z$. One knows that one can choose a lift~$\wt{z}$ of~$z$ in~$X$ such that the Mumford-Tate group~$G'$ of~$\wt{z}$ contains~$M$. We write~$X'=G'(\RR)\cdot \wt{z}$. We have a Shimura morphism
\[
\Psi:Sh_{K\cap G'(\AAA_f)}(G',X')\to Sh_K(G,X).
\]
% morphism of algebraic varieties.
(The smallest special subvariety of~$Sh_K(G,X)$ containing~$Z$ is the image of one component of~$Sh_{K\cap G'(\AAA_f)}(G',X')$.) Let~$Z'$ be the inverse image of~$Z$ by~$\Psi$. It is known that~$Z$ is weakly special if and only if any component of~$Z'$ is weakly special. 

In the notations of Prop.~\ref{prop passage to Hodge generic}, %defining
%\[
%\Sigma:=\cH([x_0,1])\cap Z\qquad
%\Sigma':=\cH'([x_0,1])\cap Z',
%\]
we have
\[
\Sigma':=\stackrel{-1}{\Psi}(\Sigma)=\cH'([x_0,1])\cap Z'.
\]
Because~$Sh_{K\cap G'(\AAA_f)}(G',X')\to \Psi\bigl(Sh_{K\cap G'(\AAA_f)}(G',X')\bigr)$ is flat, and because~$Z$ is in the image of~$\Psi$, we deduce that~$\Sigma'$ is dense in~$Z'$,
and hence dense in every component of~$Z'$.

Thus, in proving the conclusion of the theorem we may replace~$Z$ by a component of~$Z'$, and~$(G,X)$ by~$(G',X')$, and~$K$ by~$K\cap G'(\AAA_f)$.

On the other hand, the Mumford-Tate hypothesis depends only on~$M$, and thus is insensitive to such substitutions.

In other words, \emph{we can, and will, assume that~$Z$ is Hodge generic in~$Sh_K(G,X)$.}

\subsubsection{Reduction to the adjoint datum}
We will reduce the theorem to the case where~$G=G^{ad}$ is of adjoint type. Here we use geometric Hecke orbits.

Using Theorem~\ref{theo general union finie geometrique}, we write
our generalised Hecke orbit
\[
\cH([x_0,1])=
\cH^g([x_0,1])\cup\ldots\cup \cH^g([x_k,1])
\]
as a finite union of geometric Hecke orbits. We define accordingly
\[
\Sigma_i=Z\cap \cH^g([x_i,1]).
\]
As~$Z$ is irreducible there at least one~$i\in\{0;\ldots;k\}$
such that~$\Sigma_i$ is Zariski dense in~$Z$.

Because the Galois representations~$\rho_{x_1},\ldots,\rho_{x_k}$ of~$x_1,\ldots,x_k$ can be deduced from~$\rho_{x_0}$ using \S~\ref{section functoriality}, the Mumford-Tate hypothesis will still be valid even if we replace~$x_0$ by~$x_i$. We assume for simplicity that~$x_i=x_0$.

We choose an open compact subgroup~$K'\leq G^{ad}(\AAA_f)$
so that we can consider the Shimura morphism
\[
\Psi:Sh_K(G,X)\to Sh_{K'}(G^{ad},X^{ad}).
\] 
Let~$Z'$ be the image of~$Z$. One knows that~$Z$ is weakly special in~$Sh_{K}(G,X)$ if and only if~$Z'$ is weakly special in~$Sh_{K'}(G^{ad},X^{ad})$.

Then~$\Psi(\Sigma_0)$ is dense in~$Z'$. Denote~$x_0^{ad}$ the image of~$x_0$ in~$X^{ad}$, and define
\[
\Sigma':=\cH^g([x_0^{ad},1]).
\]
Using \S~\ref{subsection passing to adjoint}, we get
\[
\Psi(\Sigma_0)\subset\Sigma'\subset Z'
\]
and thus~$\Sigma'$ is Zariski dense in~$Z'$.

Let~$M'$ be the image of~$M$ by~$ad_G:G\to G^{ad}$.
Then~$M^{ad}\simeq M'^{ad}$ because~$\ker(ad_G)$ is
commutative and central in~$G$. In view of \S~\ref{section MT adelic}, the Mumford-Tate hypothesis will still hold for~$x_0^{ad}\in X^{ad}$.

Thus, \emph{we can, and will, assume~$G=G^{ad}$.}

\subsubsection{Induction argument for factorisable subvarieties}\label{partition factorisation}
The following reduction will be useful at the very end of the whole proof.

We recall that~$G$ is a direct product~$G_1\times \ldots \times G_f$ of its~$\QQ$-simple subgroups. 

It can be easily proved that in the Theorem~\ref{thm intro} we can replace~$K$ by any other compact open subgroup. After possibly replacing~$K$ by the open subgroup~$\prod_{i=1}^{f}K_i:=\prod_{i=1}^{f} K\cap G_i(\AAA_f)$, there are factorisations
\(
X=\prod_{i=1}^{f} X_i\) and 
\begin{equation}\label{facto proof}
Sh_K(G,X)=\prod_{i=1}^{f} Sh_{K_i}(G_i,X_i).
\end{equation}
The factorisation~\eqref{facto proof} is defined over the reflex field~$E(G,X)$,
hence over~$E$.
Consider a nontrivial partition~$\{1;\ldots;f\}=I\sqcup J$
and the corresponding nontrivial factorisation of Shimura data
\[
(G,X)\xrightarrow{(p_I,p_J)	} (G_I,X_I)\times (G_J,X_J)
\]
with
\[
(G_I,X_I)=\prod_{i\in I}(G_i,X_i)\text{ and }
(G_J,X_J)=\prod_{i\in I}(G_j,X_j).
\]
By functoriality~\S\ref{subsection functoriality} for~$\phi=p_I\circ \phi_0$ (resp.~$\phi=p_J\circ \phi_0$), we will have
\[
\rho_{p_I(x_0)}=p_I\circ \rho_{x_0}\text{ and }
\rho_{p_J(x_0)}=p_J\circ \rho_{x_0}.
\]
As explained in \S~\ref{section MT adelic}, the Mumford-Tate hypothesis will hold for~$p_I(x_0)$ and for~$p_J(x_0)$.

Suppose that~$Z$ factors as a Cartesian product
\begin{equation}\label{product form}
Z_I\times Z_J \subseteq Sh_{K_I}(G_I,X_I)\times Sh_{K_J}(G_J,X_J)
\end{equation}
in the corresponding factorisation of Shimura varieties. From~\S\ref{section compatibility products}, we have
\[
\cH^g(x_0)=\cH^g(p_I(x_0))\times \cH^g(p_J(x_0))
\]
and
\[
\cH^g([x_0,1])=\cH^g([p_I(x_0,1)])\times \cH^g([p_J(x_0),1]).
\]
Recall that the partition~$\{1;\ldots;f\}=I\sqcup J$ is not trivial. Arguing by induction on~$f$, 
we can assume that the Theorem~\ref{thm intro} proven for~$Z_I$ and~$Z_J$.
Then~$Z_I\times Z_J$ is also a weakly special subvariety and we are done.

\emph{Henceforth we assume that for every non trivial partition~$\{1;\ldots;f\}=I\sqcup J$, the variety~$Z$ is not a product of the form~\eqref{product form}.}

%NOTE :
%
%FAIRE: PASSAGE À HODGE GENERIQUE (ORBITES GEN). PUIS
%PASSAGE À ORBITE GEOMETRIQUE. PUIS PASSAGE À ADJOINT.
%
%As seen in section~\ref{} we may assume that~$G$ is adjoint. 
%
%We may assume that~$Z$ is Hodge generic in~$Sh_K(G,X)$, possibly passing to the Shimura
%datum~$(G',G'(\RR)\cdot x_Z)$ where~$[x_Z,1]$ is a Hodge generic point of~$Z$. 
\subsection{Central arguments.}
%\subsubsection{}	
Let us recollect some of the notations and notions we will be using.

We have an irreducible subvariety~$Z$ of $\Sh_K(G,X)$ containing a Zariski dense subset~$\Sigma$ 
contained in the generalised Hecke orbit $\cH([x_0,1])$ of the point $[x_0,1]$. Let $E$ be a field of finite type over $\QQ$ 
such that $Z$ and $[x_0,1]$ are defined over $E$, and passing to a finite extension we have a Galois
representation~$\rho:\Gal(\ol{E}/E)\to M(\AAA_f)\cap K$ as in Def.~\ref{terminologie galois repr} and our main hypothesis
is that its image~$U:=\rho(\Gal(\ol{E}/E))$ satisfies Def.~\ref{MT hypothesis}. 
Passing to a finite extension we also assume that~$E$ is a field of definition for every geometric component of~$Sh_K(G,X)$. 
%We choose a basis 
%of the Lie algebra~$\mathfrak{g}$ so, with respect to that the associated~$GL(\mathfrak{g})\to GL(d)$,
%integral  structures are meaningful and we have a~$Gal(\ol{E}/E)$-invariant height function~$\cH^g(x_0)$.

We reduce the Theorem~\ref{thm intro} to the case where~$\Sigma$ is contained in a single geometric Hecke orbit. According to Th.~\ref{theo general union finie geometrique} the generalised Hecke orbit is a 
finite union of geometric orbit, with~$\phi_0:M\to G$ the identity map, 
\begin{equation}\label{gen into geo}
\cH([x_0,1])=\cH^g([x_0,1])\cup\cH^g([\phi_1\circ x_0,1])\cup\ldots\cup\cH^g([\phi_k\circ x_0,1]).
\end{equation}
As~$Z$ is irreducible, at least one of the intersections~$Z\cap \cH^g([\phi_i\circ x_0,1])$
is Zariski dense in~$Z$. From~\S\ref{subsection functoriality}, we obtain~$\rho_{\phi_i\circ x_0}=\phi_i\circ\rho_{x_0}$
and the MT Hypothesis is still valid for~$\phi(U)$ in~$\phi(M)=M_{\phi_i\circ x_0}$. 
Without loss of generality may assume~$\phi_i=\phi_0$, that is~$\phi_i\circ x_0=x_0$.

We may also assume that~$[x_0,1]\in Z$ and thus that~$Z$ is contained in the image of~$X\times\{1\}$ 
in~$Sh_K(G,X)$.

\subsubsection{Covering by Siegel sets}
We choose a minimal $\QQ$-parabolic subgroup~$P$ of~$G$ 
and a maximal compact subgroup~$K_\infty$ of~$G(\RR)^+$, for 
instance~$K_{x_0}=Z_{G(\RR)}(x_0)$. We define~
\[
X^+=G(\RR)^+\cdot x_0\subset X
\]
and denote by
\[
S^+\subset Sh_K(G,X)
\]
the geometric component of~$Sh_K(G,X)$ which is the image of~$X^+\times\{1\}$.

See Def.~\ref{defi Siegel in G} for the definition of a Siegel set associated to~$P$ and~$K_\infty$.
It is known that there is a finite set~$\{g_1;\ldots;g_c\}\subseteq G(\QQ)$ and Siegel 
sets~$\FS_1,\ldots,\FS_c$ associated to~$g_1P{g_1}^{-1},\ldots,g_cP{g_c}^{-1}$ and~$K_\infty$
such that~$S^+$ is the image of~$\FS:=\FS_1\cup\ldots\cup \FS_c$.

%We recalled in~\ref{} a definition of Siegel sets associated to~$P$ and~$K$ 
%in~$G(\RR)$. If~$G$ is of~$\QQ$-rank zero, we allow~$P=G$ and
%Siegel sets will all be bounded subsets of~$G(\RR)$. This corresponds to the
%case the Shimura variety~$Sh_K(G,X)$ is compact (projective as an algebraic
%variety).

%The motivational property of Siegel sets is that there 
%are finitely many Siegel sets~$\FS_1,\ldots,\FS_c$ in~$G$, associated to finitely 
%many $G(\QQ)$-conjugates~$P_1,\ldots, P_c$ of~$P$ and~$K$, 
%such that~$S^+$ is the image of~$\FS=\FS_1\cup\ldots\cup \FS_c$.

%Arguing as in~\eqref{gen into geo} we may assum that for some~$i\in\{1;\ldots;c\}$ there is one~$i\in\{1;\ldots;c\}$
%such that, writing~$\FS:=\FS_i$, the following subset
%\[
%Z\cap [\FS:1]\cap \cH^g([x_0,1])\subset Z
%\]
%is Zariski dense in~$Z$.

For each~$\FS_i$, it is assumed that~$\Omega$ from Def.~\ref{defi Siegel in G} is a bounded \emph{semialgebraic} subset.

Let~$\FS_W=\FS/Z_{G(\RR)}(M)$ be the image of~$\FS$ in~$W^+(\RR)$.

The maps
\[
G(\RR)\xrightarrow{g\mapsto g Z_{G(\RR)}(M)} G(\RR)/Z_{G(\RR)}(M)\xrightarrow{gZ_{G(\RR)}(M)\mapsto g K_\infty} X=G(\RR)/K_\infty
\]
are real algebraic and thus semialgebraic. It follows that~$\FS_W$ is semi-algebraic, that its image~$\FS_X$ in~$X$ 
is semi-algebraic and that the map
\begin{equation}\label{p W X}
p_{W,X}:\FS_W\to \FS_X
\end{equation}
is semi-algebraic.
\subsubsection{O-minimality}
We use the theory of o-minimal structures and recall that the map
\[
\pi_{\FS,X}:\FS_X\to S^+
\]
is definable in the o-minimal structure~$\RR_{an,\exp}$ by~\cite{KUY}. As~\eqref{p W X} is semi-algebraic, it is definable in~$\RR_{an,\exp}$, and the 
following is definable in~$\RR_{an,\exp}$ as well
\[
\pi_{\FS,W}:=p_{W,X}\circ \pi_{\FS,X}:\FS_W\to \FS_X\to S^+.
\]
The algebraic variety~$Z$ is definable in~$\RR_{an,\exp}$ and its inverse image
\[
\wt{Z}_W=\stackrel{-1}{\pi_{\FS,W}}(Z) 
\]
is definable in~$\RR_{an,\exp}$ as well.

Because~$E$ is a field of definition for~$Z$, for every~$\sigma\in Gal(\ol{E}/E)$ 
and~$z\in Z(\ol{E})$ we have~$\sigma(z)\in Z$, and finally
\[
Gal(\ol{E}/E)\cdot z\subset Z.
\]
Assume now that~$z$ also belongs to~$\cH^g(x_0)$. For every
\[
z'\in Gal(\ol{E}/E)\cdot z
\]
we have~$z'\in Z\subset S^+$ and we can find~$\phi_{z'}\in W(\QQ)$ 
such that
\[
z'=[\phi_{z'}\circ x_0,1].
\]
Because~$\FS_X$ maps onto~$S^+$ we may assume that~$\phi_{z'}\circ x_0\in \FS_X$.
Equivalently, we have
\[
\phi_{z'}\in \FS_W.
\]
The set
\[
Q(z)=W(\QQ)\cap \stackrel{-1}{\pi_{\FS,W}}(Gal(\ol{E}/E)\cdot z)
\]
maps onto~$Gal(\ol{E}/E)\cdot z$ and we deduce
\begin{equation}\label{lift orbit}
\abs{Q(z)}\geq \abs{Gal(\ol{E}/E)\cdot z}.
\end{equation}
\subsubsection{Height bounds}\label{preuve height bounds}
We consider the affine embedding~$\iota:W\to\AAA^{\dim(M)\cdot \dim(G)}$ 
of~\S\ref{subsection specific height function}. Let~$H_W$ and~$H_f$ be as in~\eqref{H sans iota}.

%We choose some height function~$H_W$ on~$W(\QQ)$ and denote by~$H_{W,f}$
%the product of its local components at finite places, as in~\S\ref{subsection height generalities}. 
%We further assume that we have chosen the height function described in~\S\ref{subsection specific height function}.

We can of course assume that~$Z$ is infinite, and because
\[
\Sigma:=Z\cap \cH^g(x_0)
\]
is Zariski dense, it is infinite as well, and we can choose an infinite
sequence
\(
(z_n)_{n\in\ZZ_{\geq 1}}
\)
of pairwise distinct~$z_n\in\Sigma$. We also assume that this sequence is Zariski generic in~$Z$.

By hypothesis, Def.~\ref{defi MT in M} and~\ref{MT hypothesis} apply, and thus we invoke Th.~\ref{Galois lower bound} and, by~
Prop.~\ref{adelic galois orbits prop}, use it for Galois orbits.
We have
\[
H_{f}(\phi)\dom \abs{Gal(\ol{E}/E)\cdot [\phi\circ x_0,1]}\text{ on }W(\QQ).
\]
Thanks to the height comparison Th.~\ref{theorem type Orr}, we have
\begin{equation}\label{main proof height comparison}
H_{W}(\phi)\dom H_{f}(\phi)\text{ on }W(\QQ)\cap \FS_W.
\end{equation}
It follows
\[
H_{W}(\phi)\dom \abs{Gal(\ol{E}/E)\cdot [\phi\circ x_0,1]}\text{ on }W(\QQ)\cap\FS_W.
\]
More precisely, there are~$a,b\in\RR_{>0}$ such that
\[
\forall\phi\in W(\QQ)\cap\FS_W,
a+H_{W}(\phi)^b\leq \abs{Gal(\ol{E}/E)\cdot [\phi\circ x_0,1]}.
\]
Using~\eqref{lift orbit} we deduce
\begin{equation}\label{main proof precise Q bound}
a+H_{W}(\phi_{z_n})^b\leq \abs{Q(z_n)}.
\end{equation}

From~Prop.~\ref{prop inv galois height} we have
\[
\forall z'\in Gal(\ol{E}/E)\cdot z_n, H_f(\phi_{z'})=H_f(\phi_{z_n})
\]
and because~$H_f(\phi)$ only depends on~$[\phi\circ x_0,1]$ we have
\[
\forall \phi\in Q(z_n),~ H_f(\phi)=H_f(\phi_{z_n}).
\]

We make~\eqref{main proof height comparison} precise by choosing~$a',b'$ such that
\begin{equation}\label{main proof height comparison precise}
\forall\phi\in W(\QQ)\cap\FS_W, H_{W}(\phi)\leq a'+H_f(\phi)^{b'}.
\end{equation}
For~$\phi\in Q(z_n)\subset W(\QQ)\cap\FS_W$ we get
\[
H_W(\phi)\leq a'+H_f(\phi)^{b'}=a'+H_f(\phi_{z_n})^{b'}.
\]

Writing~$k(n)=H_f(\phi_{z_n})$, we deduce from the above that
the subset~$Q(z_n)\subseteq \wt{Z}\cap W(\QQ)$ 
contains at least~$a+k(n)^b$ points of~$H_W$-height at most~$a'+k(n)^{b'}$.

Because the~$z_n$ are distinct, so are the inverse images~$\phi_{z_n}$,
and by Northcott theorem we deduce that~$H_W(\phi_{z_n})\to+\infty$,
and thus~$k(n)\to+\infty$.

We are ready to use the Pila-Wilkie theorem.

\subsubsection{Pila-Wilkie theorem}
We use the form Theorem~\ref{PilaWilkie} of the Pila-Wilkie theorem.
We denote~$K^{\RR}_\infty$ the real algebraic group corresponding to~$K_\infty$,
and~$X_\RR$ the algebraic variety~$G_\RR/K^{\RR}_\infty$ over~$\RR$
(we have~$X\subset X_\RR(\RR)$). We apply Theorem~\ref{PilaWilkie}
to the morphism~$p:W=G_\RR/Z_{G_\RR}(M)\to X_\RR= G_\RR/K^{\RR}_\infty$
and the definable subset
\[
\wt{Z}_X:=\stackrel{-1}{\pi_{\FS,X}}(Z)\subset X\subset X_\RR(\RR).
\]
We deduce for every~$n$ that
\[
\abs{Q(z_n)\cap (\wt{Z}_X\smallsetminus {\wt{Z}_X}^{alg})}=(a'+H_f(\phi_{z_n})^{b'})^{o(1)}=o(\abs{Q(z_n)}) .
\]
Thus, for~$n\gg 0$, we have
\[
Q(z_n)\cap {\wt{Z}_X}^{alg}\neq \emptyset.
\]
In other terms, for almost every~$n$, there exist~$\phi\in Q(z_n)$, and a non-zero dimensional semialgebraic subset~$A_n\subset \wt{Z}_X$, such that~$\phi\circ x_0\in A_n$.

We will now use the Hyperbolic Ax-Lindemann-Weierstrass theorem.

\subsubsection{Functional transcendance}
According to Ax-Lindemann-Weierstrass theorem (see \cite{KUY}), that for~$n\gg0$, there exists a weakly special
subvariety~$S'_n$ of~$S^+$ such that
\[
z'_n\in \pi_{\FS,X}(A_n)\subset S'_n\subset Z.
\]
One can check that a weakly special subvariety containing a~$\ol{E}$-valued point
is defined over~$\ol{E}$. It follows that this~$S'_n$
is defined over~$\ol{E}$, and applying~$\sigma\in Gal(\ol{E}/E)$
such that~$\sigma(z'_n)=z_n$, the conjugated subvariety~$S_n=\sigma(S'_n)$
will be: weakly special, contained in~$Z$ and containing~$z_n$.

Because the sequence~$z_n$ is generic in~$Z$, the family~$(S_n)_{n\geq 0}$
is Zariski dense in~$Z$. 

Because~$A_n$ has non-zero semialgebraic dimension, and~$\pi_{\FS,X}$ 
has finite fibers, the image~$\pi_{\FS,X}(A_n)$ has non-zero semialgebraic dimension,
and~$S'_n$ has non-zero dimension as a variety, and~$S_n$ also.

We are ready to use the so-called geometric part of André-Oort conjecture.

\subsubsection{Geometric André-Oort}
We reuse the notations of~\S\ref{partition factorisation}

From the Geometric part of André-Oort conjecture from~\cite{UllmoAxLin,RU}, there exists a partition~$\{1;\ldots;c\}=I\sqcup J$,
with~$I\neq \emptyset$, but possibly~$J=\emptyset$, such that
we have a factorisation
\[
Z=S_1\times Z_J\subset Sh_{K_I}(G_I,X_I)\times Sh_{K_J}(G_J,X_J),
\]
where~$S_1$ is a geometric component of~$Sh_{K_I}(G_I,X_I)$, and~$Z_J$ is a subvariety of~$Sh_{K_J}(G_J,X_J)$.

Because we assumed that~$Z$ has no non trivial factorisation, the partition~$\{1;\ldots;c\}=I\sqcup J$
is trivial. We must have~$J=\emptyset$, $I=\{1;\ldots;c\}$.
Equivalently~$Z=S_1$. In other words~$Z$ is special, and in particular is weakly special.

This finishes the proof of Theorem~\ref{thm intro}.
 
\subsection{Refined Pila-Wilkie theorem.} \label{PW}The following is a variant of Pila-Wilkie Theorem,
which replaces the ``block version'' of Pila-Wilkie Theorem used by Orr. We believe this variant is easier to understand
and use, and will be of independent interest.

We deduce the following from~\cite[Th.~1.7]{Pila}.
\begin{theorem}\label{PilaWilkie}
Let~$W$ be an affine algebraic variety defined over~$\QQ$, let~$X$ be an affine algebraic variety over~$\RR$
and let~$p:W_\RR\to X$ be a morphism of algebraic varieties defined over~$\RR$.

Let~$Z\subset X(\RR)$ be a definable subset, and denote~$Z^{alg}$
be the union of the semialgebraic subsets of~$X(\RR)$ which are contained in~$Z$ and of non-zero dimension.

We consider a height function~$H_W$ on~$W(\QQ)$ associated to some affine embedding.
Then 
\[
\abs{(Z\smallsetminus Z^{alg})\cap p(\{w\in W(\QQ) : H_W(w)\leq T\})}=T^{o(1)}.
\]
Explicitly, for every~$\epsilon\in\RR_{>0}$, there exists~$C(\epsilon,Z)\in\RR_{>0}
$, such that
\[
\forall T\gg0,\abs{(Z\smallsetminus Z^{alg})\cap p(\{w\in W(\QQ) : H_W(w)\leq T\})}\leq C(\epsilon,Z)\cdot T^{\epsilon}.
\]
\end{theorem}
\subsubsection*{Comment} The theorem still holds with a semi-algebraic map~$p:W(\RR)\to X(\RR)$ instead of the real algebraic~$p:W_\RR\to X$. This slight generalisation will not be needed.

The height function we use here is denoted~$H^{\text{proj}}$ by Pila, and is not the Height function he uses in his statements.
As mentioned in the introduction of~\cite{Pila}, it is possible to invoke his statements with~$H^{\text{proj}}$ instead.

\begin{proof}We choose affine embeddings
\[
W\subseteq \AAA^n\text{ and }X\subseteq \AAA^m
\]
defined over~$\QQ$ and~$\RR$. We can then write the morphism
\[
p(w_1,\ldots,w_n)=(P_1(w_1,\ldots,w_n),\ldots,P_m(w_1,\ldots,w_n))
\]
with polynomials~$P_1,\ldots,P_m\in\RR[T_1,\ldots,T_n]$.
Let~$E$ be the finite dimensional $\QQ$-vector subspace of~$\RR$
generated by the coefficients of these polynomials.

We have 
\[
p(W(\QQ))\subseteq E^m.
\]
We choose an isomorphism~$\iota:E\to \QQ^d$ of~$\QQ$ vector spaces. For every~$P_i$ the map
\[
\iota\circ P_i:W(\QQ)\to E\to \QQ^d
\]
is polynomial with coefficients in~$\QQ$. This can be checked for every monomial of~$P_i$. The height on~$E^m$ considered in~\cite[Th.~1.7]{Pila} can be written, with our notations, 
\[
H_E=H\circ (\iota,\ldots,\iota)
\]
where~$H$ is the usual height on~$\QQ^{d\cdot m}$.
It follows from the general ``functoriality'' properties of heights of~\S\ref{sec: funct heights} that
\[
H_E\circ p\dom H_W\text{ on }W(\QQ).
\]
Explicitly, for some~$a,b\in\RR_{>0}$ we have
\[
p(\{W\in W(\QQ) : H_W(w)\leq T\})\subseteq \{e\in E^m : H_E(e)\leq a+T^b\}.
\]
We apply~\cite[Th.~1.7]{Pila} and obtain
\begin{multline*}
\abs{(Z\smallsetminus Z^{alg})\cap p(\{w\in W(\QQ) : H_W(w)\leq T\})}\\ \leq 
\abs{(Z\smallsetminus Z^{alg})\cap \{e\in E^m : H_E(e)\leq a+T^b\}}=T^{o(1)}.\qedhere
\end{multline*}
\end{proof}

%We assume the following:
%\begin{itemize}
%\item
%$Z$ is strongly Hodge generic.
%\item
%The group $G$ is of adjoint type.
%\item
%$Z$ is not of the form
%$Z= S_1 \times Z' \subset S_1 \times S_2$
%where $S_1 \times S_2$ is a decomposition of $S$.
%\end{itemize}
%
%
%Recall that we have a map
%$$
%W(\RR) \lto X
%$$
%given by 
%$$
%g \mapsto g x_0 g^{-1}
%$$
%inducing an injection
%$$
%W(\QQ) \hookrightarrow X.
%$$
%
%Recall that we have $\FS \subset W(\RR)$ a Siegel set. 
%Its image in $X$ is a Siegel set for the action of $\Gamma$,
%we still denote it by $\FS$.
%
%Let $\wt{Z}$ be $\pi^{-1}Z \cap \FS$ and $\wt{Z}_W$ be the intersection of the preimage of $\pi^{-1}Z$ in $W(\RR)$
%with $\FS$.
%
%Both $\wt{Z}$ and $\wt{Z}_W$ are definable in $\RR_{an,exp}$ because $\pi$ restricted to $\FS$ is 
%and because the map $W(\RR) \lto X$ is semialgebraic.

\appendix

\section{Exponentials of~$p$-adic matrices}\label{section exponential}

In this section we fix a prime~$p$, an integer~$d\in\ZZ_{\geq1}$ and
denote by~$M_d(\QQ_p)$ the space of square matrices of size~$d$ with entries in~$\QQ_p$. 
% and~$1_d$ the matrix of the identity. 
For~$Z\in M_d(\QQ_p)$ we denote by~$\chi_Z(T)=\det(TZ-1)\in\QQ_p[T]$ its characteristic polynomial.
Let~$\abs{~}$ be the normalised absolute value on~$\QQ_p$, extended to~$\overline{\QQ_p}$: we have~$\abs{p}=1/p$ and~$\abs{1/d}\leq d$ for~$d\in\ZZ_{\geq 1}$.
We denote the \emph{norm} of~$Z$, and the \emph{local height} of~$Z$ by
\[
\norm{Z}=\max_{1\leq i,j\leq d}\abs{Z_{i,j}}\text{ and }H_p(Z)=\max\{1;\norm{Z}\}=H_p(1+Z).
\]
We define, whenever the corresponding series converges in~$M_d(\QQ_p)$,
\[
\exp(Z)=\sum_{n\in\ZZ_{\geq0}} \frac{1}{n!}\cdot Z^n\text{ and }\log(1+Z)
=
-\sum_{n\in\ZZ_{\geq 1}}\frac{(-1)^n}{n}\cdot Z^n.
\]
It is well known (see \cite[Ch.\,5.\,\S4.1]{Robert}) that, on~$\CC_p$, the series~$\exp(T)$ has radius of convergence~$\abs{p}^{\frac{1}{p-1}}$ and the series~$\log(1+T)$ has radius of convergence~$1$. It is also true that~$\exp(Z)$, resp.~$\log(1+Z)$ converges if and only if the eigenvalues
of~$Z$ are in the open disc of convergence of~$\exp(T)$ resp.~$\log(1+T)$. (For the archimedean case, see~\cite[\S1]{FunctionMatrices}. The relevant arguments carry over to ultrametric fields.)

\begin{proposition}\label{Proposition bornes log matrice}
Let~$Y\in M_d(\QQ_p)$ be such that~$\log(1+Y)$ converges.
Then
\begin{equation}\label{conclusion chi mod p}
\chi_Y(T)\in T^d+p\ZZ_p[T].
\end{equation}
Let~$Y\in M_d(\QQ_p)$ be such that
\begin{equation}\label{hypothese chi mod p}
\chi_Y(T)\in T^d+p\ZZ_p[T].
\end{equation}
Then~$\log(1+Y)$ converges and we have
\begin{itemize}
\item in general,
\begin{equation}\label{log bound general}
\norm{\log(1+Y)} \leq d\cdot H_p(Y)^{d-1},
\end{equation}
\item and for~$p>d$, the sharper estimate
\begin{equation}\label{log bound precise}
\norm{\log(1+Y)} \leq H_p(Y)^{d-1}.
\end{equation}
\end{itemize}
\end{proposition}
We deal with the first conclusion~\eqref{conclusion chi mod p}.
\begin{proof}Let~$\lambda_1,\ldots,\lambda_d$ be the eigenvalues of~$Y$, with repetitions. 
As can be seen on a Jordan form after passing to~$\CC_p$, the series~$\log(1+Y)$ converges 
if and only if every~$\log(\lambda_1),\ldots,\log(\lambda_d)$ converges.
As the radius of convergence of~$\log(1+T)$ is~$1$, this means
\begin{equation}\label{spectre converge log}
\forall i\in\{1;\ldots;d\}, \abs{\lambda_i}<1.
\end{equation}
Let~$K=\QQ_p(\lambda_1,\ldots,\lambda_d)$, let~$O_K$ be its ring of integers, and~$\mathfrak{m}_K$
be the maximal ideal of~$O_K$. Then~\eqref{spectre converge log} means
\[
\{\lambda_1;\ldots;\lambda_d\}\subseteq \mathfrak{m}_K.
\]
We deduce that the non leading coefficients of
\[
\chi_Y(T)=\prod_{i=1}^d (T-\lambda_i)
\]
are in~$\mathfrak{m}_K$. 
We recall that~$\QQ_p\cap \mathfrak{m}_K=p\ZZ_p$ and~$\chi_Y(T)\in\QQ_p[T]$.
We conclude that
\[
\chi_Y(T)\in \QQ_p[T]\cap\bigl( T^d+ \mathfrak{m}_K\cdot O_K[T]\bigr)=T^d +p\cdot \ZZ_p[T].\qedhere
\]
\end{proof}
We have proved~\eqref{conclusion chi mod p} and before proving the rest of Prop.~\ref{Proposition bornes log matrice}, 
we prove an estimate on~$\norm{Y^n}$ for~$n\in\ZZ_{\geq0}$.
\begin{proof} We consider
\[
A:=\ZZ_p+\ZZ_p \cdot Y+\ldots+\ZZ_p\cdot Y^{d-1}.
\]
By hypothesis, we have~$\chi_Y(T)=c_0+\ldots+c_{d-1}T^{d-1}+T^d$ with~$c_0,\ldots,c_{d-1}\in p\ZZ_p$. Let us first check that
\begin{equation}\label{A stable par X}
Y A\subseteq A 
\end{equation}
on a generating family: for~$0\leq i<d-1$ we have~$Y\cdot Y^i\in A$ by construction;
for~$i=d-1$ the identity~$\chi_Y(Y)=0$ can be rearranged into
\begin{equation}\label{Xd dans pA}
Y^d=-c_0+\ldots-c_{d-1}Y^{d-1}\in p A.
\end{equation}

Repeated use of~\eqref{A stable par X} implies that, for~$i\in\ZZ_{\geq 0}$, we have~$Y^i A\subseteq A$.
We deduce~$Y^i pA\subseteq pA$. But~$Y^d\in pA$ by~\eqref{Xd dans pA},
hence~$Y^d\cdot Y^i= Y^i\cdot Y^d\in pA$. Applied to~$i=0,\ldots,d-1$ it implies
\(
Y^d A\subseteq p A
\)
and by induction
\(
(Y^d)^k A\subseteq p^k A.
\)
We deduce again that~$Y^i\cdot (Y^d)^k\in p^k A$.
Writing~$n=k\cdot d+i$ with~$k=[\frac{n}{d}]$, we get the formula
\[
Y^n\in p^{[\frac{n}{d}]}A
\]
and the bound
\begin{equation}\label{YN 1}
\norm{Y^n}\leq \abs{p}^{[\frac{n}{d}]} \cdot \norm{A}
\text{ where }
\norm{A}:=\max_{a\in A}\norm{a}.
\end{equation}
Using the ultrametric inequality~$\norm{X+Z}\leq \max\{\norm{X};\norm{Z}\}$ and submultiplicativity~$\norm{X\times Z}\leq\norm{X}\cdot\norm{Z}$
of the norm, we get
\begin{equation}\label{YN 2}
\norm{A}\leq \max\{\norm{Y^0};\ldots;\norm{Y^{d-1}}\}\leq \max\{1;\ldots;\norm{Y}^{d-1}\}=H_p(Y)^{d-1}.
\end{equation}
\end{proof}
We apply our estimate to the series~$\log(1+T)$ and finish the proof
of Proposition~\ref{Proposition bornes log matrice}.
\begin{proof}
For the series
\(
\log(1+Y)
\)
the above~\eqref{YN 1} and~\eqref{YN 2} imply the bound
\[
\norm*{\frac{(-1)^n}{n}\cdot Y^n}\leq \abs*{\frac{1}{n}}\cdot \abs{p}^{[\frac{n}{d}]} \cdot H_p(\abs{Y})^{d-1}.
\]
We note that
\(
\lim_{n\to\infty}\abs*{\frac{1}{n}}\cdot \abs{p}^{[\frac{n}{d}]}=0
\)
which implies that~$\log(1+Y)$ converges, and that
\[
\max_{n\in\ZZ_{\geq1}}\abs*{\frac{1}{n}}\cdot \abs{p}^{[\frac{n}{d}]}=\abs*{\frac{1}{d-1}}\cdot \abs{p}^{[\frac{d-1}{d}]}=\abs*{\frac{1}{d-1}}.
\]
By the ultrametric inequality and previous estimates, 
\begin{equation}\label{borne log padique}
\log(1+Y)\leq 
\sup_{n\in\ZZ_{\geq1}}
\norm*{\frac{(-1)^n}{n}\cdot Y^n}
\leq 
\abs*{\frac{1}{d-1}}\cdot H_p(Y)^{d-1}.
\end{equation}
As we used the normalised~$p$-adic norm, we have
\(
\abs*{\frac{1}{d-1}}\leq d-1\leq d
\)
in general, and~$\abs*{\frac{1}{d-1}}=1$ if~$p\geq d$. This gives~\eqref{log bound general} and~\eqref{log bound precise} respectively.
\end{proof}

The main statement of this section will require the following observation.
\begin{lem}\label{lemme composee log exp matrice}
Let $Z \in M_d(\QQ_p)$ be such that $\exp(Z)$ converges and let us write~$\exp(Z)=1+Y$. Then~$\log(1+Y)$ converges and
$$
\log(1+Y) = Z.
$$
\end{lem}
\begin{proof}
For $d=1$, it is \cite{Robert}, \S 5, prop. 3.

For $d>1$, it is \cite{Robert}  \S 6.1.1 applied to
$(\partial/\partial Y)^i \log(1+Y)\circ \exp$.
\end{proof}

The following statement is one of our main tools for proving lower bounds for Galois orbits.
\begin{theorem}[Lemma of the exponentials] \label{lemme_exp}
Let $X \in M_d(\QQ_p)$ be such that $\exp(X)$ converges
and denote by~$\exp(X)^\ZZ$ the subgroup generated by~$\exp(X)$
in~$GL_d(\QQ_p)$.

Then
\begin{itemize}
\item in general, we have
\begin{equation}\label{borne orbite general}
[\exp(X)^\ZZ:\exp(X)^\ZZ\cap GL_d(\ZZ_p)]\geq H_p(X)/d
\end{equation}
\item and, if $p>d$, we have more sharply 
\begin{equation}\label{borne orbite precise}
[\exp(X)^\ZZ:\exp(X)^\ZZ\cap GL_d(\ZZ_p)]\geq H_p(X).
\end{equation}
\end{itemize}
\end{theorem}
\begin{proof} For every~$i\in\ZZ$, we know that if~$\exp(X)$ converge, then~$\exp(iX)$ converges as well, and
we have 
\[
\exp(iX) = \exp(X)^i.
\]
By~\ref{lemme composee log exp matrice}, with~$Y_i=\exp(i\cdot X)-1$, we have convergence and identity
\[
\log(1+Y_i)=i\cdot X.
\]
The Proposition~\ref{Proposition bornes log matrice} gives
\begin{align}
\norm{i\cdot X}
&=\norm{\log(1+Y_i)}\leq d \cdot H_p(1+Y_i)^{d-1}
\label{calcul general}
\\
\text{ and, if~$d\leq p$, }
\norm{i\cdot X}
&=\norm{\log(1+Y_i)}\leq  H_p(1+Y_i)^{d-1}.
\label{calcul precis}
\end{align}
Assume that
\[
i=[\exp(X)^\ZZ:\exp(X)^\ZZ\cap GL_d(\ZZ_p)]<+\infty.
\]
Then~$H_p(1+Y_i)=H_p(\exp(X)^i)=1$, and~\eqref{calcul general}, resp.~\eqref{calcul precis}, specialises to
\[
\abs{i}\cdot \norm{X}\leq d \text{, resp. }\abs{i}\cdot \norm{X}\leq 1.
\]
Recall that~$\abs{i}\leq \frac{1}{i}$ as we use the normalised~$p$-adic absolute value. The conclusions~\eqref{borne orbite general}, resp.~\eqref{borne orbite precise}, follow.
\end{proof}
We finish with a sufficient criterion for~$\exp(X)$ to converge.
\begin{theorem}\label{criter exp converge}
Let~$X$ be a matrix in~$M_d(\QQ_p)$ and~$b\in\ZZ_{\geq1}$ be such that
$$
\chi_X(T) \in T^d+p^k\ZZ_p[T]\text{ and }d <k(p-1).
$$
Then $\exp(X)$ converges.
\end{theorem}

\begin{proof}By the usual criterion, it is sufficient to prove that every eigenvalue~$\lambda$ of~$X$ is in the open disc of convergence for~$\exp(T)$. This amounts to proving the inequality~$\abs{\lambda}<\abs{p}^{\frac{1}{p-1}}$. 

For any eigenvalue~$\lambda$ of~$X$, we have~$\chi_X(\lambda)=0$ hence~$\lambda^d \in p^k \ZZ_p[\lambda]$  by assumption.
It follows~$|\lambda|^d \leq |p|^k$, that is~$\abs{\lambda}\leq \abs{p}^{\frac{k}{d}}$.
Using the inequality~$d <k(p-1)$, it implies~$\abs{\lambda}<\abs{p}^{\frac{1}{p-1}}$. 
%As~$\abs{p}^{\frac{1}{p-1}}$ is the radius of convergence of~$\exp(T)$, and~$\exp(\lambda)$ converges. As it is equally true every eigenvalue in the spectrum, we deduce that~$\exp(X)$ converges.
\end{proof}

\section{Heights  bounds for  adelic orbits of linear groups}
\label{section adelic bounds}

Our bound on~$p$-adic exponentials is combined with structure theory of linear algebraic groups to obtain the
following general lower bound. It is applied to Galois orbits in section~\ref{section MT adelic}.
\begin{theorem}\label{general global bounds}
Let~$M\leq GL(N)$ be a linear algebraic subgroup defined over~$\QQ$, denote by~$\phi_0: M\to GL(N)$
the identity morphism and~$W$ the~$GL(N)$-conjugacy class of~$\phi_0$.
We define
\[
M(\widehat{\ZZ})=M(\AAA_f)\cap GL(N,\widehat{\ZZ})\text{ and }\mathfrak{m}_{\widehat{\ZZ}}=\mathfrak{m}\tens \AAA_f\cap \mathfrak{gl}(N,\widehat{\ZZ}).
\]

We consider the standard  Weil $\AAA_f$-height function, see~\eqref{notation:AF height},
\[
H_f:\Hom(\mathfrak{m}\tens\AAA_f,\mathfrak{gl}(N)\tens\AAA_f)\to \ZZ_{\geq1}
\]
given by~$H_f(\Phi)=\min\{n\in\ZZ_{\geq1} : n\Phi(\mathfrak{m}_{\widehat{\ZZ}})\subset \mathfrak{gl}(N,\widehat{\ZZ})\}.$

There exists~$c=c(\phi_0)\in\RR_{>0}$ such that, as~$\phi$ ranges through~$W(\AAA_f)$, we have
\begin{align}\label{explicit global adelic bound}
[\phi(M(\widehat{\ZZ})):\phi(M(\widehat{\ZZ}))\cap GL(N,\widehat{\ZZ})]
\geq \frac{1}{c^{\omega(H_f(d\phi))}}\cdot H_f(d\phi).
\end{align}
(Where~$\omega(n)$ counts the number of prime factors of~$n$.)
\end{theorem}
The proof of Th.~\ref{general global bounds} will start in~\S\ref{sec:proofB1}.
We deduce from Th.~\ref{general global bounds} the following.
%In the other sections of this article we will only retain from the explicit bound~\eqref{explicit global adelic bound} the qualitative consequence~\eqref{coro global bounds dom} below. The observation~\eqref{invariance Hf} will be used in section~\ref{section invariant height}.
\begin{cor}\label{B1} We have
\begin{equation}\label{coro global bounds o}
[\phi(M(\widehat{\ZZ})):\phi(M(\widehat{\ZZ}))\cap GL(N,\widehat{\ZZ})]\geq H_f(d\phi)^{1-o(1)}
\end{equation}
and, if~$M$ is reductive and connected and~$\iota:W\to \AAA^d$ is an affine embedding, then, as~$\phi$ ranges through~$W(\AAA_f)$,
%$H_{W,f}$ is a Height function on~$W(\AAA_f)$,
\begin{equation}\label{coro global bounds dom}
H_{\iota,f}(\phi)\approx H_f(d\phi) \dom [\phi(M(\widehat{\ZZ})):\phi(M(\widehat{\ZZ}))\cap GL(N,\widehat{\ZZ})].
\end{equation}

Furthermore, for every~$\Phi\in\Hom(\mathfrak{m}\tens\AAA_f,\mathfrak{gl}(N)\tens\AAA_f)$, we have
\begin{equation}\label{invariance Hf}
\forall m\in M(\widehat{\ZZ}), g\in G(\widehat{\ZZ}), H_f(g\circ \Phi\circ m)=H_f(\Phi).
\end{equation}
\end{cor}
\begin{proof}[Proof of Cor.~\ref{B1}]One passes from~\eqref{explicit global adelic bound} to~\eqref{coro global bounds o}
by recalling the known estimate (see~\cite[22.10]{HW})
\[
c^{\omega(n)}\leq n^{\abs{\log(c_2)}\cdot \frac{1+o(1)}{\log\log n}}=n^{o(1)}.
\]
As for~\eqref{coro global bounds dom},
we know that~$W$ is affine as~$M$ is reductive, and~$\phi\mapsto d\phi$ is an affine embedding because~$M$ is connected.
Lastly, two heights functions on~$W$ are polynomially equivalent, so we may replace~$H_{W,f}(\phi)$ by~$H_{f}(d\phi)$ and
this follows from~\eqref{coro global bounds o}.

The identity in~\eqref{invariance Hf} follows from the observations~
\[m\cdot \mathfrak{m}_{\widehat{\ZZ}}=\mathfrak{m}_{\widehat{\ZZ}}\text{, and }g^{-1}\cdot \mathfrak{gl}(N,{\widehat{\ZZ}})=\mathfrak{gl}(N,{\widehat{\ZZ}})\]
and the defining property we provided: we have~$n\cdot g\cdot\Phi(m\mathfrak{m}_{\widehat{\ZZ}})\subset \mathfrak{gl}(N,\widehat{\ZZ})$ if and only if
\[
n\Phi(\mathfrak{m}_{\widehat{\ZZ}})=
n\Phi(m\mathfrak{m}_{\widehat{\ZZ}})\subset g^{-1}\mathfrak{gl}(N,\widehat{\ZZ})=\mathfrak{gl}(N,\widehat{\ZZ}).
\qedhere\]
\end{proof}

The combination of Th.~\ref{C1} \eqref{Extra bound 1} with~\eqref{coro global bounds dom} gives the following.
\begin{theorem}
Let~$M\leq GL(N)$ be a connected reductive linear algebraic subgroup defined over~$\QQ$, denote~$\phi_0: M\to GL(N)$
the identity morphism and~$W$ the~$GL(N)$-conjugacy class of~$\phi_0$,  and let~$\iota:W\to \AAA^d$ be an affine embedding. Then, as~$\phi$ ranges through~$W(\AAA_f)$,
\begin{equation}
H_{\iota,f}(\phi)\approx H_f(d\phi) \approx [\phi(M(\widehat{\ZZ})):\phi(M(\widehat{\ZZ}))\cap GL(N,\widehat{\ZZ})].
\end{equation}
\end{theorem}
\subsection{Proof of Theorem~\ref{general global bounds}.}\label{sec:proofB1}
The global theorem~\ref{general global bounds} will follow directly from~\eqref{explicit local adelic bound} in the analogous local theorem below.

\begin{theorem}\label{general local bounds}
We keep~$M$, $\phi_0$, $W$ and~$H_f$ as in Theorem~\ref{general global bounds}.

For every prime~$p$, let~$H_p:\Hom(\mathfrak{m}\tens\AAA_f,\mathfrak{gl}(N)\tens\AAA_f)\to \ZZ_{\geq1}$ be given by~$H_p(\Phi)=\min\{p^k\in p^{\ZZ_{\geq1}} : p^k\Phi(\mathfrak{m}_{\ZZ_p})\subset \mathfrak{gl}(N,\ZZ_p)\}$.

There exists~$c=c(\phi_0)\in\RR_{>0}$ such that, for every prime~$p$, and every~$\phi\in W(\QQ_p)$,
\begin{equation}\label{explicit local adelic bound}
[\phi(M({\ZZ}_p)):\phi(M({\ZZ}_p))\cap \GL(N,{\ZZ}_p)]\geq\frac{H_p(d\phi)}{c}
\end{equation}
and if~$\mathfrak{m}_{\ZZ_p}$ is generated over~$\ZZ_p$ by nilpotent elements and~$p> N$,
\begin{equation}\label{explicit local adelic bound nilpotent basis}
[\phi(M({\ZZ}_p)):\phi(M({\ZZ}_p))\cap \GL(N,{\ZZ}_p)]\geq H_p(d\phi).
\end{equation}
\end{theorem}
Here is how to deduce Th.~\ref{general global bounds} from Th.~\ref{general local bounds}.
\begin{proof} Let us multiply the inequalities~\eqref{explicit local adelic bound} for the~$\omega(H_f(d\phi))$ primes dividing~$H_f(d\phi)$
with the trivial inequalities
\[
[\phi(M({\ZZ}_p)):\phi(M({\ZZ}_p))\cap \GL(N,{\ZZ}_p)]\geq 1
\]
for all the other primes. Then one can identify the product on both sides with the corresponding sides of~\eqref{explicit global adelic bound}.
\end{proof}

Theorem~\ref{general local bounds} will follow from different cases gathered in Theorem~\ref{cases local bounds}.
%The conclusion~\ref{local bound 1} and~\ref{local bound 3} would suffice for~\eqref{coro global bounds dom} and our need in this article. 
%The conclusion~\ref{local bound 2}
%is used to cover non reductive algebraic groups~$M$, and in prevision for future use. These conclusions could suffice for
%achieving~\eqref{coro global bounds o} with a lower estimate~$H_f(d\phi)^{\frac{1}{2}-o(1)}$. We will anyway derive
%the conclusion~\ref{local bound 3} from conclusion~\ref{local bound 4}, and in the given form,~\ref{local bound 4} will
%take part in getting Theorem~\ref{general global bounds} with estimates we have given.
\begin{theorem}\label{cases local bounds}We keep the notations from Theorem~\ref{general local bounds}.
For every prime~$p$, let~$K_p:=GL(N,\ZZ_p)$ and, for any~$U\leq G(\QQ_p)$, let~$[U]_p:=[U:U\cap K_p]$.
We write~$N^*=\mathrm{lcm}(1,\ldots,N)$ so that~$\abs{1/N^*}_p=p^{[log_p(N)]}$ and~$\abs{N^*}_p=1$ 
if~$p>N$.

\begin{enumerate}
\item \label{local bound 1}
 For every prime~$p$ we have~$\exp(2p\mathfrak{m}_{\ZZ_p})\leq M(\ZZ_p)$ and
\begin{equation}\label{conclusion 1}
[\phi(\exp(2p\mathfrak{m}_{\ZZ_p}))]_p\geq \abs{2pN^*}_p\cdot H_p(d\phi)\geq \frac{1}{2Np}\cdot H_p(d\phi).
\end{equation}
\item \label{local bound 2} Assume that~$M$ is unipotent or more generally that~$\mathfrak{m}_{\ZZ_p}$ is generated over~$\ZZ_p$
by nilpotent elements, then
\begin{equation}\label{conclusion 2}
[\phi(M({\ZZ}_p))]_p\geq \abs{N^*}_p\cdot H_p(d\phi).
\end{equation}
%\item \label{local bound 3} Assume that~$M$ is reductive or more generally generated by algebraic tori. There is~$c_2=c_2(\phi_0)\in\RR_{>0}$ such that for every prime~$p$, and every~$\phi\in W(\QQ_p)$,
%such that 
%\[
%\text{if~$H_p(d\phi)\neq 1$ then }[\phi(M({\ZZ}_p))]_p\geq
%\frac{H_p(d\phi)}{c}.
%\]
\item \label{local bound 4} Assume that~$M$ is an algebraic torus. There is~$c_2=c_2(\phi_0)\in\RR_{>0}$ such that for every prime~$p$, and every~$\phi\in W(\QQ_p)$,
\begin{equation}\label{conclusion 3}
\text{if~$H_p(d\phi)\neq 1$ then }
\abs*{\frac{\phi(M(\ZZ_p))}{\phi(\exp(2p\mathfrak{m}_{\ZZ_p}))\cdot\phi(M(\ZZ_p))\cap K_p}}
\geq
\frac{p}{c_2}.
\end{equation}
\end{enumerate}
\end{theorem}
We deduce the Theorem~\ref{general local bounds} from Theorem~\ref{cases local bounds}.
\begin{proof} The bound~\eqref{explicit local adelic bound nilpotent basis} follows from~\eqref{conclusion 2}, and the observation that~$\abs{N^*}_p=p^{[\log_p(N)]}=p^0=1$ for~$p>N$.

Let~$U$ be the unipotent radical of~$M^0$ and~$L$ be a reductive Levi subgroup of~$M^0$ so that we have the Levi decomposition~$\mathfrak{m}=\mathfrak{u}+\mathfrak{l}$. By the principle of~\S\ref{Subgroup principle} we may assume~$M=U$ or~$M=L$. 

In the first case~$M=U$, one deduces~\eqref{explicit local adelic bound}, with~$c=N^*\geq\abs{1/N^*}_p$, from \eqref{conclusion 2}.

In the second case,~$M=L$ is reductive, and thus generated by algebraic tori. By the principle~\ref{Subgroup principle}
we may assume that~$M$ is a torus.

Let us mention a simpler argument giving the following weaker conclusion, which is sufficient for the purpose of this article:
\begin{equation}\label{explicit local adelic bound square root}
[\phi(M({\ZZ}_p)):\phi(M({\ZZ}_p))\cap GL(N,{\ZZ}_p)]\geq\frac{H_p(d\phi)^{1/2}}{c_2}.
\end{equation}
\begin{proof}We know that~$H_p(d\phi)$ is a power~$p^k$ of~$p$. For~$k=0$, we may take~$c=1$.
For~$k=1$ we deduce from conclusion~\eqref{local bound 4} of Th.~\ref{cases local bounds} that\footnote{The bound~\eqref{EY weak form} is from~\cite[Prop.~4.3.9]{EdYa}.}
\begin{equation}\label{EY weak form}
[\phi(M(\ZZ_p))]_p\geq p/c_2=H_p(d\phi)/c_2.
\end{equation}
For~$k\geq 2$, we have~$H_p(d\phi)/p\geq \sqrt{H_p(d\phi)}$ and we take~$c_2=2N$ and use~\eqref{conclusion 1}.
\end{proof}
We now explain how to improve upon the exponent~$1/2$.

We suppose that~$p$ is large enough, that~$p\neq 2$, and that the reduction~$T_{\FF_p}$ of the torus~$T=M$ is a torus over~$\FF_p$. 
Then~$T_{\FF_p}(\FF_p)$ is diagonalisable over~$\ol{\FF_p}$ and its elements have order prime to~$p$, and
thus the order~$\abs{T_{\FF_p}(\FF_p)}$ is prime to~$p$.

From the exact sequence
\[
0\to p\cdot\mathfrak{t}_{\ZZ_p}\xrightarrow{\exp} T(\ZZ_p)\to T_{\FF_p}(\FF_p)
\]
we deduce that~$U_p:=\exp(p\mathfrak{t}_{\ZZ_p})\leq T(\ZZ_p)$ is a topological~$p$-group 
and~$\frac{T(\ZZ_p)}{U_p}\hookrightarrow T(\FF_p)$ has order prime to~$p$.

For any open subgroup~$H\leq T(\ZZ_p)$, we have
\begin{equation}\label{torus orbit factorisation}
[T(\ZZ_p):H]=[T(\ZZ_p):U_p\cdot H]\cdot [U_pH:U_p\cap H].
\end{equation}
We now choose~$H$ defined by~$\phi(H)=K_p\cap \phi(T(\ZZ_p))$. We have
\begin{equation}\label{substitut T and U}
[T(\ZZ_p):H]=[T(\ZZ_p)]_p,\qquad [U_p:U_p\cap H]=[U_p]_p
\end{equation}
and
\begin{equation}\label{substitut T mod U}
[T(\ZZ_p):U_p\cdot H]=
\abs*{\frac{\phi(T(\ZZ_p))}{\phi(\exp(2p\mathfrak{t}_{\ZZ_p}))\cdot\phi(M(\ZZ_p))\cap K_p}}.
\end{equation}
Substituting~\eqref{substitut T and U} and~\eqref{substitut T mod U} in~\eqref{torus orbit factorisation} yields
\begin{equation}\label{torus orbit factorisation phi}
[T(\ZZ_p)]_p=[U_p]_p\cdot \abs*{\frac{\phi(T(\ZZ_p))}{\phi(\exp(2p\mathfrak{t}_{\ZZ_p}))\cdot\phi(M(\ZZ_p))\cap K_p}}.
\end{equation}
We now use~\eqref{torus orbit factorisation} and~\eqref{conclusion 1} and~\eqref{conclusion 3} from Theorem~\ref{cases local bounds} and conclude
\[
[T(\ZZ_p)]_p\geq \frac{1}{2Np}\cdot H_p(d\phi)\cdot \frac{p}{c}= \frac{1}{2cN}H_p(d\phi).\qedhere
\]
\end{proof}
We now prove Theorem~\ref{cases local bounds}.
\begin{proof}[Proof of conclusion~\ref{local bound 1}]
Assume for now the claim that~$\exp$ converges on~$2p\mathfrak{m}_{\ZZ_p}$ and~$U:=\exp(2p\mathfrak{m}_{\ZZ_p})\leq M(\ZZ_p)$. 
Let~$X_1,\ldots,X_k$ be generators of~$\mathfrak{m}_{\ZZ_p}$, then
\begin{equation}\label{B4 1 eq1}
H_p(d\phi)=\max\{H_p(d\phi X_1);\ldots;H_p(d\phi X_k)\}.
\end{equation}
As~$U_i:=\exp(2pX_i)^\ZZ\leq U$ for every~$i\in\{1;\ldots;k\}$ we have
\[
[\phi(U)]_p=\abs{\phi(U)\cdot K_p/K_p}\geq \abs{\phi(U_i)\cdot K_p/K_p}=[\phi(U_i)]_p,
\]
and thus
\begin{equation}\label{B4 1 eq2}
[\phi(U)]_p\geq \max_{i=1,\ldots,k}[\exp(2p\cdot d\phi(X_i))^\ZZ]_p.
\end{equation}
According to Th.~\ref{lemme_exp} for~$X=2p\cdot d\phi(X_i)$ we have
\begin{equation}\label{B4 1 eq3}
[\exp(2pd\phi(X_i))^\ZZ]_p\geq \abs{N^*}_p\cdot H_p(2p\cdot d\phi(X_i)).
\end{equation}
We remark
\begin{multline}\label{B4 1 eq4}
H_p(2p\cdot d\phi(X_i))=\max\{1;\norm{2p\cdot d\phi(X_i)}\}\\\geq\abs{2p}_p\cdot \max\{1;\norm{d\phi(X_i)}\}
=\abs{2p}_p\cdot \norm{d\phi(X_i)}. 
\end{multline}
Substituting~\eqref{B4 1 eq4} into~\eqref{B4 1 eq3} and~\eqref{B4 1 eq3} into~\eqref{B4 1 eq2}, we get
\[
[\phi(U_p)]_p\geq \abs{2pN^*}_p\cdot \max_{i=1,\ldots,k} H_p(d\phi(X_i))=\abs{2pN^*}_p\cdot H_p(d\phi).
\] 

We now recall why, for~$2pX\in 2p\mathfrak{m}_{\ZZ_p}$, the series~$\exp(2pX)$ converges and~$\exp(2pX)\in M(\ZZ_p)$ for~$2pX\in 2p\mathfrak{m}_{\ZZ_p}$.
\begin{proof}
We remark that~$\exp(2pT)\in\ZZ_{(p)}[[T]]$ and recall that the~$p$-adic radius of convergence  of~$\exp(2pT)$
is~$2\cdot p/ p^{\frac{1}{p-1}}>1$. For~$2pX\in 2p\mathfrak{m}_{\ZZ_p}$, we have~$\norm{X}\leq 1$ and so~$\exp(2pX)$ 
converges. We have~$\exp(2pX)\in M(d,\ZZ_p)$ because~$ \exp(2pT)\in\ZZ_p[[X]]$ has~$\ZZ_p$
entries.
Likewise~ and~$\exp(2pX)^{-1}=\exp(-2pX)\in M(d,\ZZ_p)$ and we conclude~$\exp(2pX)\in GL(N,\ZZ_p)$.
%Let~$P=0$ with~$P\in\QQ_p[T_{i,j},\det(T_{i,j})]$ be 
%one of the equations of~$M$ in~$GL(N)$. Then, noting~$\det\circ \exp (T_{i,j})= Tr(T_{i,j})$,
%\[
%P(\exp(T_{i,j}),Tr(T_{i,j}))(=0)
%\]
%certainly converges with same radius as~$\exp$, and we have
%\[P(\exp(T_{i,j},Tr(T_{i,j}))(X)=P(\exp(X),Tr(X)).\]
%Finally we notice that~$P(\exp(T_{i,j},Tr(T_{i,j}))$ ...
\end{proof}
Conclusion~\ref{local bound 1} has been proved.
\end{proof}
\begin{proof}[Proof of conclusion~\ref{local bound 2}]
Let~$X_1,\ldots,X_k$ be a nilpotent basis of~$\mathfrak{m}_{\ZZ_p}$. Then the~$d\phi(X_1),\ldots,d\phi(X_k)$ generate~$d\phi(\mathfrak{m}_{\ZZ_p})$
and there exists an~$i\in\{1;\ldots;k\}$ such that~$H_p(d\phi)=H_p(d\phi(X_i))$. Because~$X_i$ is nilpotent, we have
\[\exp(N^*\cdot X_i)=1+N^*\cdot X_i+\ldots+\frac{1}{(N-1)!}(N^*\cdot X_i)^{N-1}\]
and thus~$\exp(N^*\cdot X_i)\in M(\ZZ_p)$.

Thus
\[
[\phi(M(\ZZ_p))]_p\geq[\phi(\exp(N^*\cdot X_i))^{\ZZ}]_p.
\] 
Finally, by~\eqref{borne orbite general}, we have
\[
[\phi(\exp(N^*\cdot X_i\cdot \ZZ_p))]_p\geq H_p(d\phi(N^*\cdot X_i))/N.
\]
Because~$H_p(d\phi(N^*\cdot X_i))$ and~$[\phi(\exp(N^*\cdot X_i\cdot \ZZ_p))]_p$ are powers of~$p$, we actually have
\[
[\phi(\exp(N^*\cdot X_i\cdot \ZZ_p))]_p\geq\abs{N}_p\cdot H_p(d\phi(d^*\cdot X_i))\geq \abs{N^*}_p\cdot H_p(d\phi).\qedhere
\]

%We argue as above... attention exp entiere (ok si p grand, facile car nilpotent).
\end{proof}
Conclusion~\ref{local bound 4} is due to~\cite{EdYa} and we detail how their formulation \cite[Prop.~4.3.9]{EdYa} relates to ours.
\begin{proof}[Proof of conclusion~\ref{local bound 4}] We can discard finitely many primes and
assume~$p$ is big enough so that~\cite[Prop.~4.3.9]{EdYa} and its proof applies.

We first notice that, in the matrix algebra~$M(N,\QQ)$, the subalgebra~$\QQ[T(\QQ)]$ contains~$\mathfrak{t}$. 
\begin{proof} The inclusion of vector spaces can be checked after passing to~$\RR/\QQ$.
We know that
\[
\RR[T(\QQ)]=\RR[T(\RR)]
\]
because, by  weak approximation,~$T(\QQ)$ is dense in~$T(\RR)$. Let~$t$ be a sufficiently small element in~$\mathfrak{t}\tens\RR$,
so that~$\log(\exp(t))$ converges and~$\log(\exp(t))=t$. Then~$t\in\RR[\exp(t)]$, as is seen using Jordan forms, and~$\exp(t)\in T(\RR)$. Because~$\mathfrak{t}\tens\RR$ admits a basis of such elements, we can conclude.
%
%or~$\RR$ as follows. If~$t\in T(\QQ)$ is close enough to~$1$ in the topology of~$K$, then~$\log(t)$ is
%defined and belongs to~$\mathfrak{t}\tens K$. By weak approximation,~$T(\QQ)$ is dense in~$T(K)$ and
%the closure in~$\mathfrak{t}\tens K$ of the converging~$\log(t)$ with~$t\in T(\QQ)$ will contain the
%converging~$\log(t)$ with~$t\in T(K)$, and the latter contains a neighbourhood of~$0$ in~$ \mathfrak{t}\tens K$.
%Certainly converging~$\log(t)$ with~$t\in T(\QQ)$ will span~$\mathfrak{t}\tens K$. We also know that
%\[
%\log(t)\in K[t]
%\]
%by looking at a Jordan form (cf.~\cite{}). This establishes
%\[
%\mathfrak{t}\tens K \subset \QQ[T(\QQ)]\tens K.
%\]
%This is an inclusion of~$K$-vector spaces which are defined over~$\QQ$, and thus an inclusion
%of the underlying~$\QQ$-vector spaces.
\end{proof}

We can choose~$t_1,\ldots,t_k$ in~$\QQ[T(\QQ)]$ so that
\[
\mathfrak{t}\subset t_1\cdot \QQ+\ldots+t_k\cdot \QQ
\]
and thus~$t_1\cdot \ZZ+\ldots+t_k\cdot \ZZ$ contains a lattice of~$\mathfrak{t}$.
It will hence contain~$n\cdot(\mathfrak{t}\cap \mathfrak{gl}(N,\ZZ))$
for some commensurability index~$n\in\ZZ_{\geq1}$.

As we discard finitely many primes~$p$, we may assume that~$p$ do not divide the denominators of the~$t_i$ 
and 
do not divide~$n$. We will then have
\[
t_1,\ldots,t_k\in T(\ZZ_{(p)})
\]
and
\begin{equation}\label{taupe}
\mathfrak{t}_{\ZZ_p}=\mathfrak{t}\cap \mathfrak{gl}(N,\ZZ_{(p)})\subset t_1\cdot \ZZ_{(p)}+\ldots+t_k\cdot \ZZ_{(p)},
\end{equation}
and, applying~$\tens_{\ZZ_{(p)}}\ZZ_p$, we may replace~$\ZZ_{(p)}$ by~$\ZZ_p$.

Let~$\phi \in W(\QQ_p)$.  
Using Th.~\ref{rational conjugation representations}, we can write
\[
\phi= g\phi_0 g^{-1}
\]
for some~$g\in GL(N,\QQ_p)$. We assume~$H_p(d\phi)\neq 1$, that is~
\[
g\mathfrak{t}_{\ZZ_p}g^{-1}\not \subset \mathfrak{gl}(N,\ZZ_p),
\]
and, by~\eqref{taupe}, there is at least one~$i\in\{1;\ldots;k\}$ such 
that
\[
gt_ig^{-1}\not\in\mathfrak{gl}(N,\ZZ_p).
\] 
Equivalently~$gt_ig^{-1}\not\in GL(N,\ZZ_p)$, which also means
\[
t_i\cdot g{\ZZ_p}^d\neq g{\ZZ_p}^d.
\]
As~$t_i\in T(\ZZ_p)$, this implies, in the sense of~\cite[Prop.~4.3.9]{EdYa} 
(for $W_{\ZZ_p}=g{\ZZ_p}^d$),
\[
\text{``$T_{\ZZ_p}$ does not fix~$\{W_{\ZZ_p}\}$''}.
\]
% We can then use their conclusion which already 
%gives~\eqref{EY weak form}, which was seem to suffice for the purpose of this article.

Looking into the proof of~\cite[Prop.~4.3.9]{EdYa} we notice that their lower bound
is given by a lower bound of some orbit of~$T(\FF_p)$, thus, in~\eqref{EY weak form},
there exists~$n\in\ZZ_{\geq1}$ such that~$n$ divides~$\abs{T(\FF_p)}$ and
\[
[\phi(T(\ZZ_p))]_p\geq n\geq p/c.
\]
In the factorisation~\eqref{torus orbit factorisation phi} the first factor in the right-hand side
is a power of~$p$ and prime to~$n$. Thus the inequality~$[\phi(T(\ZZ_p))]_p\geq n$ comes from the second factor,
i.e. we have inequality of Conclusion~\ref{local bound 4}.
\end{proof}
\subsubsection{Subgroup principle}\label{Subgroup principle}
The following elementary lemmas were useful in passing to subgroups in the proofs of Theorems~\ref{general global bounds},~\ref{general local bounds} and~\ref{cases local bounds}. Proofs are left to the reader.
\begin{lem}[Global subgroup principle]
\begin{subequations}
	Let~$M_1,\ldots,M_k\leq M \leq GL(N)$ be algebraic groups over~$\QQ$
such that~$\mathfrak{m}_1+\ldots+\mathfrak{m}_k=\mathfrak{m}$. 

\begin{enumerate}
\item Then
\begin{equation}
\Lambda:=\mathfrak{m}_1\cap \mathfrak{gl}(N,\ZZ)+\ldots+\mathfrak{m}_k\cap \mathfrak{gl}(N,\ZZ)\leq\mathfrak{m}\cap \mathfrak{gl}(N,\ZZ)
\end{equation}
and the index
\begin{equation}
c=[\mathfrak{m}\cap \mathfrak{gl}(N,\ZZ):\Lambda]
\end{equation}
is finite. For every prime~$p$, we have
\begin{equation}
\Lambda\tens\ZZ_p=\mathfrak{m}_1\tens\QQ_p\cap \mathfrak{gl}(N,\ZZ_p)+\ldots+\mathfrak{m}_k\tens\QQ_p\cap \mathfrak{gl}(N,\ZZ_p)\leq\mathfrak{m}\tens\QQ_p\cap \mathfrak{gl}(N,\ZZ_p)
\end{equation}
and
\begin{equation}
[\mathfrak{m}\tens\QQ_p\cap \mathfrak{gl}(N,\ZZ_p):\Lambda\tens\ZZ_p]=\abs{1/c}_p
\end{equation}
with~$\abs{1/c}_p\leq c$ and~$\abs{1/c}_p=1$ if~$gcd(c,p)=1$.
\item Assume moreover that, for some morphism~$\phi:M\to GL(d)$ defined over~$\QQ$, we have
\begin{equation}\label{explicit local adelic bound for subgroup}
[\phi(M_i(\widehat{\ZZ})):\phi(M_i(\widehat{\ZZ}))\cap GL(d,\widehat{\ZZ})]\geq\frac{H_f(d\phi)}{c_i}.
\end{equation}
Then we have, with~$c=n\cdot \max\{c_1;\ldots;c_k\}$,
\begin{equation}%\label{explicit local adelic bound}
[\phi(M(\widehat{\ZZ})):\phi(M(\widehat{\ZZ}))\cap GL(d,\widehat{\ZZ})]\geq\frac{H_f(d\phi)}{c}.
\end{equation}
\end{enumerate}
\end{subequations}
\end{lem}
\begin{lem}[Local subgroup principle]
\begin{subequations}
Let~$p$ be a prime and~$M_1,\ldots,M_k\leq M \leq GL(N)$ be algebraic groups over~$\QQ_p$.
\begin{enumerate}
\item Then
\begin{equation}
[M(\ZZ_p)]_p\geq \max_{i\in\{1;\ldots;k\}}[M_i(\ZZ_p)]_p.
\end{equation}
\item Assume that~$\mathfrak{m}_1+\ldots+\mathfrak{m}_k=\mathfrak{m}$, then the index
\begin{equation}
[\mathfrak{m}_{\ZZ_p}:\Lambda]=n
\end{equation}
is a finite power of~$p$.
\item With~$n$ as above, for any~$\QQ_p$ linear map~$\Phi:\mathfrak{m}\to\mathfrak{gl}(d,\QQ_p)$, we have
\begin{equation}\label{90C}
\frac{1}{n} H_p(\Phi)\leq  \max_{i\in\{1;\ldots;k\}} H_p(\Phi|_{\mathfrak{m}_i})\leq H_p(\Phi).
\end{equation}
\item Assume moreover, for some morphism~$\phi:M\to GL(N)$ defined over~$\QQ_p$, that we have~\eqref{90C}
for~$\Phi=d\phi$ and that
\begin{equation}
\forall i\in\{1;\ldots;k\}, [M_i(\ZZ_p)]_p\geq \frac{1}{c_i}\cdot H_p(d\phi).
\end{equation}
Then we have, with~$c=n\cdot \max\{c_1;\ldots;c_k\}$,
\begin{equation}
[M(\ZZ_p)_p]\geq \frac{1}{c}\cdot H_p(d\phi).
\end{equation}
\end{enumerate}
\end{subequations}
\end{lem}

\section{Upper bound on Adelic orbits}\label{app c}
In this appendix, we prove upper bounds on adelic orbits. 
Combined with Prop.~\ref{adelic galois orbits prop} this implies corresponding
upper bounds on Galois orbits. This is not used in the proof of our
main result but we believe can be useful in other contexts. 

\begin{theorem}\label{C1} Let~$M\leq G$ be reductive groups over~$\QQ$,
and~$K\leq G(\AAA_f)$ be a compact open subgroup and~$K_M\leq K\cap M(\AAA_f)$ be a compact subgroup.

Let~$\phi_0:M\to G$ be the inclusion monomorphism, and~$W=G\cdot \phi_0$ be the conjugacy class
of~$\phi_0$, as an algebraic variety.

Let~$\iota:W\hookrightarrow \AAA^N$ be an affine embedding, and let~$H_f$ be as defined in~\eqref{notation:AF height}.
Then we have, as~$\phi$ describes~$W(\AAA_f)$, 
\begin{equation}\label{Extra bound 1}
[\phi(K_M):\phi(K_M)\cap K]\dom H_{\iota,f}(\phi).
\end{equation}
\end{theorem}

We prove a more precise version. Let~$\rho:G\hookrightarrow GL(d)$ be a faithful representation
and let us identify~$G$ with~$\rho(G)$. In the associative algebra~$\End(\QQ^d)$, we denote
the subalgebras linearly generated by~$M(\QQ)$ and~$G(\QQ)$ by
\[
B_M:=\sum_{m\in M(\QQ)} \QQ\cdot m\text{ and }B_G:=\sum_{g\in G(\QQ)} \QQ\cdot g.
\]
Let~$\Phi_0:B_M\to B_G$ denote the inclusion. We have~$M(\QQ)\subseteq B_M$, and~$G(\QQ)\subseteq B_G$,
and~$\phi_0:M(\QQ)\to G(\QQ)$ is the restriction of~$\Phi_0$.

For every field extension~$L/\QQ$, and~$\phi=g\cdot \phi_0\cdot g^{-1}\in W(L)$, with~$g\in G(\overline{L})$, the map
\[
B_\phi=g\cdot \Phi_0\cdot g^{-1}:B_M\tens L\to B_G\tens L
\]
is a $L$-linear extension of~$\phi$ to~$B_M\tens L$, and is the unique $L$-linear extension.

We choose linear bases of~$B_M$ and~$B_G$ generating~$B_M\cap \End(\ZZ^d)$ and~$B_G\cap \End(\ZZ^d)$ respectively, and we consider the corresponding isomorphism~$\Hom(B_M,B_G)\simeq \QQ^{\dim(B_M)\cdot \dim(B_G)}$.
Then~$\phi\mapsto B_\phi$ induces an affine embedding~$\iota_\rho:W\hookrightarrow \Hom(B_M,B_G)\simeq \QQ^{\dim(B_M)\cdot \dim(B_G)}$.
\begin{theorem}\label{C2} Define~$G(\wh{\ZZ}):=G(\AAA_f)\cap GL(d,\wh{\ZZ})$ and~$M(\wh{\ZZ}):=M(\AAA_f)\cap GL(d,\wh{\ZZ})$.
Then, for every~$\phi\in W(\AAA_f)$, we have
\[
[\phi(M(\widehat{\ZZ})):\phi(M(\widehat{\ZZ}))\cap G(\widehat{\ZZ})]\leq  H_{\iota_\rho,f}(\phi)^{2+d^2}.
\]
\end{theorem}
We note that if~$G$ is of adjoint type, we can use the adjoint representation and pick~$d=\dim(G)$.

Let us prove Th.~\ref{C2}.
\begin{proof}We endow~$\Hom(B_M\tens\QQ_p,B_G\tens\QQ_p)$ with the norm
\begin{multline}\label{C2eq}
\norm{\Phi}=\min \{ p^k\in p^{\ZZ}:
\\
\forall m\in B_M\tens\QQ_p\cap \End({\ZZ_p}^d), p^k\cdot\Phi(m)\in  B_G\tens\QQ_p\cap \End({\ZZ_p}^d)\}.
\end{multline}
We note that~$H_{\iota_\rho,p}(\phi)=\max\{1;\norm{B_\phi}\}$.

It suffices to prove that, for every prime~$p$, and~$\phi\in W(\QQ_p)$, we have
\begin{equation}\label{App C; proof 2; eq 1}
[\phi(M(\ZZ_p)):\phi(M(\ZZ_p))\cap G(\ZZ_p)]\leq  {\norm{B_\phi}_p}^{2+d^2}.
\end{equation}
Let us write~$\norm{B_\phi}_p=p^k$. Then, in the notations of Lemma~\ref{App C Lemma}, we have
\[
\phi(M(\ZZ_p))\subseteq S(d,p,p^k).
\]

Thus~\eqref{App C; proof 2; eq 1} follows from~\eqref{App C Lemma eq}.
\end{proof}
We deduce Th.~\ref{C1} from~Th.~\ref{C2}.
\begin{proof}

The assumptions imply the finiteness of
\[
C_M:=[K_M:K_M\cap M(\wh{\ZZ})]=[\phi(K_M):\phi(K_M)\cap \phi(M(\wh{\ZZ}))]
\]
and
\[
C_G:=[G(\wh{\ZZ}):K\cap G(\wh{\ZZ})].
\]
We have
\[
[\phi(K_M):\phi(K_M)\cap K]\leq C_M\cdot C_G\cdot [\phi(M(\widehat{\ZZ})):\phi(M(\widehat{\ZZ}))\cap G(\widehat{\ZZ})].
\]
By Prop.~\ref{prop:functoriality heights}, we have~$H_f\approx H_{\iota_\rho,f}$. Using~\eqref{C2eq}, we conclude 
\[
[\phi(K_M):\phi(K_M)\cap K]\leq C_M\cdot C_G \cdot H_{\iota_\rho,f}(\phi)^{2+d^2}\approx H_{f}(\phi).\qedhere
\]
\end{proof}

\begin{lem}\label{App C Lemma}
Let~$p$ be a prime,~$d$ be in~$\ZZ_{\geq0}$, and~$k$ be in~$\ZZ_{\geq0}$.

Define~$S(d,p,p^k)=\{b\in \End({\QQ_p}^d):\norm{b}\leq p^k, \det(b)\in {\ZZ_p}^\times \}$.

Then~$S(d,p,p^k)=S(d,p,p^k)\cdot GL(d,\ZZ_p)$ and
\begin{equation}\label{App C Lemma eq}
\#S(d,p,p^k)/GL(d,\ZZ_p)\leq (p^k)^{2+d^2}.
\end{equation}
\end{lem}
\begin{proof}
We endow~$\End({\QQ_p}^d)$ with the additive Haar measure~$\mu$ normalised by~$\mu(B(1))=1$,
where~$B(p^k)$, for~$k\in\ZZ_{\geq1}$ is the ball of radius~$p^k$. One knows that the Haar measure satisfies~$\mu(g\cdot A)=\abs{\det(g)}\cdot \mu(A)$.

For~$A=B(1)$ and~$g=p^k\cdot \mathrm{Id}$ this yields
\[
\mu(S(d,p,p^k))\leq \mu(B(p^k))=(p^k)^{d^2}.
\]
For~$b\in GL(N,\QQ_p)$ such that~$\det(b)\in {\ZZ_p}^\times$ this yields
\begin{equation}\label{App C Lemma proof eq}
\mu(b\cdot GL(d,\ZZ_p))=\mu(GL(d,\ZZ_p)).
\end{equation}
One can also check
\begin{multline}
\mu(GL(d,\ZZ_p))=\frac{\#GL(d,\FF_p)}{\#\End({\FF_p}^d)}=\prod_{i=1}^d1-\frac{1}{p^i}\\
\geq \prod_{i=1}^\infty 1-\frac{1}{2^i}\geq 0.25\geq 1/p^2.
\end{multline}

The norm multiplicativity~$\norm{b\cdot g}=\norm{b}\cdot\norm{g}$ implies the right invariance
\begin{equation}\label{App C eq inv}
S(d,p,p^k)=S(d,p,p^k)\cdot GL(d,\ZZ_p).
\end{equation}
Equivalently, we can write~$S(d,p,p^k)=b_1\cdot GL(d,\ZZ_p)\sqcup\ldots\sqcup b_c\cdot GL(d,\ZZ_p)$,
with~$c=\#S(d,p,p^k)/GL(d,\ZZ_p)$.

Using~\eqref{App C Lemma proof eq}, we deduce
\[
\#S(d,p,p^k)/GL(d,\ZZ_p)=\mu(S(d,p,p^k))/\mu(GL(d,\ZZ_p)).
\]

Assume~$k=0$. Then~\eqref{App C Lemma eq} follows from
\[
S(d,p,p^k)=GL(d,\ZZ_p)\text{ and }\#S(d,p,p^k)/GL(d,\ZZ_p)=1\leq 1^{2+d^2}.
\]
We may now assume~$k\geq 1$. Then~\eqref{App C Lemma eq} follows from
\[
\#S(d,p,p^k)/GL(d,\ZZ_p)\leq p^2\cdot (p^k)^{(d^2)}\leq (p^k)^{2+d^2}.\qedhere
\]
\end{proof}

\end{document}